\documentclass[12pt]{amsart}

\usepackage{upref}
\usepackage{enumitem,dsfont,amssymb,bm}
\usepackage{amsmath}
\usepackage{amsthm}
\usepackage{graphicx,color}
\usepackage{mathtools}
\usepackage{mathrsfs}
\usepackage{braket}
\usepackage[english]{babel}
\usepackage[colorlinks=true, pdftitle={TITLE}, pdfauthor={Timothy 
	Allen Collier and Daniel Hauer}]{hyperref}

\numberwithin{equation}{section}

\newtheorem{theorem}{Theorem}[section]
\newtheorem{proposition}[theorem]{Proposition}
\newtheorem{lemma}[theorem]{Lemma}
\newtheorem{corollary}[theorem]{Corollary}

\newtheorem{definition}[theorem]{Definition}

\newtheorem{remark}[theorem]{Remark}

\newcommand\R{{\mathbb{R}}}
\newcommand\N{\mathbb{N}}

\newcommand\E{\mathcal{E}}

\newcommand\td{\mathrm{d}}
\newcommand\dx{\mathrm{d}x }
\newcommand\dtau{\mathrm{d}\tau }
\newcommand\dy{\mathrm{d}y }
\newcommand\dr{\mathrm{d}r }
\newcommand\dmu{\mathrm{d}\mu }
\newcommand\dz{\mathrm{d}z }

\newcommand\dt{\mathrm{d}t }

\def\1{\raisebox{2pt}{\rm{$\chi$}}}
\DeclareMathOperator*{\sign}{sign}

\newcommand\abs[1]{\lvert#1\rvert}
\newcommand\norm[1]{\lVert#1\rVert}
\newcommand\lnorm[1]{\left\lVert#1\right\rVert}

\definecolor{darkred}{rgb}{0.7,0.1,0.1}

\begin{document}

%\begin{frontmatter}
    \title[Barenblatt solution of a doubly nonlinear nonlocal eq.]{The Barenblatt solution of an evolution problem governed by a doubly nonlinear nonlocal operator}
        
    \author{Timothy A. Collier}
    \author{Daniel Hauer}
    \address{School of Mathematics and Statistics, The University of Sydney, NSW 2006, Australia}
    \email{\href{mailto:timothyc@maths.usyd.edu.au}
    {timothyc@maths.usyd.edu.au}}
    \email{\href{mailto:daniel.hauer@sydney.edu.au}{daniel.hauer@sydney.edu.au}}

    \begin{abstract}
      In this article, we prove existence and uniqueness of the
      Barenblatt solution of the evolution equation on the whole
      Euclidean space where the principle part is the nonlocal
      fractional $p$-Laplacian composed with a power function. Our
      proof generalizes methods developped by V\'azquez [Nonlinear Anal., 199
      (2022), Calc. Var. Partial Differential Equations, 60 (2021)]
      for the evolution equation driven by the fractional
      $p$-Laplacian on the whole Euclidean space. In particular, we
      required an Aleksandrov symmetry principle, which can be applied
      to the mild solutions of the evolution equation in $L^1$
      governed by the doubly nonlinear nonlocal operator, and the
      construction of global barrier functions. The Aleksandrov
      symmetry principle might be of independent interest. % As an
      % application, we determine the long-time asymptotic behavior of
      % every mild solution in $L^1$.
    \end{abstract}
    


\maketitle
\tableofcontents

\section{Introduction}
	
	% PROBLEM
Let $1<p<\infty$, $0<s<1$, $m>0$, $d\in \N$, and $L^1$ denotes the
classical Lebesgue space $L^1(\R^{d})$. Then, the aim of this
paper is to determine the self-similarity properties and the
long-time asymptotic behaviour of the mild solution $u$ in $L^1$ (see
Definition~\ref{def:mild-solutions}) of the doubly nonlinear nonlocal diffusion equation
\begin{equation}
  \label{eq:heat-eq}
  \partial_{t}u(t)+(-\Delta_{p})^{s}u^m(t)=0\qquad 
  \text{in $\R^d \times(0,\infty)$}
\end{equation}
for every given initial data $u(0)=u_{0}\in L^1$. In the diffusion
equation~\eqref{eq:heat-eq}, the principal part $(-\Delta_{p})^{s}u^m$
models a nonlocal diffusion and is the composi\-tion operator of the
\emph{(variational) fractional $p$-Laplacian} on $\R^{d}$ (with
vanishing conditions at infinity) defined by
\begin{equation}
  \label{eq:11}
  \langle
  (-\Delta_{p})^{s}u,\xi\rangle:= \frac{1}{2}\int_{\R^{2d}}\tfrac{|u(x)-u(y)|^{p-2}(u(x)-u(y))
    (\xi(x)-\xi(y))}{|x-y|^{d+ps}}\;
  \dx\dy
\end{equation}
for $\xi \in W^{s,p}(\R^d)$ and the power function
$u^m:=\abs{u}^{m-1}u$.\medskip

In particular, we focus on 
the existence and uniqueness of the so-called \emph{Barenblatt
  solution} or also \emph{fundamental solution} $\Gamma$ to
equation~\eqref{eq:heat-eq}, which is a self-similar (see
Section~\ref{sec:selfsimilar}), strong distributional solution (see
Definition~\ref{def:distributional-L1}) of equation~\eqref{eq:heat-eq}
with a Dirac delta (of mass $M>0$) as initial value in the sense of
distribution; that is,
\begin{equation}
  \label{eq:initial-condition-Gamma}
  \lim_{t\to
    0+}\int_{\R^{d}}\Gamma(x,t;M)\,\xi(x)\,\dx=M\,\xi(0)=
  \langle M\delta_{0},\xi\rangle_{\mathcal{D}',\mathcal{D}}
\end{equation}
for every test-function $\xi\in C^{\infty}_{c}(\R^{d})$. More
precisely, we show that the following result holds.

\begin{theorem}
  \label{thm:BarenblattSolutions}
  For $d\ge 1$, $0<s<1$, $m\ge 1$, let 
 \begin{equation}
   \label{eq:pcDef}
   p_{m,c}:= \frac{d(1+m)}{md+s},
 \end{equation}
 and $p_{c}<p<\infty$ such that $m(p-1) \neq 1$. Then for every
 $M > 0$, there exists a unique strong distributional solution
 $\Gamma$ in $L^1$ of the doubly nonlinear nonlocal diffusion
 equation~\eqref{eq:heat-eq} having $M\delta_0$ as initial datum in
 the sense~\eqref{eq:initial-condition-Gamma}. The function
 $\Gamma : \R^{d}\times (0,\infty)\to \R$ is given by
  \begin{displaymath}
    %\label{eq:selfSimilarForm}
    \Gamma(x,t;M) = M^{sp\beta}t^{-d\beta}F(M^{-(m(p-1)-1)\beta}xt^{-\beta})
 \end{displaymath}
for every $(x,t)\in \R^{d}\times (0,\infty)$ with self-similar exponents
 \begin{equation}
    \label{eq:scalingFactors}
    \alpha=d\beta\quad\text{ and }\quad \beta = \frac{1}{d(m(p-1)-1)+sp},
 \end{equation}
 and the profile function $F: \R^d \rightarrow (0,\infty)$ is continuous,
$F(r)$ with $r=\abs{x} t^{-\beta}$ is  radially symmetric, decreasing
and has decay at infinity
\begin{displaymath}
  F(r)\sim \begin{cases}
        r^{-d-sp}& \text{if $p > p_1$,}\\
        r^{-d-sp}\log(r)& \text{if $p = p_1$,}\\
        r^{-\frac{sp}{1-m(p-1)}}& \text{if $p_{m,c}< p < p_1$,}
        \end{cases}
\end{displaymath}
where $p_1$ is the positive solution of the quadratic equation (in $q$)
\begin{displaymath}
  (s q+d) m (q-1)-d=0.
\end{displaymath}
Moreover, $\Gamma(x,t;M)$ decays uniformly with respect to $x$ as
$t\to\infty$.
\end{theorem}

It is important to note that Theorem~\ref{thm:BarenblattSolutions}
extends the ones without power ($m=1$) and for $2d/(d+s)<p<\infty$
obtained in the excellent works~\cite{MR4114983,MR4280519} by
V\'azquez. Note the critical parameter $p_c$ given in~\eqref{eq:pcDef}
reduces to the critical parameter $p_{1,c}=2d/(d+s)$ found by V\'azquez. The
requirement $p_{m,c}<p<\infty$ ensures that the self-similar
transformations presented in Section~\ref{sec:selfsimilar} and the
global barrier functions as stated in Section~\ref{sec:barrier}
exist. The method of our proofs follows the one elaborated in
\cite{MR4114983,MR4280519} by V\'azquez. But we also refine here some
of his arguments. We give the details of the proof of
Theorem~\ref{thm:BarenblattSolutions} in
Section~\ref{section:existence-of-Barenblatt} (existence of $\Gamma$)
and Section~\ref{sec:barenblattUniqueness} (uniquiness of $\Gamma$).\medskip

In particular, for the proof of Theorem~\ref{thm:BarenblattSolutions},
we require the following Aleksandrov symmetry
principle~\cite{MR0147776} for mild solutions $u$ in $L^1$
of~\eqref{eq:heat-eq}.

\begin{theorem}[Aleksandrov's symmetry principle]
  \label{cor:aleksandrov}
  Let $0 <s<1$, $1<p<\infty$, $d \ge 1$, and $m > s$. Then every mild
  solutions $u$ in $L^1$ of the diffusion equation~\eqref{eq:heat-eq}with
  non-negative, compactly supported initial data $u_{0}$ in a ball
  $B_R(0)$ are radially decreasing in space for all $|x| \ge R$ and
  $t \ge 0$.
\end{theorem}

Our proof of Theorem~\ref{cor:aleksandrov} relies on the radial
symmetry of the fractional $p$-Laplacian $(-\Delta_{p})^{s}$ and on
sharpening an estimate from~\cite{MR4114983} in the case of the
fractional $p$-Laplacian.  For this, we prove the following comparison
principle between a mild solution $u$ of \eqref{eq:heat-eq} and its
reflection $\Pi$ on a half-space. More precisely, we deal here with the
following type of projections.

\begin{definition}
  For given parameters $a_{1}, \dots, a_{d}, b\in \R$, the set
  \begin{displaymath}
    H:=\{x\,\vert\,a_{1}x_{1}+\dots +a_{d}x_{d}=b\} 
  \end{displaymath}
  is called an \emph{affine hyperplane} in $\R^{d}$ splitting $\R^{d}$
  into the two \emph{half spaces}
  \begin{align*}
    H_{+}&:=\{x\,\vert\,a_{1}x_{1}+\dots +a_{d}x_{d}>b\},\\
    H_{-}&:=\{x\,\vert\,a_{1}x_{1}+\dots +a_{d}x_{d}<b\}
  \end{align*}
  if there exits $i\in \{1,\dots, d\}$ such that $a_{i}\neq
  1$. Further, let $\Pi_{0}:\R^d\rightarrow \R^d$ be the reflection
  around the hyperplane $\{x\in \R^{d}\,\vert\,x_{d}=0\}$ given by
  \begin{displaymath}
    \Pi_{0}x:=(x',-x_{d})\in \R^{d}\qquad\text{for every
      $x=(x',x_{d})\in \R^{d-1}\times \R$.}
  \end{displaymath}
  Then, for a given affine hyperplane $H$, we call $\Pi$ the
  \emph{reflection around $H$ splitting $\R^d$ into two half spaces
    $H_{+}$ and $H_{-}$} if after rotation and translation $\Pi$
  coincides with $\Pi_{0}$.
  % ; that is there is an rotation matrix $O\in
  % SO(d)$ and a vector $\vec{b}\in \R^{d}$ such that $(O+\vec{b})\Pi=\Pi_{0}$
\end{definition}

The method applied here modifies
a standard comparison principle idea, used for example
in~\cite{MR2286292}, obtaining estimates on the whole domain.  This
differs to the method used for the fractional Laplacian
in~\cite{MR3191976} so that we may use the variational formulation of
the fractional $p$-Laplacian (see
also~\cite{MR1260981}and~\cite{MR4114983}).

\begin{theorem}[Comparison on reflections]
  \label{thm:aleksandrovReflection}
  Let $1<p<\infty$, $0 <s <1$, $d \ge 1$, and $m > s$. Further, for a
  given affine hyperplane $H$ in $\R^{d}$, let
  $\Pi:\R^d\rightarrow \R^d$ be the reflection around $H$ splitting
  $\R^d$ into two half spaces $H_{+}$ and $H_{-}$. Then, for every
  mild solution $u$ of the diffusion equation~\eqref{eq:heat-eq} with initial data
  $u_0 \in L^1$, if
  \begin{displaymath}
    u_0(x) \le u_0(\Pi(x))\quad \text{for a.e. }x\in H_{+},
  \end{displaymath}
  then one has that
  \begin{displaymath}
    u(x,t) \le u(\Pi(x),t)\quad \text{for a.e. $x\in H_{+}$ and every $t>0$.}
  \end{displaymath}
\end{theorem}

Nonlinear integro-differential operators such as the doubly
nonlinear non\-local operator $(-\Delta_{p})^{s}\phi$ for $\phi$ being
a monotone function have recently received much interest for their role in the mathematical
analysis of anomalous diffusion, as well as their applications
in such fields as statistical mechanics, physics, finance,
fluid dynamics and image processing. We refer the interested
reader, for example,
to~\cite{MR1406564,MR2042661,MR2480109,MR3821542,MR2763032}.\medskip

Well-posedness in $L^1$ and regularity properties of mild solutions
$u$ of the parabolic boundary-value problem
\begin{displaymath}
  \begin{cases}
    \begin{alignedat}{2}
	\partial_{t}u(t)+(-\Delta_{p})^{s}u^m(t)&=0 &&\quad 
	\text{in $\Omega\times(0,\infty)$,}\\
	 u(t)&=0 \quad&&\text{on $\partial\Omega\times (0,\infty)$,}
    \end{alignedat}
  \end{cases}
\end{displaymath}
where $\Omega$ is an an open subset of $\R^{d}$
have recently been studied by Giacomoni et al.~\cite{MR4230963} in the range
$\frac{1}{2p-1} \le m < 1$ with $1<p < \frac{d}{s}$ and for bounded $\Omega$, 
and by the authors in~\cite{doublyNonlinearPaper} for all $m > 0$,
$1<p<\infty $, and on general open subsets of $\R^d$.

The general theory of accretive operators in $L^1$ (see, for example,
\cite{MR2582280} and, in particular,
\cite{CoulHau2016,BenilanNonlinearEvolutionEqns}) provides existence
in the sense of mild solutions which enjoy further
regularity/smoothing properties due to the homogeneity properties
(see~\cite{BenilanRegularizingMR648452}) and sub-differential
structure of the fractional $p$-Laplacian (see~\cite{BenilanGariepy,
  doublyNonlinearPaper}). Hence we will retain this $L^1$ setting,
expanding upon the regularity results of~\cite{doublyNonlinearPaper}
in Section~\ref{sec:notion-of-solution} to prove the existence of
strong distributional solutions to the diffusion
equation~\eqref{eq:heat-eq} on the whole Euclidean space
$\R^{d}$.\medskip
 
% Our focus in this article is to establish an existence result for
% \emph{Barenblatt solutions} to equation~\eqref{eq:heat-eq}.  In
% particular, these are self-similar solutions with a Dirac delta as
% initial data in a limiting sense as $t \rightarrow 0$.
We introduce
self-similar solutions to~\eqref{eq:heat-eq} in
Section~\ref{sec:selfsimilar}, in particular the appropriate scaling
transformation under which these solutions are invariant.
Furthermore, we prove an Aleksandrov symmetry principle in
Section~\ref{sec:Aleksandrov} and a global barrier for solutions
to~\eqref{eq:heat-eq} in Section~\ref{sec:barrier} which we use in the
proof of our main result.

% Barenblatt solutions
% are often otherwise referred to as \emph{fundamental} or
% \emph{source-type} solutions due to their role in the theory of linear
% equations.

% In the case of $m = 1$ and $p > \frac{2d}{d+s}$, Barenblatt
% solutions for~\eqref{eq:heat-eq} have been studied by V\'azquez
% in~\cite{MR4114983,MR4280519}.

\section{Self-similar solutions}
    \label{sec:selfsimilar}
    We are interested in scaling transformations of time and space under which~\eqref{eq:heat-eq} is invariant. 
    Moreover, we consider self-similar solutions, that is, classes of solutions which are invariant under such a one-parameter group of scaling symmetries.
    Such solutions plays a key role in the asymptotic behaviour as well as the existence and regularity of such evolution equations. 
     
	We first introduce two important scaling transformations for~\eqref{eq:heat-eq}. Scaling time, space and amplitude, we have that~\eqref{eq:heat-eq} is invariant under the transformation
	\begin{equation}
	\label{eq:transform1}
	    T_k u(x,t) = k^d u(kx,k^{dm(p-1)-d+sp}t)
	\end{equation}
	for all $k > 0$. By the amplitude scaling, this transformation also conserves mass. We can scale only the time and amplitude with
	\begin{equation}
	\label{eq:transform2}
	    \hat{T}_M u(x,t) = M u(x,M^{m(p-1)-1}t)
	\end{equation}
	which will scale the mass of $u$ by the factor $M > 0$.
	
	Since~\eqref{eq:heat-eq} is invariant under the scaling transformation~\eqref{eq:transform1}, we look for self-similar solutions of the form
	\begin{equation}
	\label{eq:selfSimilarForm}
	U(x,t;M) = t^{-d\beta}F(xt^{-\beta};M)
	\end{equation}
	for some $\beta \in \R$. Here $F(z;M)$ is chosen to have mass $M > 0$ and the factor $t^{-d\beta}$ conserves the mass of $U$ in time. Substituting~\eqref{eq:selfSimilarForm} into~\eqref{eq:heat-eq}, we have that setting
	\begin{equation}
	\tag{\ref{eq:scalingFactors}}
	    \beta = \frac{1}{dm(p-1)-d+sp}
	\end{equation}
	we obtain such a self-similar solution where the profile $F$ satisfies
	\begin{align*}
	    -d\beta F(z) -\beta\nabla F(z)\cdot z+(-\Delta_p)^s (F(z))^m = 0.
	\end{align*}
    Equivalently,
	\begin{equation}
	\label{eq:selfSimilarEq}
	    (-\Delta_p)^s (F(z))^m = \beta \nabla\cdot(zF)
	\end{equation}
    which we refer to as the self-similar profile equation. 
    We will also use the radial form with $r = |z|$, given by
    \begin{equation}
    \label{eq:selfSimilarProfileRadial}
        (-\Delta_p)^s (F(z))^m = \beta r^{1-d}\left(r^d F(z)\right)_r.
    \end{equation}

	Note that we can scale the profile function with $M = 1$, $F(z;1)$, to have unit mass. Then by the rescaling~\eqref{eq:transform2}, we can deduce the profile of mass $M$ to be
	\begin{displaymath}
    	F(z;M) = M^{sp\beta}F(M^{-(m(p-1)-1)\beta}z;1).
	\end{displaymath}
 
\section{Notion of solution and well-posedness}
	\label{sec:notion-of-solution}
	In this section, we introduce different notions of solutions related to the initial-value problem
	\begin{equation}
	\label{eq:1}
	\begin{cases}
    \begin{alignedat}{2}
	\partial_{t}u(t)+(-\Delta_{p})^{s}u^m(t)&=0 &&\quad 
	\text{in $\R^d \times(0,\infty)$,}\\
	u(0)&=u_{0} &&\quad \text{on $\R^d$}
    \end{alignedat}
	\end{cases}
	\end{equation}
	with different regularity properties. The main result of this section provides well-posedness of problem~\eqref{eq:1} for given $u_{0}\in L^{1}(\R^d)$.\medskip 
	
	It is worth mentioning that for $1\le q\le \infty$, we will write $L^{q}$ and $L^{q}_{loc}$ to denote the classic Lebesgue spaces $L^{q}(\R^{d})$ and $L^{q}_{loc}(\R^{d})$, $\norm{\cdot}_{q}$ denotes the standard $L^{q}$-norm on $L^{q}$ with respect to the Lebesgue measure, and we write $L^{1\cap\infty}$ for the intersection $L^{1}\cap L^{\infty}$. Further, for $1<q<\infty$ and $0<s<1$, we write $W^{s,p}$ for the fractional Sobolev space $W^{s,p}(\R^d)$ equipped with the norm
	\begin{math}
	\norm{u}_{W^{s,p}} := \sqrt[p]{\norm{u}_p^p+[u]_{s,p}^p}.
	\end{math}
	Finally, we denote by $q^{\mbox{}_{\prime}}$ the H\"older-conjugate of $1\le q\le \infty$ given by $1/q +1/q^{\mbox{}_{\prime}}=1$.

    For an open domain $\Omega$, we define $W_0^{s,p}(\Omega)$ as the closure in $W^{s,p}$ of the set $C_c^{\infty}(\Omega)$ of test functions. By~\cite[Theorem~10.1.1]{MR1411441}, for $1<p<\infty$, this space can be characterized as
	\begin{displaymath}
	W^{s,p}_{0}(\Omega)=\left\lbrace u\in
	W^{s,p}\,\bigg|
	\begin{array}[c]{l}
	\exists \text{ } \overline{u} : \R^{d}\to \R\text{
		s.t. }\overline{u}=u\text{ a.e.~on $\R^{d}$}\\
	\text{ and } \overline{u}=0\text{ quasi-everywhere on }\R^{d}\setminus\Omega
	\end{array}
	\right\rbrace,
	\end{displaymath}
	where $\overline{u}$ denotes a (quasi-continuous) representative of $u$. While we focus on~\eqref{eq:1} on $\R^d$ in this paper, we will use a formulation in $\Omega \subset \R^d$ for approximation methods.
 
\subsection{\texorpdfstring{The fractional $p$-Laplacian}{}}
\label{sub:fractionalpLap}

In order to establish existence of solutions to the Cauchy problem~\eqref{eq:1}, we realize the doubly nonlinear operator $(-\Delta_{p})^{s}\cdot^{m}$ as an operator $A$ in $L^{1}$. To do this, we first realize the fractional $p$-Laplacian $(-\Delta_{p})^{s}$ on $\R^{d}$ as an operator in $L^2$ by using the concept of the sub-differential operator $\partial_{\! L^2}\E$ in $L^2$ defined (through its graph) by the set
\begin{displaymath}
    \partial_{\! L^2}\E=\Bigg\{
    (u,v)\in D(\E)\times L^2\,\Bigg\vert\,\int_{\R^{d}}v\,\xi\,\dx\le \liminf_{h\downarrow 0}
    \tfrac{\E(u+h\xi)-\E(u)}{h}\,\forall\,\xi\in D(\E)
    \Bigg\}
\end{displaymath}
of the energy functional $\E:L^2(\Omega) \rightarrow (-\infty,+\infty]$ defined by
\begin{equation}
    \label{eq:energyFunctional}
    \E(u)=
    \begin{cases}
        \frac{1}{2p}[u]_{s,p}^{p} & \text{if $u\in W^{s,p}\cap L^2$,}\\
        +\infty & \text{otherwise,}
    \end{cases}
\end{equation}
where $[\cdot]_{s,p}$ denotes the Gagliardo semi-norm 
	\begin{equation}
	\label{eq:seminormWsp}
	[u]_{s,p}:=\left(\int_{\R^d}\int_{\R^d}\frac{\abs{u(x)-u(y)}^{p}}{\abs{x-y}^{d+sp}}\dy\dx\right)^{1/p}.
\end{equation}
Since the energy functional $\E$ is G\^ateaux-differentiable in every direction $\xi\in D(\E)$, the sub-differential $\partial_{\! L^2}\E$ is, in fact, a well-defined mapping $\partial_{\! L^2}\E : D(\partial\E)\to L^{2}$ given by 
\begin{equation}
    \label{def:sub-diff-fract-plalacian}
    \partial_{\! L^2}\E(u):=v\qquad\text{for every $u\in D(\partial_{\! L^2}\E)$}
\end{equation}
where $v$ is uniquely defined via the integral equation
\begin{equation}
  \label{eq:variationalEq}
   \int_{\R^{d}} v\xi\dx= \frac{1}{2}\int_{\R^{2d}}\frac{|u(x)-u(y)|^{p-2}(u(x)-u(y))
		(\xi(x)-\xi(y))}{|x-y|^{d+ps}}\,\td(x,y)
\end{equation}
for every $\xi \in W^{s,p}\cap L^2$ (for more details see, for example,~\cite{doublyNonlinearPaper,CoulHau2016}). 

\begin{definition}
  \label{def:fractional-p-laplacian-in-L2}
  For $1<p<\infty$ and $0<s<1$, let $\partial_{\! L^2}\E$ be the sub-differential operator given by \eqref{def:sub-diff-fract-plalacian} and \eqref{eq:variationalEq}. Then, call the operator 
  \begin{displaymath}
      (-\Delta_{p})^{s}u:=\partial_{\! L^2}\E(u)\qquad
      \text{for every $u\in D(\partial_{\! L^2}\E)$}
  \end{displaymath}
  the \emph{fractional $p$-Laplacian} in $L^2$.
\end{definition}

It is interesting to note that it was recently shown in~\cite{MR4030247} (see also \cite{MR4303657}) that for $u$ sufficiently smooth, the fractional $p$-Laplacian $(-\Delta_{p})^{s}$ in $L^2$ can be rewritten as the Cauchy principal value of the singular integral
\begin{displaymath}
\label{eq:pLapSingularIntegral}
   (-\Delta_{p})^{s}u(x) = P.V.\int_{\R^d}\frac{|u(x)-u(y)|^{p-2}(u(x)-u(y))}{|x-y|^{d+sp}}\dy
\end{displaymath}
for every $x\in \R^{d}$.\medskip

%Note that if $u$ is defined on $\Omega \subset \R^d$, we extend $u$ by zero to be defined %on $\R^d$ in order to evaluate~\eqref{eq:seminormWsp}.

\subsection{The doubly nonlinear nonlocal operator}

Now, we are in the position to define a realization of the nonlocal doubly nonlinear operator $(-\Delta_{p})^{s}\cdot^{m}$ in $L^1$. 

\begin{definition}
  \label{def:realization-in-L1-of-NDN}
  Let $(-\Delta_{p})^{s}$ be the fractional $p$-Laplacian in $L^2$ and $(-\Delta_p)^s_{1\cap\infty}$ be the \emph{part of $(-\Delta_{p})^{s}$ in $L^{1\cap\infty}\times L^{1\cap\infty}$} given by
  \begin{equation}
     \label{A-operator-L1Linfty}
    (-\Delta_{p})^{s}_{1\cap\infty}:=\Big\{(u,v)\in L^{1\cap\infty}\times L^{1\cap\infty}\,\Big\vert\, v=(-\Delta_{p})^{s}u
  %  \begin{array}[c]{l}
  %  u\in W^{s,p} \text{ and } (u,v) \text{ satisfies \eqref{eq:variationalEq}}\\
  %  \text{for all } \xi \in L^ \cap W^{s,p}
  %  \end{array}
    \Big\}.
\end{equation}
 Further, for $m>0$, let $\varphi(r):=\abs{r}^{m-1}r$ for every $r\in \R$. Then, we call
 the operator $A$ given by
 \begin{displaymath}
    A:=\Bigg\{(u,v)\in L^{1}\times L^{1}\,\Bigg\vert\, 
   \begin{array}[c]{l}
    \exists\, ((u_{k},v_{k}))_{k\ge 1}
    \text{ s.t. } v_{k}= (-\Delta_{p})^{s}_{1\cap\infty}\varphi(u_{k}),\\
    \displaystyle \lim_{k\to\infty}u_{k}=u \text{ in } L^1 \text{ and }
    \lim_{k\to\infty}v_{k}=v \text{ in } L^1
   \end{array}
    \Bigg\}    
 \end{displaymath}
 the \emph{nonlocal doubly nonlinear operator} in $L^1$, since it provides a realization of the operator $(-\Delta_{p})^{s}\cdot^{m}$ in $L^1$.
\end{definition}

After introducing $A$, the realization in $L^{1}$ of the composed operator $(-\Delta_{p})^{s}\cdot^{m}$, we now aim to establish well-posedness of the following Cauchy problem (in $L^1$)
 \begin{equation}
	\label{eq:abstract-CP}
	\begin{cases}
    \begin{alignedat}{2}
	\tfrac{\td u}{\dt}(t)+Au(t)&\ni 0 &&\quad 
	\text{in $(0,\infty)$,}\\
	u(0)&=u_{0} &&\mbox{}
    \end{alignedat}
	\end{cases}
\end{equation}
with respect to the notion of \emph{mild solutions} and study their regularizing properties.

\subsection{\texorpdfstring{Well-posedness in $L^1$ \& Regularity}{}}

We begin by introducing the notion of solutions required in this paper.

\begin{definition}[Mild solutions]\label{def:mild-solutions}
 Let $A$ be an operator on $L^1$, where $L^1$ is equipped with a $\sigma$-finite measure space $(\Sigma,\mu)$. For given $u_{0}\in L^{1}$, a function $u \in C([0,\infty);L^{1})$ is called a \emph{mild} solution in $L^1$ of Cauchy problem~\eqref{eq:abstract-CP} if $u(0)=u_{0}$ in $L^{1}$, and for every $T>0$ and $\varepsilon>0$, there is a partition $\pi_{\varepsilon} : 0=t_0<t_1<\cdots <t_{N_{\varepsilon}}=T$ satisfying $\max_{i=1,\dots, N_{\varepsilon}}t_{i}-t_{i-1}<\varepsilon$, and there is a finite recursively defined sequence $(u_{i,\varepsilon})_{i=0}^{N_{\varepsilon}}\subseteq D(A)$ given by $u_{0,\varepsilon}=u_0$, and for every $i=1,\dots, N_{\varepsilon}-1$, $u_{i,\varepsilon}$ is defined by
 \begin{displaymath}
	    u_{i,\varepsilon}+(t_{i,\varepsilon}-t_{i-1,\varepsilon})Au_{i,\varepsilon}\ni u_{i-1,\varepsilon}
\end{displaymath}
such that the sequence $(u_{\varepsilon})_{\varepsilon>0}$ of step functions
	\begin{displaymath}
	    u_{\varepsilon}(t):=\sum_{i=0}^{N_{n}-1}u_{i,\varepsilon}\,\mathds{1}_{[t_{i,\varepsilon},t_{i+1,\varepsilon})}
	\end{displaymath}
	converges to $u$ in $C([0,T];L^{1})$ as $\varepsilon\to 0+$.
\end{definition}

Since mild solutions are just limits of step functions, we also introduce stronger notions of solutions (cf.~\cite{MR2582280}).

\begin{definition}[Strong solutions]\label{def:strong-solutions}
Let $A$ be an operator on $L^1$, where $L^1$ is equipped with a $\sigma$-finite measure space $(\Sigma,\mu)$. A function $u \in C((0,\infty);L^{1})$ is called a \emph{strong} solution in $L^{1}$ of
    \begin{equation}
      \label{eq:ODE}
      \tfrac{\td u}{\dt}(t)+Au(t)\ni 0 \qquad 
	\text{in $(0,\infty)$,}
    \end{equation}
    if for a.e.~$t>0$, $u(t)$ is differentiable in $L^{1}$, $u(t) \in D(A)$ and satisfies $-\tfrac{\td u}{\dt}(t)\in Au(t)$. 
    Further, a function $u \in C([0,\infty);L^{1})$ is called a \emph{strong} solution in $L^{1}$ of the Cauchy problem~\eqref{eq:abstract-CP} if $u(0)=u_{0}$ in $L^{1}$ and $u$ is a strong solution in $L^{1}$ of~\eqref{eq:ODE}.
\end{definition}

Since the nonlocal doubly nonlinear operator $A$ in $L^1$ is the closure
in $L^1\times L^1$, it is \emph{a priori} not clear whether strong
solutions $u$ of \eqref{eq:abstract-CP} satisfy the Cauchy
problem~\eqref{eq:1} in a distributional sense. Hence we still require
the following definition.
    
\begin{definition}[Distributional solutions]\label{def:distributional-L1}
 A function $u \in C([0,\infty);L^{1})$ is called a
 \emph{distributional} (in space) solution of
 Cauchy problem~\eqref{eq:1} if $u(0)=u_{0}$ in $L^{1}$, $u^{m}$ belongs
 to $L^{p}_{loc}((0,\infty);W^{s,p})$, and one has that
 \begin{align*}
    & -\int_{\R^{d}}u\,\xi\,\dx\Bigg\vert_{t_{1}}^{t_{2}}
     -\int_{t_{1}}^{t_{2}}\int_{\R^d}u\,\xi_t\,\dx\dt\\ 
     &\hspace{1cm} +\int_{t_{1}}^{t_{2}}\int_{\R^{2d}}\frac{(u^m(t,x)-u^m(t,y))^{p-1}
     (\xi(t,x)-\xi(t,y))}{|x-y|^{d+sp}}\td(x,y)\,\dt=0
 \end{align*}
for everz $0< t_{1}<t_{2}<\infty$, and every test function $\xi \in 
 C_c^{\infty}([0,\infty)\times\R^d)$. Furthermore, if $u$ is also differentiable with
$\partial_t u(t):=\tfrac{\td u}{\dt}(t) \in L^1$ for a.e.~$t > 0$, then we call this a
\emph{strong distributional} solution in $L^1$.
\end{definition}

The following theorems provide the well-posedness
of~\eqref{eq:abstract-CP} in the sense of mild solutions
and~\eqref{eq:1} in the sense of strong distributional solutions. Here
we write $[u]^{+}$ to denote the positive part $\max\{u,0\}$ of $u$, and
$[u]^{1}$ for $u$. This has been proven
in~\cite[Theorem~1.3]{doublyNonlinearPaper}.
 
\begin{theorem}[Well-posedness \& Comparison]
	\label{thm:well-posedness}
    Let $1 < p < \infty$, $0 < s < 1$ and $m>s$. Then, the following statements hold.
    \begin{enumerate}
        \item For every $u_{0}\in L^{1}$, there is a unique mild solution $u$ of the Cauchy problem~\eqref{eq:abstract-CP}.
        \item For every $u_{0,1}$, $u_{0,2}\in L^{1}$, the corresponding mild solutions $u_{1}$ and $u_{2}$ of~\eqref{eq:abstract-CP} satisfy
        \begin{equation}
        \label{eq:comparisonPrinciple}
            \norm{[u_{1}(t)-u_{2}(t)]^{\nu}}_{1}
            \le \norm{[u_{0,1}-u_{0,2}]^{\nu}}_{1}\qquad\text{for every $t\ge 0$}
       \end{equation}
        and for $\nu\in \{+,1\}$.
        \item For $\phi(r) := r^m$ for all $r \in \R$ with $m > s$,
        \begin{equation}
        \label{eq:densityResult}
            \overline{D((-\Delta_p)^s_{1\cap \infty}\varphi)}^{\mbox{}_{L^1}} = L^1.
        \end{equation}
      \end{enumerate}
\end{theorem}

Now, we turn to the smoothing effect of mild solutions of~\eqref{eq:abstract-CP}. This
result comes from the homogeneity of the doubly nonlinear operator and has been proven in \cite[Theorem~1.7]{doublyNonlinearPaper}.

\begin{theorem}
  \label{thm:strong-regularisation-effect}
  Let $1 < p < \infty$, $0 < s < 1$ and $m>0$. If $m(p-1) \neq 1$, then every mild solution $u$ in $L^1$ of the Cauchy problem~\eqref{eq:abstract-CP} is locally Lipschitz on $(0,\infty)$, satisfying
  \begin{displaymath}
      \limsup_{h\rightarrow 0^+}\frac{\norm{u(t+h)-u(t)}_{1}}{h}\le\frac{m(p-1)+2}{|m(p-1)-1|}\frac{\norm{u_0}_1}{t}
  \end{displaymath}
  for all $t > 0$.
\end{theorem}

Finally, we turn to the well-posedness of the Cauchy problem~\eqref{eq:1} in the sense of strong distributional solutions in $L^1$.

\begin{theorem}[Strong distributional solutions]
  \label{thm:regularisation-effect}
  Let $1 < p < \infty$, $0 < s < 1$ and $m>s$. Suppose that either $m(p-1) \neq 1$ or there exists $(v_{0,n})_{n \ge 1} \subset D((-\Delta_{p})^{s}_{1\cap\infty})$ such that $v_{0,n} \rightarrow u_0^m$ in $L^1$ as $n \rightarrow \infty$ and $\norm{(-\Delta_{p})^{s}_{1\cap\infty}v_{0,n}}_1$ is uniformly bounded. If either $u_{0}\in L^{1\cap\infty}$ or $u_{0} \in L^1$, $m\ge 1$ and $p$, $s$ and $m$ satisfy
    \begin{equation}
     \label{eq:conds-on-p-s-m}
        m(p-1) > 1-\frac{ps}{d},
    \end{equation}
    then every mild solution in $L^1$ of~\eqref{eq:abstract-CP} is also a strong distributional solution in $L^1$ of the Cauchy problem~\eqref{eq:1}.
\end{theorem}

\subsection{Proof of Theorem~\ref{thm:regularisation-effect}}
\label{subsec:proof-of-distributional-and-regularity}

To establish the existence of strong
distributional solutions of Cauchy problem~\eqref{eq:1}, we approximate the
mild solution $u$ of the Cauchy problem~\eqref{eq:abstract-CP} governed by the nonlocal doubly nonlinear operator $A$ in $L^1$ by strong distributional solutions $u_{n}$ in $L^{1}(\Omega_{n})$ of the initial boundary-value problem 
\begin{equation}
   \label{eq:1bis}
   \begin{cases}
    \begin{alignedat}{2}
    \partial_{t}u_{n}(t)+(-\Delta_{p})^{s}u^{m}_{n}(t)&=0 \quad
        &&\text{in $\Omega_{n} \times(0,\infty),$}\\
        u_{n}(t)&=0 \quad &&\text{in $\R^{d}\setminus\Omega_{n} \times[0,\infty),$}\\
        u_{n}(0)&=u_{0,n} \quad &&\text{on $\Omega_n$,}
    \end{alignedat}
    \end{cases}
\end{equation}
for $1 < p <\infty$ and $0 < s < 1$.
Here $(\Omega_{n})_{n\ge 1}$ is an increasing sequence of open and bounded subsets $\Omega_{n}$ in $\R^{d}$ and $(u_{0,n})_{n \ge 1}$ is a sequence of functions in $L^{1}(\Omega_n)$ converging to the $u_0$ of~\eqref{eq:abstract-CP} in $L^1$.
\medskip 
%$A_n: = A_{\Omega_n}$

To provide the rigorous functional analytical setting of this approximation argument, we first realize the fractional $p$-Laplacian $(-\Delta_{p})^{s}$ equipped with
homogeneous Dirichlet boundary conditions as an operator in $L^{2}(\Omega_{n})$ through
the sub-differential operator $\partial_{\! L^2(\Omega_{n})}\E_{n}$ of the energy
functional $\E_{n}:L^2(\Omega_{n}) \rightarrow (-\infty,+\infty]$ defined by
\begin{equation}
    \label{eq:energyFunctionaln}
    \E_{n}(u)=
    \begin{cases}
        \frac{1}{2p}[u]_{s,p}^{p} & \text{if $u\in W^{s,p}_{0}(\Omega_{n})\cap L^2(\Omega_{n})$,}\\
        +\infty & \text{otherwise.}
    \end{cases}
\end{equation}
It is not difficult to see (for example, refer to \cite{CoulHau2016,doublyNonlinearPaper}) that
the functional $\E_{n}$ is proper, convex, lower semicontinuous on $L^{2}(\Omega_{n})$, and 
G\^ateaux-differentiable in every direction $\xi\in D(\E_{n})=W^{s,p}_{0}(\Omega_{n})\cap
L^2(\Omega_{n})$. Thus, the sub-differential operator $\partial_{\! L^2(\Omega_{n})}\E_{n}$
can be characterized by
\begin{equation}
 \label{eq:subdiff-En}
    \partial_{\! L^2(\Omega_{n})}\E_{n}=\Big\{
    (u,v)\in D(\E_{n})\times L^2\,\Big\vert\, (u,v)\text{ satisfy \eqref{eq:variationalEq-Omega-n}}
    \Big\},
\end{equation}
where
\begin{equation}
  \label{eq:variationalEq-Omega-n}
   \int_{\Omega_{n}} v\xi\dx= \tfrac{1}{2}\int_{\R^{2d}}\tfrac{|u(x)-u(y)|^{p-2}(u(x)-u(y))
		(\xi(x)-\xi(y))}{|x-y|^{d+ps}}\,\td(x,y)
\end{equation}
for every $\xi\in D(\E_{n})=W^{s,p}_{0}(\Omega_{n})\cap
L^2(\Omega_{n})$.

\begin{definition}
  \label{def:fractional-p-laplacian-in-L2-Omega_n}
  For $1<p<\infty$ and $0<s<1$, let $\partial_{\! L^2(\Omega_{n})}\E_{n}$ be the sub-differential operator in $L^{2}(\Omega_{n})$ given by \eqref{eq:subdiff-En}. Then, we call the operator 
   \begin{displaymath}
    (-\Delta_{p})^{s}_{\Omega_{n}}u:=\partial_{\! L^2(\Omega_{n})}\E_{n}(u)\qquad 
    \text{for every $u\in D(\partial_{\! L^2(\Omega_{n})}\E_{n})$}
   \end{displaymath}
   the \emph{Dirichlet fractional $p$-Laplacian} in $L^2(\Omega_{n})$.
\end{definition}

Next, we realize of the doubly nonlinear nonlocal operator $(-\Delta_{p})^{s}\cdot^{m}$ on $L^1(\Omega_{n})$. For convenience, we write here $(L^{1\cap\infty}(\Omega_{n}))^2$ instead of $L^{1\cap\infty}(\Omega_{n})\times 
  L^{1\cap\infty}(\Omega_{n})$.

\begin{definition}
  \label{def:realization-in-L1-of-NDN-n}
  Let $(-\Delta_{p})^{s}_{\Omega_{n}}$ be the Dirichlet fractional $p$-Laplacian in $L^2(\Omega_{n})$
  and $(-\Delta_p)^s_{\Omega_{n},1\cap\infty}$ be the \emph{part of
  $(-\Delta_{p})^{s}_{\Omega_{n}}$ in $L^{1\cap\infty}(\Omega_{n})\times 
  L^{1\cap\infty}(\Omega_{n})$} given by
  \begin{equation}
     \label{A-operator-L1Linftyn}
    (-\Delta_{p})^{s}_{\Omega_{n},1\cap\infty}:=\Big\{(u,v)\in (L^{1\cap\infty}(\Omega_{n}))^2 %L^{1\cap\infty}(\Omega_{n})\times %L^{1\cap\infty}(\Omega_{n})
    \,\Big\vert\, v=(-\Delta_{p})^{s}_{\Omega_{n}}u
  %  \begin{array}[c]{l}
  %  u\in W^{s,p} \text{ and } (u,v) \text{ satisfies \eqref{eq:variationalEq}}\\
  %  \text{for all } \xi \in L^ \cap W^{s,p}
  %  \end{array}
    \Big\}.
\end{equation}
 Further, for $m>0$, let $\varphi(r):=\abs{r}^{m-1}r$ for every $r\in \R$. Then, we call
 the operator $A_{n}$ given by
 \begin{displaymath}
    A_{n}:=\Bigg\{(u,v)\in L^{1}(\Omega_{n})^2\,\Bigg\vert\, 
   \begin{array}[c]{l}
    \exists\, ((u_{k},v_{k}))_{k\ge 1}\;\text{ s.t. } v_k = (-\Delta_{p})^{s}_{\Omega_{n},1\cap\infty}\varphi(u_k),\\
    \displaystyle L^1-\lim_{k\to\infty}u_{k}=u\text{ and }
    L^1-\lim_{k\to\infty}v_{k}=v
   \end{array}
    \Bigg\}    
 \end{displaymath}
 the \emph{Dirichlet nonlocal doubly nonlinear operator} in $L^1(\Omega_{n})$, since it provides a realization in $L^1(\Omega_n)$ of the operator $(-\Delta_{p})^{s}\cdot^{m}$ equipped with homogeneous Dirichlet boundary condition.
\end{definition}

It was shown by the authors \cite{doublyNonlinearPaper} that for every $n\ge 1$, the
operator Dirichlet nonlocal doubly nonlinear operator $A_{n}$ is $m$-$T$ accretive in
$L^1(\Omega_{n})$, which means that for every $\lambda>0$, the resolvent operator
$J_{\lambda}^{A_{n}}:=(I+\lambda A_{n})^{-1}$ is a contraction mapping on 
$L^1(\Omega_{n})$ satisfying
\begin{equation}
    \label{eq:Taccretivity}
    \int_{\Omega}[w_1-w_2]^{+}\dx \le \int_{\Omega}[f_1-f_2]^{+}\dx
\end{equation}
for every $f_{1}$, $f_{2}\in L^1(\Omega_{n})$, where $w_1 = J_{\lambda}^{A_{n}}f_1$
and $w_2 = J_{\lambda}^{A_{n}}$, and
\begin{displaymath}
    \norm{J_{\lambda}^{A_{n}}f}_{q}\le \norm{f}_{q}
\end{displaymath}
for every $f\in L^q(\Omega_{n})$. In particular, the statements of 
Theorem~\ref{thm:well-posedness} and Theorem~\ref{thm:strong-regularisation-effect} hold for the Cauchy problem~\eqref{eq:abstract-CP}
(in $L^1$) replaced by the following one (in $L^1(\Omega_{n})$)
\begin{equation}
   \label{eq:abstract-CPn}
    \begin{cases}
     \begin{alignedat}{2}
	\tfrac{\td u}{\dt}(t)+A_{n}u(t)&\ni 0 &&\quad 
	\text{in $(0,\infty)$,}\\
	u(0)&=u_{0}. &&\mbox{}
    \end{alignedat}
   \end{cases}
\end{equation}
Moreover, we have a restricted version of Theorem~\ref{thm:regularisation-effect}. 
In particular, suppose $u_0 \in L^{1\cap\infty}$ such that there exists a sequence $(v_{0,n})_{n \ge 1}$ with $\lim\limits_{n\rightarrow\infty}v_{0,n} = u_0^m$ in $L^1$ and $\norm{(-\Delta_{p})^{s}_{\Omega_n,1\cap\infty}v_{0,n}}_1$ uniformly bounded.
Then the mild solution to \eqref{eq:abstract-CPn} is a strong distributional solution of the initial boundary-value problem~\eqref{eq:1bis} in $L^1$.

\begin{remark}
 It is worth noting that the notion of \emph{mild} and \emph{strong} solutions $u$ of Cauchy problem~\eqref{eq:abstract-CPn} as given in 
 Definition~\ref{def:mild-solutions} and Definition~\ref{def:strong-solutions} describe different differentiability properties with respect to the \emph{time} $t>0$. But since the operator $A_n$ given in Definition~\ref{def:realization-in-L1-of-NDN-n} is in fact the closure in $L^1(\Omega_n)$ of the doubly nonlinear operator $(-\Delta_{p})^{s}_{\Omega_{n},1\cap\infty}\varphi(\cdot)$, we require another notion of solution describing weak differentiability properties with respect to the \emph{space} variable.
\end{remark}

\begin{definition}[Distributional (in space) solutions]\label{def:distributional-n}
 For given $n\ge 1$ and $u_{0,n}\in L^{1}(\Omega_{n}))$, a function $u \in C([0,\infty);L^{1}(\Omega_{n}))$ is called a \emph{distributional} solution of
 Cauchy problem~\eqref{eq:1bis} if $u(0)=u_{0,n}$ in $L^{1}(\Omega_{n})$, $u^{m}$ belongs to
 $L^{p}_{loc}((0,\infty);W^{s,p}_{0}(\Omega_{n}))$ and for every test function $\xi \in 
 C_c^{\infty}([0,\infty)\times\Omega_{n})$, one has that
 \begin{align*}
    & -\int_{\Omega_{n}}u\,\xi\,\dx\Bigg\vert_{t_{1}}^{t_{2}}
     -\int_{t_{1}}^{t_{2}}\int_{\R^d}u\,\xi_t\,\dx\dt\\ 
     &\hspace{1cm} +\int_{t_{1}}^{t_{2}}\int_{\R^{2d}}\frac{(u^m(t,x)-u^m(t,y))^{p-1}
     (\xi(t,x)-\xi(t,y))}{|x-y|^{d+sp}}\td(x,y)\,\dt=0
 \end{align*}
for all $0< t_{1}<t_{2}<\infty$. Furthermore, if $u$ is also differentiable at a.e.~$t>0$
with $\partial_t u(t):=\tfrac{\td u}{\dt}(t) \in L^1(\Omega_{n})$, then we call this a
\emph{strong distributional} solution in $L^1(\Omega_{n})$.
\end{definition}

In order to compare solutions to~\eqref{eq:abstract-CPn} with solutions to~\eqref{eq:abstract-CP} on $\R^d$, we introduce the following extension to an operator on $L^{1}(\R^d)$. Let
\begin{equation}
    \tilde{A}_{n} u(x) = \begin{cases}
    [A_{n} (u\mathds{1}_{\Omega_n})](x)& \text{for $x \in \Omega_n$,}\\
	0& \text{otherwise}.
    \end{cases}
\end{equation}
This simply extends the operator by zero. So solutions to the Cauchy problem~\eqref{eq:abstract-CPn} can be extended by zero to satisfy
\begin{equation}
\label{eq:abstract-CPn-Rd}
\begin{cases}
\begin{alignedat}{2}
    \frac{\td u_{n}}{\dt}(t)+\tilde{A}_{n}u_{n}(t)&=0\quad && 
\text{in $(0,\infty)$,}\\
u_{n}(0)&=v_n\quad &&
\end{alignedat}
\end{cases}
\end{equation}
and solutions to~\eqref{eq:1bis} can be extended by zero to satisfy
\begin{equation}
\label{eq:1nRd}
\begin{cases}
\begin{alignedat}{2}
    \partial_{t}u_{n}(t)+(-\Delta_{p})^{s}_{\Omega_{n},1\cap\infty}(u_{n}(t))&=0\quad && 
\text{in $\R^d \times(0,\infty)$,}\\
u_{n}(0)&=v_n\quad && \text{in $\R^d$.}
\end{alignedat}
\end{cases}
\end{equation}

To prove the convergence of $(u_n)_{n \ge 1}$ as $\Omega_n$ grows to cover $\R^d$, we require the following two results. 
The following lemma gives us a comparison principle for the fractional $p$-Laplacian on subdomains which we will then apply in order to approximate solutions to~\eqref{eq:abstract-CP} by solutions to~\eqref{eq:1bis} or equivalently~\eqref{eq:1nRd}.

\begin{lemma}
\label{lem:comparisonSubdomain}
     Let $1 < p <\infty$, $0 < s <1$, $m > s$, $f \in L^1(\R^d)$ and $\Omega_{1}$, $\Omega_{2}$ be open subsets of $\R^d$ satisfying $\Omega_{1} \subseteq \Omega_{2}$. 
     Define $w_1 :=(I+\lambda\tilde{A}_{1})^{-1}(f\mathds{1}_{\Omega_{1}})$ and $w_2 :=(I+\lambda\tilde{A}_{2})^{-1}(f\mathds{1}_{\Omega_{2}})$ for any $\lambda > 0$. Suppose $f(x) \ge 0$ a.e.~on $\R^d$. 
     Then $w_1 \le w_2$ a.e.~on $\R^d$.
\end{lemma}

\begin{proof}
    By a standard density argument it is sufficient to consider $f \in L^1\cap L^\infty$ with $f \ge 0$ so that $w_1$, $w_2 \in L^1 \cap L^\infty$. In particular, then $w_1 \in D((-\Delta_p)^s_{\Omega_{1},1\cap\infty})$ and $w_2 \in D((-\Delta_p)^s_{\Omega_{2},1\cap\infty})$ when viewed as functions on $H_{+}$ and $H_{-}$ respectively, since they are both zero outside their respective domains. Moreover, by~\eqref{eq:Taccretivity}, since $f \ge 0$ we have that $w_1 \ge 0$ and $w_2 \ge 0$. We aim to estimate $(w_1-w_2)^{+}$, hence we introduce a sequence of $q \in C^1(\R)$ to approximate $\sign_0^{+}(w_1-w_2) := \mathds{1}_{\set{w_1 > w_2}}$. In particular, we require that $0 \le q \le 1$, $q(r) = 0$ for $r \le 0$ and that there exists $M > 0$ such that $0 < q'(r) < M$ for all $r > 0$. 
    For convenience, we define $Q := q(w_1-w_2)$. Note that $Q = 0$ in $\R^d\setminus H_{+}$ and
    \begin{align*}
        [Q]_{s,p}& = \frac{1}{p}\left(\int_{\R^d}\int_{\R^d} \frac{|q(w_1(x)-w_2(x))-q(w_1(y)-w_2(y))|^{p}}{|x-y|^{d+sp}}\dy\dx\right)^{\frac{1}{p}}\\
        & \le \frac{M}{p} \left(\int_{\R^d}\int_{\R^d} \frac{|w_1(x)-w_2(x)-(w_1(y)-w_2(y))|^{p}}{|x-y|^{d+sp}}\dy\dx\right)^{\frac{1}{p}}\\
        & \le M([w_1]_{s,p}+[w_2]_{s,p})
    \end{align*}
    so that $Q \in W_0^{s,p}(H_{+})$ and hence also in $W_0^{s,p}(H_{-})$. We have that
    \begin{align*}
        &\int_{\R^d}(w_1-w_2)q(w_1-w_2)\dx = -\lambda\int_{\R^d}(\tilde{A}_{1} w_1-\tilde{A}_{2}w_2)Q\dx\\
        & = -\lambda\int_{H_{+}}(-\Delta_p)^s_{\Omega_{1},1\cap\infty} (w_1^m) Q\dx+\lambda\int_{H_{-}}(-\Delta_p)^s_{\Omega_{2},1\cap\infty} (w_2^m)Q\dx\\
        & = -\lambda\int_{\R^{2d}}\frac{\left((w_1^m(x)-w_1^m(y))^{p-1}-(w_2^m(x)-w_2^m(y))^{p-1}\right)(Q(x)-Q(y))}{|x-y|^{d+sp}}\td(x,y)\\
        & \le 0.
    \end{align*}
    By a standard approximation for $\sign_0^{+}$, we then have that
    \begin{align*}
        \int_{\R^d}(w_1-w_2)^{+}\dx = 0
    \end{align*}
    so that $w_1 \le w_2$ a.e.~in $\R^d$.
\end{proof}

We can now prove the following approximation result for~\eqref{eq:abstract-CP}.

\begin{proposition}
   \label{prop:approxStrongSols}
    Let $d \ge 1$, $p > 1$, $0 < s < 1$, $m > s$ and suppose $(\Omega_n)_{n \ge 1}$ is the sequence of balls with radius $n$. 
    Then, for every $u_{0}\in L^{1}$, there is a sequence of functions $(u_{0,n})_{n \ge 1}$ with $u_{0,n} \in L^1(\Omega_n)$ and $u_{0,n}^m \in D((-\Delta_p)^s_{\Omega_{n},1\cap\infty})$ for each $n \ge 1$ such that
    \begin{displaymath}
    \lim_{n\rightarrow\infty}u_{0,n}\mathds{1}_{\Omega_n} = u_0\quad \text{in $L^1$}.
    \end{displaymath}
    Moreover, for each $n \ge 1$, there is a unique strong distributional solution $u_n$ of the Cauchy problem
    \begin{equation}
	\label{eq:1n}
	\begin{cases}
    \begin{alignedat}{2}
	\partial_t u_n(t)+(-\Delta_p)^s_{\Omega_{n},1\cap\infty} u_{n}(t)&=0 &&\quad 
	\text{in $\Omega_{n} \times(0,\infty)$,}\\
    u_{n}(0)&=u_{0,n} &&\quad \text{in $\Omega_{n}$}
    \end{alignedat}
	\end{cases}
    \end{equation}
    such that extending to $\R^d$ by zero,
    \begin{equation}
     \label{eq:limit-CTL1-of-un}
        \lim_{n\to\infty}u_{n}\mathds{1}_{\Omega_{n}}=u\qquad\text{in $C([0,T];L^{1})$}
    \end{equation}
    for every $T>0$, where $u$ is the unique mild solution of Cauchy problem~\eqref{eq:abstract-CP}.
\end{proposition}

\begin{proof}
 By~\cite[Theorem~1.3]{doublyNonlinearPaper}, for each $n \ge 1$, we can approximate $u_0\mathds{1}_{\Omega_n} \in L^1(\Omega_n)$ by a sequence $(u_{0,k})_{k \in \N}$ with $u_{0,k}^m \in D((-\Delta_p)^s_{\Omega_{n},1\cap\infty})$. 
 Taking a diagonal subsequence and relabelling, we obtain the sequence $(u_{0,n})_{n \ge 1}$ of the proposition. 
 Moreover, by \cite[Theorem 1.6]{doublyNonlinearPaper}, for each $n \ge 1$ there is a unique strong distributional solution $u_{n}$ of the initial value
 problem~\eqref{eq:1n} with initial data $u_{0,n}$. 
 Then, when extended by zero to be defined on $\R^d$, $u_{n}$ is also a strong distributional solution of~\eqref{eq:1nRd} and therefore a strong solution of~\eqref{eq:abstract-CPn-Rd}.

 Next, by \cite[Corollary~3.2]{doublyNonlinearPaper} (see also~\cite[Section~2]{CoulHau2016}), for every $\lambda>0$, the resolvent operators $J_{\lambda}^{\tilde{A}_{n}}:=(I+\lambda \tilde{A}_{n})^{-1}$ and $J_{\lambda}^{A}:=(I+\lambda A)^{-1}$ are contractions on $L^{1}$. Furthermore, by the $m$-accretivity of $A_n$ on $L^1(\Omega_n)$, $\tilde{A}_n$ is also $m$-accretive on $L^1$. 
 Then since $u_{n,0} \rightarrow u_0$ in $L^1$, applying \cite[Proposition~4.4 and Theorem~4.14]{MR2582280}, \eqref{eq:limit-CTL1-of-un} is equivalent to showing that for every $f\in L^{1\cap \infty}$ and every $\lambda>0$, for
 \begin{displaymath}
    w:=J_{\lambda}^{A}f\quad\text{and}\quad
    w_{n}:= J_{\lambda}^{\tilde{A}_{n}}f,
 \end{displaymath}
 one has that
 \begin{equation}
    \label{eq:resolvent-convergence}
        \lim_{n\to\infty} w_{n}=w\qquad\text{in $L^{1}$.}
 \end{equation}
 
 Hence we now aim to prove~\eqref{eq:resolvent-convergence}. 
 By definition of $\tilde{A}_n$, we can write
 \begin{displaymath}
    w_{n}(x)=
    \begin{cases}
    [J_{\lambda}^{A_{n}}(f\mathds{1}_{\Omega_{n}})](x) & \text{if $x\in \Omega_{n}$,}\\
    0 & \text{otherwise.}
   \end{cases}
 \end{displaymath}
 By \cite[Proposition~2.19]{CoulHau2016}, one has that $w^m \in D((-\Delta_p)^s_{1\cap\infty})$,  
 $w_{n}^m\in D((-\Delta_p)^s_{\Omega_{n},1\cap\infty})$ satisfying $\norm{w}_{1\cap \infty}\le \norm{f}_{1\cap\infty}$, and
    \begin{equation}
        \label{eq:L1-Linfty-bound}
        \norm{w_{n}}_{1\cap \infty}\le \norm{f}_{1\cap\infty}\qquad\text{for every $n\ge 1$.}
    \end{equation}
    Thus, multiplying 
    \begin{equation}
        \label{eq:approx-An}
        w_{n}\mathds{1}_{\Omega_{n}}+\lambda\, (-\Delta_p)^s_{\Omega_{n},1\cap\infty}(w_{n}\mathds{1}_{\Omega_{n}})^m=f\mathds{1}_{\Omega_{n}}
    \end{equation}
    by $w_{n}^{m}$ and subsequently applying~\eqref{eq:L1-Linfty-bound}, one sees that
    \begin{displaymath}
        \norm{w_{n}}^{m+1}_{m+1}+\lambda\,[w_{n}^{m}]_{s,p}^{p}
        \le \norm{f}_{1}\norm{f}_{\infty}^{m}. 
    \end{displaymath}
    From this and by~\eqref{eq:L1-Linfty-bound}, we can conclude that
    $(w_{n}^{m})_{n\ge 1}$ is bounded in $W^{s,p}$. Thus, there is a $v\in
    W^{s,p}$ such that, after possibly passing to a subsequence, and by setting $\tilde{w}^{m}:=v$, 
    \begin{equation}
        \label{eq:weak-limit-wmn}
        \lim_{n\to\infty}w_{n}^{m}=\tilde{w}^{m}\qquad
        \text{weakly in $W^{s,p}$.}
    \end{equation}
    Moreover, the fractional Rellich-Kondrachov theorem
    (see~\cite[Theorem~2.1]{MR4108436}) yields that after possibly passing to a subsequence, $w_{n}^{m}\to \tilde{w}^{m}$
    in $L^{q}(K)$ for every compact subset $K$ of $\R^{d}$, where $q \in [1,\infty]$ depends on $s$, $p$ and $d$. In particular, again possibly passing to a subsequence, one has that $w_{n}(x)\to \tilde{w}(x)$
    for a.e.~$x\in \R^{d}$.
    
    By Lemma \ref{lem:comparisonSubdomain}, supposing first that $f \ge 0$, we have that $w_n \le w_{n+1}$ a.e.~on $\R^d$. Thus and by~\eqref{eq:L1-Linfty-bound}, Beppo-Levi's
    theorem of monotone convergence yields that
    \begin{equation}
    \label{eq:resolvent-convergence-tilde}
        \lim_{n\to\infty} w_{n}=\tilde{w}\qquad\text{in $L^{1}$}
    \end{equation}
    holds. Similarly, if $f\le 0$ then for every integer $n \ge 1$, $0\le w_{1}-w_{n}\le w_{1}-w_{n+1}$ a.e. on $\R^{d}$ and by~\eqref{eq:L1-Linfty-bound}, it follows again from Beppo-Levi's theorem of monotone convergence that~\eqref{eq:resolvent-convergence-tilde} holds. Now, let $f\in L^{1\cap \infty}$ be general. Then we decompose $f=f^{+}-f^{-}$ into the positive part $f^{+}:=f\vee 0$ and negative part $f^{-}:=(-f)\vee 0$ of $f$, and set for every integer $n\ge 1$,
    \begin{displaymath}
       w_{n}^{+}(x):=\begin{cases}
            [J_{\lambda}^{A_{n}}[f^{+}\mathds{1}_{\Omega_{n}}]](x) & \text{if $x\in \Omega_{n}$,}\\
            0 & \text{otherwise,}
        \end{cases}
    \end{displaymath}
    and
    \begin{displaymath}
        w_{n}^{-}(x):=\begin{cases}
            [J_{\lambda}^{A_{n}}[-f^{-}\mathds{1}_{\Omega_{n}}]](x) & \text{if $x\in \Omega_{n}$,}\\
            0 & \text{otherwise.}
        \end{cases}
   \end{displaymath}
   Since $-f^{-}\le f\le f^{+}$, it follows from \eqref{eq:Taccretivity} that
    \begin{displaymath}
        w_{n}^{-}\le w_{n}\le w_{n}^{+}\qquad\text{a.e. on $\R^{d}$ for all
        $n\ge 1$}.
    \end{displaymath}
    Moreover, by the previous monotone convergence arguments with $f^{+}\ge 0$ and $-f^{+} \le 0$, $w_{n}^{+}\uparrow \tilde{w}^{+}$ and $w_{n}^{-}\downarrow \tilde{w}^{-}$ in $L^1$
    for some limits $\tilde{w}^{-}$, $\tilde{w}^{+}\in L^{1\cap\infty}$. Hence,
    $\tilde{w}^{+}+\abs{\tilde{w}^{-}}\in L^{1\cap\infty}$ and one has that
    $\abs{w_{n}}\le \tilde{w}^{+}+\abs{\tilde{w}^{-}}$ a.e. in $\R^d$ for all $n\ge 1$. Therefore, Lebesgue's dominated 
    convergence theorem implies that~\eqref{eq:resolvent-convergence-tilde} holds
    for every $f\in L^{1\cap \infty}$.
    
    Since $v^{\frac{1}{m}}=\tilde{w}$, it
    remains to show that $v^{\frac{1}{m}}=w:=J_{\lambda}^{A}f$. To see this, we first note that
    by~\eqref{eq:L1-Linfty-bound} and \eqref{eq:resolvent-convergence-tilde}, one has that for every $1\le q<\infty$, 
    \begin{equation}
        \label{eq:lim-in-Lq}
        \lim_{n\to\infty}w_{n}=\tilde{w}\qquad
        \text{in $L^{q}$}
    \end{equation}
    and so, multiplying~\eqref{eq:approx-An} by $w_{n}^{m}$ and subsequently, letting $n\to\infty$ gives that
    \begin{equation}
        \label{limit-wnm-wnm}
        \lim_{n\to \infty}[w_{n}^{m}]^{p}_{s,p}
        %=\lim_{n\to\infty}
        %\frac{\langle f,w_{n}^{m}\rangle_{W^{-s,p^{\mbox{}_{\prime}}},W^{s,p}}
        %-\norm{w_{n}}_{m}^{m}}{\lambda} 
        =\frac{\langle f,\tilde{w}^m\rangle_{W^{-s,p^{\mbox{}_{\prime}}},W^{s,p}}
        -\norm{\tilde{w}}_{m+1}^{m+1}}{\lambda}. 
    \end{equation}
    Further, since
    $(w_{n}^{m})_{n\ge 1}$ is bounded in $W^{s,p}$, one has that
    $((-\Delta_p)^{s}(w_{n}^{m}))_{n\ge 1}$ given by
    \begin{align*}
        &\langle (-\Delta_p)^{s}(w_{n}^{m}),\xi\rangle_{W^{-s,p^{\mbox{}_{\prime}}},W^{s,p}}\\
        &\qquad 
        =\int_{\R^{2d}}\frac{\abs{w_{n}^{m}(x)-w_{n}^{m}(y)}^{p-2}
        (w_{n}^{m}(x)-w_{n}^{m}(y))(\xi(x)-\xi(y))}{\abs{x-y}^{d+sp}}\dy\dx
    \end{align*}
    for every $\xi\in W^{s,p}$, is bounded in $W^{-s,p^{\mbox{}_{\prime}}}$.
    Therefore, there is a $\chi\in W^{-s,p^{\mbox{}_{\prime}}}$ such that after
    possibly passing to a subsequence, one has that 
    \begin{equation}
        \label{eq:limit-of-chi}
        \lim_{n\to\infty}(-\Delta_p)^{s}(w_{n}^{m})=\chi\qquad
        \text{weakly$\mbox{}^{\ast}$ in $W^{-s,p^{\mbox{}_{\prime}}}$.}
    \end{equation}
   By the two limits~\eqref{eq:lim-in-Lq} and~\eqref{eq:limit-of-chi}, multiplying the equation~\eqref{eq:approx-An} by $\xi \in C_c^{\infty}$ and subsequently, sending $n\to \infty$ yields that
    \begin{displaymath}
        \braket{\tilde{w}+\lambda\, \chi-f,\xi}_{W^{-s,p^{\mbox{}_{\prime}}},
        W^{s,p}}=0
    \end{displaymath}
    for all $\xi \in C_c^{\infty}$, from where a standard density argument yields that 
    %Approximating $\xi \in W^{s,p}$ by test functions gives that
    \begin{equation}
        \label{eq:w-lambda-chi-f}
        \tilde{w}+\lambda\, \chi=f\qquad
        \text{in $W^{-s,p^{\mbox{}_{\prime}}}$.}
    \end{equation}
    Therefore, it remains to show that
    \begin{equation}
     \label{eq:chi-equals-A1capinfty}
         \chi=(-\Delta_p)^s_{1\cap\infty}\tilde{w}^m.
    \end{equation}
    To prove this, we use the classical \emph{monotonicity trick} (see,
    e.g.,~\cite[p.~172]{MR0259693}). We begin by multiplying~\eqref{eq:w-lambda-chi-f} by $w^{m}_{n}$. Due
    to \eqref{eq:weak-limit-wmn},~\eqref{eq:lim-in-Lq} and \eqref{limit-wnm-wnm}, sending $n\to\infty$ in the resulting equation yields that
    \begin{equation}
        \label{eq:chi-wnm}
        \begin{split}
        \langle \chi,\tilde{w}^{m}\rangle_{W^{-s,p^{\mbox{}_{\prime}}},W^{s,p}}
        & =\frac{\langle f,\tilde{w}^m\rangle_{W^{-s,p^{\mbox{}_{\prime}}},W^{s,p}}
        -\norm{\tilde{w}}_{m+1}^{m+1}}{\lambda}\\
        & = \lim_{n\to \infty}[w_{n}^{m}]^{p}_{s,p}.
        \end{split}
    \end{equation}
    Next, let $\xi\in W^{s,p}$. Since the variational
    fractional $p$-Laplace operator is a monotone operator
    $(-\Delta_p)^{s} : W^{s,p}\to W^{-s,p^{\mbox{}_{\prime}}}$ in the sense that
    \begin{displaymath}
        \langle (-\Delta_p)^{s}(v)-(-\Delta_p)^{s}
        (\xi),v-\xi\rangle_{W^{-s,p^{\mbox{}_{\prime}}},W^{s,p}}\ge 0
    \end{displaymath}
    for every $v$, $\xi\in  W^{s,p}$, one has that 
    \begin{align*}
        0& \le \langle (-\Delta_p)^{s}(w_{n}^{m})-(-\Delta_p)^{s}
        (\xi),w_{n}^{m}-\xi\rangle_{W^{-s,p^{\mbox{}_{\prime}}},W^{s,p}}\\
        & = [w_{n}^{m}]^{p}_{s,p}-\langle (-\Delta_p)^{s}(w_{n}^{m}),\xi
        \rangle_{W^{-s,p^{\mbox{}_{\prime}}},W^{s,p}}-\langle (-\Delta_p)^{s}
        (\xi),w_{n}^{m}-\xi\rangle_{W^{-s,p^{\mbox{}_{\prime}}},W^{s,p}}
    \end{align*}
    for every $n\ge 1$. Hence, sending $n\to\infty$ in the previous inequality and using \eqref{eq:weak-limit-wmn}, \eqref{eq:limit-of-chi} and \eqref{eq:chi-wnm}
     gives that
    \begin{displaymath}
        0\le \langle \chi-(-\Delta_p)^{s}
        (\xi),\tilde{w}^{m}-\xi\rangle_{W^{-s,p^{\mbox{}_{\prime}}},W^{s,p}}.
    \end{displaymath}
    Since in the last inequality $\xi\in W^{s,p}$ was arbitrary, we can choose $\xi=\tilde{w}^{m}-\mu\,\zeta$ for any $\zeta\in W^{s,p}$ and $\mu>0$. It follows that
    \begin{displaymath}
        0\le \langle \chi-(-\Delta_p)^{s}
        (\tilde{w}^{m}-\mu\,\zeta),\zeta\rangle_{W^{-s,p^{\mbox{}_{\prime}}},W^{s,p}}.
    \end{displaymath}
    In this last inequality, we can send $\mu\to 0+$ and obtain that
    \begin{displaymath}
        0\le \langle \chi-(-\Delta_p)^{s}
        (\tilde{w}^{m}),\zeta\rangle_{W^{-s,p^{\mbox{}_{\prime}}},W^{s,p}}
        \qquad\text{for every $\zeta\in W^{s,p}$,}
    \end{displaymath}
    implying that $\chi=(-\Delta_p)^{s}(\tilde{w}^{m})$. Since $\tilde{w}$ and
    $f\in L^{1\cap\infty}$, it follows from~\eqref{eq:w-lambda-chi-f}
    that~\eqref{eq:chi-equals-A1capinfty} holds and so we have shown that
    $\tilde{w}+\lambda A\tilde{w}=f$, or, equivalently,
    $\tilde{w}=J_{\lambda}^A f$. Since the resolvent $J_{\lambda}^A$ is a contraction
    on $L^1$, it follows that $\tilde{w}=w$ and thereby we have shown
    that~\eqref{eq:resolvent-convergence} holds. This completes the proof of this lemma.  
 \end{proof} 

With the preceding lemma in mind, we can now give the proof of Theorem~\ref{thm:regularisation-effect}. Note, below we write $L^{p}(W^{s,p})$ instead of $L^{p}(0,T;W^{s,p})$ for $T>0$, and denote by $L^{p^{\mbox{}_{\prime}}}(W^{-s,p^{\mbox{}_{\prime}}})$ the dual space of $L^{p}(W^{s,p})$.

\begin{proof}[Proof of Theorem~\ref{thm:regularisation-effect}]
We only provide the proof for initial values $u_{0}\in L^{1\cap
\infty}$ since by Proposition~\ref{prop:L1LinfRegularity} the same conclusion holds for general $u_{0}\in L^1$ under the given restrictions on $p$, $s$, and $m$.
%Hence it remains to show that for 
%initial values $u_{0}\in L^{1\cap\infty}$, the unique mild solution $u$ of Cauchy problem~\eqref{eq:abstract-CP} is a strong
%distributional solution in $L^1$ of~\eqref{eq:1}. 
So let $u_{0}\in L^{1\cap\infty}$, take $(\Omega_n)_{n \ge 1}$ to be the sequence of balls with radius $n$ and $(u_{0,n})_{n\ge 1}$ the sequence of functions converging to $u_0$ in $L^1$ given by Proposition~\ref{prop:approxStrongSols}. 
Then define $u_n$ to be the associated sequence of strong distributional solutions to~\eqref{eq:1n}.
We now prove the strong distributional properties of the mild solution $u$ of Cauchy problem~\eqref{eq:abstract-CP} by considering the boundedness and convergence of $u_n$ and $\frac{du_n}{dt}$ in appropriate Banach spaces.

Multiplying the differential equation
in~\eqref{eq:1bis} by $u_{n}^{m}$ and subsequently, for every
$T>0$, integrating over $(0,T)$
yields that
    \begin{equation}
     \label{eq:Lp-Wsp-bound}
        \tfrac{1}{m+1}\norm{u_{n}(t)\mathds{1}_{\Omega_{n}}}_{m+1}^{m+1}
        \Big\vert_{0}^{T}
        +\int_{0}^{T}
        [u_{n}^{m}(t)]_{s,p}^{p}\,\dt=0.
    \end{equation}
    Since $u_{0}\in L^{1\cap \infty}$, Proposition~\ref{prop:approxStrongSols} and Proposition~\ref{prop:growthEstimate} imply that for every
    $1\le q<\infty$,
    \begin{equation}
        \label{eq:lim-un-in-Lq}
        \lim_{n\to\infty}u_{n}\mathds{1}_{\Omega_{n}}=u
        \qquad\text{in $C([0,T];L^{q})$}
    \end{equation}
    for every $T>0$. Thus \eqref{eq:Lp-Wsp-bound} implies
    that $(u^{m}_{n})_{n\ge 1}$ is bounded in
    $L^{p}(W^{s,p})$. Since the space
    $L^{p}(W^{s,p})$ is reflexive and
    by~\eqref{eq:lim-un-in-Lq}, we can conclude that $u^{m}
    \in L^{p}(W^{s,p})$ and after possibly passing
    to a subsequence, 
    \begin{equation}
        \label{eq:weak-lim-LpWsp}
        \lim_{n\to\infty}u_{n}^m\mathds{1}_{\Omega_{n}}=
        u^m\qquad\text{weakly in $L^{p}(W^{s,p})$.}
    \end{equation}
    In particular, one has that the sequence
    $(\mathcal{A}_{s,p}(u_{n}^{m}))_{n\ge 1}$ of linear bounded functionals
    given by
    \begin{align*}
        &\langle \mathcal{A}_{s,p}(u_{n}^{m}),\xi
        \rangle_{L^{p^{\mbox{}_{\prime}}}(W^{-s,p^{\mbox{}_{\prime}}}),
        L^{p}(W^{s,p})}\\
        &\quad 
        =\frac{1}{2}\int_{0}^{T}\int_{\R^{2d}}
        \frac{\abs{u_{n}^{m}(x)-u_{n}^{m}(y)}^{p-2}
        (u_{n}^{m}(x)-u_{n}^{m}(y))(\xi(x)-\xi(y))}{\abs{x-y}^{d+sp}}
        \,\td(x,y)\,\dt
    \end{align*}  
    for every $\xi\in L^{p}(W^{s,p})$, is bounded in
    $L^{p^{\mbox{}_{\prime}}}(W^{-s,p^{\mbox{}_{\prime}}})$.
    Therefore, there is an $\chi\in 
    L^{p^{\mbox{}_{\prime}}}(W^{-s,p^{\mbox{}_{\prime}}})$ such that after
    possibly passing to a subsequence, one has that 
    \begin{equation}
        \label{eq:limit-of-chiX}
        \lim_{n\to\infty}\mathcal{A}_{s,p}(u_{n}^{m})=\chi\qquad
        \text{weakly$\mbox{}^{\ast}$ in $L^{p^{\mbox{}_{\prime}}}(W^{-s,p^{\mbox{}_{\prime}}})$.}
    \end{equation}
    Further, it follows from the equation 
    \begin{equation}
      \label{eq:approx-evol-eq}
        \tfrac{\td u_{n}}{\dt}(t)+A_{n}u_{n}(t)=0\quad \text{in } \Omega_n\times(0,\infty)
    \end{equation}
    that $(\tfrac{\td u_{n}}{\dt}\mathds{1}_{\Omega_{n}})_{n\ge 1}$ is bounded in $L^{p^{\mbox{}_{\prime}}}(W^{-s,p^{\mbox{}_{\prime}}})$. 
    Thus and by~\eqref{eq:lim-un-in-Lq}, we can conclude that $\tfrac{\td u}{\dt}\in L^{p^{\mbox{}_{\prime}}}(W^{-s,p^{\mbox{}_{\prime}}})$ and after possibly passing to another subsequence of $(u_{n})_{n\ge 1}$ that
    \begin{equation}
        \label{eq:limit-of-partial-t-u-n}
        \lim_{n\to\infty}\tfrac{\td u_{n}}{\dt}\mathds{1}_{\Omega_{n}} =\tfrac{\td u}{\dt}\qquad
        \text{weakly$\mbox{}^{\ast}$ in $L^{p^{\mbox{}_{\prime}}}(W^{-s,p^{\mbox{}_{\prime}}})$.}
    \end{equation}

    Hence, multiplying~\eqref{eq:approx-evol-eq} by a test function $\xi \in C_c^{\infty}((0,T)\times \R^d)$ and sending $n\to \infty$ leads to
    \begin{align*}
        \left\langle\tfrac{\td u}{\dt}(t)+\chi,\xi\right\rangle_{L^{p^{\mbox{}_{\prime}}}(W^{-s,p^{\mbox{}_{\prime}}}),L^p(W^{s,p})} = 0.
    \end{align*}
    Therefore, we have shown that
    \begin{equation}
      \label{eq:evol-eq}
        \tfrac{\td u}{\dt}(t)+\chi=0\qquad\text{in $W^{-s,p^{\mbox{}_{\prime}}}$ for a.e. $t\in (0,T)$.}
    \end{equation}
    Now we can show that $u$ is a distributional solution of~\eqref{eq:1} by proving that $\chi=\mathcal{A}_{s,p}(u^{m})$. Here, $\mathcal{A}_{s,p}$ is the lifted operator $\mathcal{A}_{s,p} : 
    L^{p}(W^{s,p})\to L^{p^{\mbox{}_{\prime}}}(W^{-s,p^{\mbox{}_{\prime}}})$ given by 
    \begin{align*}
        &\langle \mathcal{A}^{s}_{p}(v),\xi\rangle_{L^{p^{\mbox{}_{\prime}}}(W^{-s,p^{\mbox{}_{\prime}}}),L^{p}(W^{s,p})}\\
        &\qquad 
        =\frac{1}{2}\int_{t_{1}}^{t_{2}}\int_{\R^{2d}}
        \frac{\abs{v(x)-v(y)}^{p-2}
        (v(x)-v(y))(\xi(x)-\xi(y))}{\abs{x-y}^{d+sp}}\,\td(x,y)\,\dt
    \end{align*}
    for every $v$, $\xi\in L^{p}(W^{s,p})$. To prove this, we again use the monotonicity trick (cf, the proof of Proposition~\ref{prop:approxStrongSols}), but to the lifted operator $\mathcal{A}^{s}_{p}$. First, we note that multiplying~\eqref{eq:evol-eq} by $u^{m}$ yields that
    \begin{displaymath}
        \tfrac{1}{m+1}\norm{u}_{m+1}^{m+1}\Big\vert_{0}^{T}+\langle \chi,u^{m}\rangle_{L^{p^{\mbox{}_{\prime}}}(W^{-s,p^{\mbox{}_{\prime}}}),L^{p}(W^{s,p})}=0.
    \end{displaymath}
    Thus \eqref{eq:Lp-Wsp-bound} and~\eqref{eq:lim-un-in-Lq} imply that
    \begin{equation}
        \label{eq:lim-lift-unm}
        \lim_{n\to\infty}\int_{0}^{T}[u_{n}^{m}(t)]_{p,s}^{s}\,\dt=
        \langle \chi,u^{m}\rangle_{L^{p^{\mbox{}_{\prime}}}(W^{-s,p^{\mbox{}_{\prime}}}),L^{p}(W^{s,p})}.
    \end{equation}
    By the monotonicity of $\mathcal{A}_{s,p}$, one has that
    \begin{align*}
      0\le &\langle \mathcal{A}^{s}_{p}(u_{n}^{m})-\mathcal{A}_{s,p}(\xi),u_{n}^{m}-\xi\rangle_{L^{p^{\mbox{}_{\prime}}}(W^{-s,p^{\mbox{}_{\prime}}}),L^{p}(W^{s,p})}\\
        =&\int_{0}^{T}[u_{n}^{m}(t)]_{s,p}^{p}\,\dt- \langle \mathcal{A}^{s}_{p}(u_{n}^{m}),\xi
        \rangle_{L^{p^{\mbox{}_{\prime}}}(W^{-s,p^{\mbox{}_{\prime}}}),L^{p}(W^{s,p})}\\
        &\hspace{3cm} -
        \langle\mathcal{A}_{s,p}(\xi),u_{n}^{m}-\xi\rangle_{L^{p^{\mbox{}_{\prime}}}(W^{-s,p^{\mbox{}_{\prime}}}),L^{p}(W^{s,p})}
    \end{align*}
    for all $n\ge 1$. Thus, sending $n\to\infty$ in the last inequality and by using~\eqref{eq:weak-lim-LpWsp}, \eqref{eq:limit-of-chiX}, and \eqref{eq:lim-lift-unm}, one finds that
    \begin{displaymath}
         0\le \langle \chi-\mathcal{A}_{s,p}(\xi),u^{m}-\xi\rangle_{L^{p^{\mbox{}_{\prime}}}(W^{-s,p^{\mbox{}_{\prime}}}),L^{p}(W^{s,p})}
    \end{displaymath}
    for every $\xi\in L^{p}(W^{s,p})$. Now, by proceeding as in the proof of Proposition~\ref{prop:approxStrongSols}, we can choose $\xi=u^{m}-\mu\,\zeta$, taking $\mu \rightarrow 0+$, to conclude that $\chi=\mathcal{A}_{s,p}(u^{m})$.
    
    Next, owing to~\cite[Theorem~1.6]{doublyNonlinearPaper}, if one
    multiplies~\eqref{eq:approx-evol-eq} by $\tfrac{\td u^{m}_{n}}{\dt}$, then one sees that 
    \begin{displaymath}
     \frac{\td}{\dt}\frac{1}{p}[u^{m}_{n}(t)]^{p}_{s,p}=-\frac{4m}{(m+1)^{2}}\norm{\tfrac{\td}{\dt}u_{n}^{
        \frac{m+1}{2}}\mathds{1}_{\Omega_{n}}}_{2}^{2}\le 0    
    \end{displaymath}
    for a.e. $t>0$. Thus, the function $t\mapsto \tfrac{1}{p}[u^{m}_{n}(t)]^{p}_{s,p}$ is monotonically decreasing along $(0,\infty)$. Subsequently integrating over $(t_{1},t_{2})$ for given $0<t_1<t_2<\infty$, one finds that
    \begin{equation}
     \label{eq:L2-bdds}
        \frac{4m}{(m+1)^{2}} \int_{t_{1}}^{t_{2}}\norm{\tfrac{\td}{\dt}u_{n}^{
        \frac{m+1}{2}}\mathds{1}_{\Omega_{n}}}_{2}^{2}\,\dt +\frac{1}{p}
        [u^{m}_{n}(t_{2})]^{p}_{s,p}=\frac{1}{p}[u^{m}_{n}(t_{1})]^{p}_{s,p}
  \end{equation}
  for every $n\ge 1$. Since $t\mapsto \tfrac{1}{p}[u^{m}_{n}(t)]^{p}_{s,p}$ is monotonically decreasing and
  by~\eqref{eq:Lp-Wsp-bound}, it follows that
  \begin{displaymath}
   [u^{m}_{n}(t_{1})]^{p}_{s,p}\le \frac{1}{t_1}\int_{0}^{t_1}[u^{m}_{n}(t)]^{p}_{s,p}\,\dt= 
   \frac{1}{m+1}\norm{u_{n}\mathds{1}_{\Omega_{n}}}_{m+1}^{m+1}
   \Big\vert_{t_1}^{0}.   
  \end{displaymath}
  From this and
  by~\eqref{eq:lim-un-in-Lq}, we can deduce that the sequence $(u^{m}_{n}(t))_{n\ge 1}$
  is bounded in $W^{s,p}$ for every $t>0$, and that
  $(\tfrac{\td}{\dt}u^{\frac{m+1}{2}}_{n}\mathds{1}_{\Omega_{n}})_{n\ge 1}$ is
  bounded in $L^{2}(t_1,t_2;L^2)$ for every $0<t_1<t_2<\infty$. Then 
  by~\eqref{eq:lim-un-in-Lq} and Fatou's lemma, we can conclude that
  $\tfrac{\td}{\dt}u^{\frac{m+1}{2}}\in L^{2}_{loc}((0,\infty);L^2)$ and $u^{m}(t)\in W^{s,p}$ for every $t>0$. In particular, one has that $w:=u^{\frac{m+1}{2}}\in 
  W^{1,1}_{loc}((0,\infty);L^{1}_{loc})$. 
  %, and after possibly passing to another subsequence that
  %\begin{align*}
  %    \lim_{n\to\infty}\partial_{t}u^{
  %\frac{m+1}{2}}_{n}\mathds{1}_{\Omega_{n}}
  %    &=\partial_{t}u^{\frac{m+1}{2}}\qquad\text{weakly in 
  %$L^{2}((0,T),t\,\dt;L^2)$, and}\\
  %    \lim_{n\to\infty}u^{m}_{n}(T)
  %    &=u^{m}(T)\qquad\text{weakly in $W^{s,p}$.}
  %\end{align*}
  %By these two limits, we can send $n\to\infty$ in~\eqref{eq:L2-bdds}. Then, %we obtain that
  By \cite[Theorem~1.7]{doublyNonlinearPaper}, if $m(p-1) \neq 1$ then $u$ is locally Lipschitz
  continuous on $(0,\infty)$ with values in $L^1$. 
  Otherwise we approximate $u$ and apply the Lipschitz estimate of \cite[Lemma 7.8]{BenilanNonlinearEvolutionEqns} to obtain local Lipschitz continuity in $L^1$ of $u$ on $(0,\infty)$.
  Therefore, and since
  $p(r):=r^{\frac{2}{m+1}-1}r\in L^{1}_{loc}(\R)$, it follows from
  \cite[Theorem~1.1]{BenilanGariepy} that $u\in W^{1,1}_{loc}
  ((0,T);L^{1}_{loc})$. 
  Now, applying the Lipschitz estimate on $u$, in the case $m(p-1) \neq 1$, one has that
  \begin{align*}
      \norm{\partial_{t}u(t)}_{L^{1}(\Omega_{n})}
      &= \limsup_{h\to 0+}\frac{\norm{u(t+h)-u(t)}_{L^{1}(\Omega_{n})}}{h}\\
      &\le \limsup_{h\to 0+}\frac{\norm{u(t+h)-u(t)}_{L^{1}}}{h}\\
      &\le \frac{\norm{u_{0}}}{t}\frac{2}{\abs{m(p-1)-1}}
  \end{align*}
 for every $t>0$ and $n\ge 1$. Sending $n\to \infty$ in this inequality shows that $u\in W^{1,1}_{loc}((0,\infty);L^{1})$, and in particular a strong solution of~\eqref{eq:1}. Similarly in the case $m(p-1) = 1$.
 The further regularity on $\partial_t u$ follows from Proposition~\ref{prop:timeDerivativeEst}.
\end{proof} 
	
\section{A priori bounds and regularity properties}
	\label{sec:apriori}
	To prove the existence of the self-similar solutions introduced in Section~\ref{sec:selfsimilar} we first mention a few useful regularity estimates for solutions to~\eqref{eq:1}. 
	
    By the complete resolvent property of $(-\Delta_p)^s \cdot^m$ (again, see~\cite{CoulHau2016} and~\cite{doublyNonlinearPaper}) we have the following growth estimate.
	\begin{proposition}[\cite{doublyNonlinearPaper}]
	    \label{prop:growthEstimate}
		Let $u$ be a (mild) solution to~\eqref{eq:abstract-CP} with initial data  $u_0 \in L^q$ for some $1 \le q \le \infty$. Then we have that
		\begin{equation}
		\label{eq:non-increasing-Lq-norm}
		\begin{split}
		    \norm{u(t)}_q& \le \norm{u(s)}_q
		\end{split}
		\end{equation}
		for all $0 \le s \le t$.
	\end{proposition}
	
	We begin by mentioning the following global $L^1-L^{\infty}$ estimate from~\cite{doublyNonlinearPaper}. Here $q_s$ comes from the usual embedding of $W^{s,p} \hookrightarrow L^{q_s}$. Note that we have a limited range of $m$ and $p$.
	\begin{proposition}[\cite{doublyNonlinearPaper}]
		\label{prop:L1LinfRegularity}
		Let $d \ge 1$,  $0 < s < 1$, $m \ge 1$ and $p > p_c$.
        Define $q_s$ by
		\begin{displaymath}
		q_s = \begin{cases}
		\left(\frac{1}{p}-\frac{s}{d}\right)^{-1}& \text{if } p < \frac{d}{s},\\
		\tilde{q}& \text{if } p = \frac{d}{s},\\
		\infty& \text{if } p > \frac{d}{s},
		\end{cases}
		\end{displaymath}
        where $\tilde{q}$ is chosen such that $\tilde{q} > \max(p,1+\frac{1}{m},\frac{p}{m(p-1)})$. Then for $u$ a mild solution to~\eqref{eq:abstract-CP}, we have the $L^1-L^\infty$ estimate
		\begin{equation}
		\label{eq:aprioriL1Linf}
		\norm{u(t)}_\infty \le Ct^{-\alpha}\norm{u_0}_1^{\gamma}
		\end{equation}
		where $\alpha = \frac{1}{m(p-1)-\frac{p}{q_s}}$ and $\gamma = \frac{1-\frac{p}{q_s}}{m(p-1)-\frac{p}{q_s}}$.
	\end{proposition}
	
	By~\cite{BenilanRegularizingMR648452} (see also~\cite{Hauer2019},~\cite{doublyNonlinearPaper}) we can estimate the time derivative using the homogeneity of the operator $(-\Delta_p)^s \cdot ^m$. In particular, since this composed operator is homogeneous of order $m(p-1)$, we have the following estimates from~\cite{doublyNonlinearPaper}.
	\begin{proposition}[\cite{doublyNonlinearPaper}]
	\label{prop:timeDerivativeEst}
		Suppose $d \ge 1$, $m > 0$, $p > 1$, $0 < s < 1$ and let $\alpha$, $\gamma$ and $q_s$ be as in Proposition~\ref{prop:L1LinfRegularity}. If $m(p-1) \neq 1$ and $u$ is a strong solution to~\eqref{eq:abstract-CP} with initial data $u_0 \in L^1(\R^d)$, then one has that
		\begin{equation}
		\label{eq:aprioriTime}
		\lnorm{\frac{\td u}{\dt}(t)}_1 \le \frac{1}{t}\frac{2}{|m(p-1)-1|}\norm{u_0}_1
		\end{equation}
		for all $t > 0$. 
        Moreover, if $p > p_c$ and $u$ is a distributional solution to~\eqref{eq:1} with initial data $u_0 \in L^1(\R^d)$, then
		\begin{equation}
		\label{eq:energyEstimate}
		    [u^{m}(t))]^{p}_{s,p} \le C\,t^{-(1+m\alpha)} \norm{u_0}_1^{1+m\gamma}
		\end{equation}
		for all $t>0$.
	\end{proposition}
	
	%We have the following local H\"older continuity result %presented in~\cite{doublyNonlinearPaper} which uses the %elliptic estimate from~\cite{IannizzottoHolderRegularity}.
	%\begin{proposition}
	%    H\"older continuity requires a bounded domain!
	%\end{proposition}

%Before presenting uniqueness of this Barenblatt solution, we introduce the following refined dissipation estimate in $L^1$ for differences of solutions to \eqref{eq:heat-eq}.

We now prove an estimate on the integral
\begin{displaymath}
    \int_{t_1}^{t_2} \int_{\R^d}(-\Delta_p)^s u^m\, u^q\dx \dt
\end{displaymath}
where $m > s$ and $q \ge 0$. We will use this to prove mass conservation in Section~\ref{sec:massConservation}. 
We note that in the case $q = m$, we have equality in \eqref{eq:generalEnergyEst} since in this case $u^m \in W^{s,p}$ so we do not need to approximate.
In the case $q = 0$ we let $r^0$ denote $\sign_0(r)$ for $r \in \R$.

\begin{lemma}
    \label{lem:generalEnergyEst}
    Let $m > 0$, $q \ge 0$ and suppose $u$ is a strong distributional solution to \eqref{eq:1} with initial data $u_0 \in L^1 \cap L^{\infty}$. 
    Then for all $0 < t_1 < t_2$,
    \begin{equation}
    \label{eq:generalEnergyEst}
    \begin{split}
        &\int_{t_1}^{t_2}\int_{\R^{2d}} \frac{(u^m(x,t)-u^m(y,t))^{p-1}(u^q(x,t)-u^q(y,t))}{|x-y|^{d+sp}}\td(x,y)\dt\\
        &\qquad \le \frac{1}{q+1}\left(\norm{u(t_1)}_{q+1}^{q+1}-\norm{u(t_2)}_{q+1}^{q+1}\right).
    \end{split}
    \end{equation}
\end{lemma}

\begin{proof}
   Due to the density result \eqref{eq:densityResult} and continuity of solutions given by \eqref{eq:comparisonPrinciple} of Theorem~\ref{thm:well-posedness}, we can approximate mild solutions with initial data in $L^1 \cap L^{q+1}$ by strong distributional solutions with initial data in $D((-\Delta_p)^s_{1\cap \infty}\cdot^m)$.
    Hence, we may assume without loss of generality that $u_0 \in L^{1}\cap L^{\infty}$ and so $u \in L^{\infty}(0,\infty;L^{1\cap\infty})$.

    We first approximate the power function $r^q: \R \rightarrow \R$ given by $r^q: r \mapsto |r|^{q-1}r$ by a sequence $(\phi_M)_{M \ge 1}$ of Lipschitz continuous functions on $\R$.
    If $q \ge 1$, we take
    \begin{displaymath}
        \phi_M(r) := \begin{cases}
            r^q & \text{ if } |r| \le M,\\
            \sign_0(r)\left(qM^{q-1}|r-M| + M^q\right) & \text{ if } |r| > M.
        \end{cases}
    \end{displaymath}
    If $0 \le q < 1$, we take
    \begin{displaymath}
        \phi_M(r) := \begin{cases}
            M^{1-q}r & \text{ if } |r| \le \frac{1}{M},\\
            r^q & \text{ if } |r| > \frac{1}{M}.
        \end{cases}
    \end{displaymath}
    Then for all $M \ge 1$, $|\phi_M(r)| \le |r|^q$ for $r \in \R$ and $\phi_M(r)$ converges to $r^q$ pointwise on $\R$ as $M \rightarrow \infty$.
    Since $\phi_M$ is Lipschitz continuous on $\R$, we have that for $v \in W^{s,p}\cap L^2$, $\phi_M(v) \in W^{s,p}\cap L^2$.
    
    For $t \in (t_1,t_2)$, we multiply \eqref{eq:1} by $\phi_M(u)$ and integrate to obtain
    \begin{displaymath}
        \int_{t_1}^{t_2}\int_{\R^d}(-\Delta_p)^s u^m\, \phi_M(u)\dx\dt = -\int_{t_1}^{t_2}\int_{\R^d}u_t\,\phi_M(u)\dx\dt.
    \end{displaymath}
    Since $u$ is a strong distributional solution and $\phi_M(u) \in W^{s,p}\cap L^2$, we can apply the variational formulation to the doubly nonlinear term, giving
    \begin{align*}
        &\int_{t_1}^{t_2}\int_{\R^{2d}} \frac{(u^m(x,t)-u^m(y,t))^{p-1}(\phi_M(u(x,t))-\phi_M(u(y,t)))}{|x-y|^{d+sp}}\td(x,y)\dt\\
        &\quad = -\int_{t_1}^{t_2}\int_{\R^d}u_t(x,t) \phi_M(u(x,t))\dx\dt.
    \end{align*}
    Since $\phi_M$ is monotonically increasing, the sign of the differences $u^m(x,t)-u^m(y,t)$ and $\phi_M(u(x,t))-\phi_M(u(y,t))$ are equal.
    Hence we have that the integrand on the left is positive. So taking $M \rightarrow \infty$ and applying Fatou's lemma, we estimate
    \begin{align*}
        &\int_{t_1}^{t_2}\int_{\R^{2d}} \frac{(u^m(x,t)-u^m(y,t))^{p-1}(u^q(x,t)-u^q(y,t))}{|x-y|^{d+sp}}\td(x,y)\dt\\
        &\quad \le \liminf\limits_{M\rightarrow\infty}-\int_{t_1}^{t_2}\int_{\R^d}u_t(x,t) \phi_M(u(x,t))\dx\dt\\
        &\quad = -\int_{t_1}^{t_2}\int_{\R^d}u_t(x,t) u^q(x,t)\dx\dt
    \end{align*}
    by dominated convergence.
    For this, we note that $u_t \in L^1(t_1,t_2;L^1)$ and $u \in L^{\infty}(t_1,t_2;L^{\infty})$, in particular so that $\norm{\phi_M(u)}_{\infty} \le \norm{u}_{\infty}^q$.
    Integrating by parts we obtain the result.
    %\comment{and applying the standard growth estimate,
    %\begin{align*}
      %  \left|\int_{t_1}^{t_2}\int_{\R^d}u_t(x,t) u^q(x,t)\dx\dt\right|& = \frac{1}{q+1}\left|\norm{u(t_2)}_{q+1}^{q+1}-
      %   \norm{u(t_1)}_{q+1}^{q+1}\right|\\
      %  & \le \frac{2}{q+1}\norm{u(t_1)}_{q+1}^{q+1}.
    %\end{align*}}
\end{proof}

\section{Proof of Aleksandrov's symmetry principle}
\label{sec:Aleksandrov}

In this section, we outline the proof of the comparison principle for
reflections (Theorem~\ref{thm:aleksandrovReflection}) and as a
consequence Aleksandrov's symmetry principle
(Theorem~\ref{cor:aleksandrov}) for mild solutions in $L^1$ of
Cauchy problem~\eqref{eq:abstract-CP}.\medskip

% Using the radial symmetry of the fractional $p$-Laplacian, we prove an
% Aleksandrov symmetry principle~\cite{MR0147776} for~\eqref{eq:1},
% sharpening an estimate from~\cite{MR4114983} in the case of the
% fractional $p$-Laplacian.  For this we prove a comparison principle
% between a solution and its reflection on the half-space
% $H_{+}$. The method applied here modifies a standard comparison
% principle idea, used for example in~\cite{MR2286292}, obtaining
% estimates on the whole domain.  This differs to the method used for
% the fractional Laplacian in~\cite{MR3191976} so that we may use the
% variational formulation of the fractional $p$-Laplacian (see
% also~\cite{MR1260981}and~\cite{MR4114983}).

% \begin{theorem}
%   \label{thm:aleksandrovReflection} Suppose $0 <s <1$, $p > 1$, $d \ge
% 1$ and $m > s$. Let $u$ be the mild solution of~\eqref{eq:abstract-CP}
% with initial data $u_0 \in L^1$. Let $\Pi:\R^d\rightarrow \R^d$ denote
% reflection around a hyperplane $H$, splitting $\R^d$ into $H_{+}$
% and $H_{-}$. If
%   \begin{align*} u_0(x) \le u_0(\Pi(x))\quad \text{in } H_{+},
%   \end{align*} then
%   \begin{align*} u(x,t) \le u(\Pi(x),t)\quad \text{in }
% H_{+}\times[0,T].
%   \end{align*}
% \end{theorem}

Before proving Theorem~\ref{thm:aleksandrovReflection}, we prove two
lemmas. The first allows us to use the variational formulation of the
fractional $p$-Laplacian when restricting to the half-space in
Theorem~\ref{thm:aleksandrovReflection}.

\begin{lemma}
  \label{lem:restrictionToHalfSpace}
  Let $\Pi:\R^d\rightarrow \R^d$ be the reflection around a hyperplane
  $H$ in $\R^d$ and $H_{+}$ be a half-space associated with $H$. If
  a measurable function $u$ satisfies $[u]_{s,p} < \infty$ and
  \begin{equation}
    \label{eq:reflectionProperty}
    |u(x)| \le |u(x)-u(\Pi(x))|
  \end{equation}
  for all $x \in H_{+}$ then
  $[u\mathds{1}_{H_{+}}]_{s,p} \le 2[u]_{s,p}$.
\end{lemma}
	
\begin{proof}
  We let $H_{-}$ be the other half-space given by reflecting
  $H_{+}$ around $H$. To estimate the seminorm
  \begin{equation}
    \label{eq:uIndOmega1}
    [u\mathds{1}_{H_{+}}]_{s,p}^p = \int_{\R^d}\int_{\R^d}\frac{|u(x)\mathds{1}_{H_{+}}(x)-u(y)\mathds{1}_{H_{+}}(y)|^p}{|x-y|^{d+sp}}\dx\dy,
  \end{equation}
  we consider the integrand for $(x,y) \in \R^d\times\R^d$, taking $x$
  and $y$ in each half-space given by the hyperplane $H$.
	    
  In the region $H_{+}\times H_{+}$, the integrand
  of~\eqref{eq:uIndOmega1} will correspond to the integrand of
  $[u]_{s,p}$ in the same region. Moreover, on
  $H_{-}\times H_{-}$ the integrand of~\eqref{eq:uIndOmega1} is
  zero. Hence we consider the case where $x$ and $y$ are in opposite
  regions, assuming without loss of generality that $x \in H_{+}$
  and $y \in H_{-}$. Then we can
  apply~\eqref{eq:reflectionProperty} to the integrand
  of~\eqref{eq:uIndOmega1}, with
	    \begin{align*}
	        \int_{H_{-}}&\int_{H_{+}}\frac{|u(x)\mathds{1}_{H_{+}}(x)-u(y)\mathds{1}_{H_{+}}(y)|^p}{|x-y|^{d+sp}}\dx\dy
	        = \int_{H_{-}}\int_{H_{+}}\frac{|u(x)|^p}{|x-y|^{d+sp}}\dx\dy\\
	        &\hspace{1cm} \le \int_{H_{-}}\int_{H_{+}}\frac{|u(x)-u(y)-(u(\Pi(x))-u(y))|^p}{|x-y|^{d+sp}}\dx\dy\\
	        &\hspace{1cm} \le 2^{p-1}\int_{H_{-}}\int_{H_{+}}\frac{|u(x)-u(y)|^p+|u(\Pi(x))-u(y)|^p}{|x-y|^{d+sp}}\dx\dy.
	    \end{align*}
	    We can then estimate the second term by a change of variables,
	    \begin{align*}
	        \int_{H_{-}}\int_{H_{+}}\frac{|u(\Pi(x))-u(y)|^p}{|x-y|^{d+sp}}\dx\dy& = \int_{H_{-}}\int_{H_{-}}\frac{|u(x)-u(y)|^p}{|\Pi(x)-y|^{d+sp}}\dx\dy\\
	        & \le \int_{H_{-}}\int_{H_{-}}\frac{|u(x)-u(y)|^p}{|x-y|^{d+sp}}\dx\dy.
	    \end{align*}
	    Note that for the last inequality, choosing coordinates such that the hyperplane $H$ is given by the set $\left\{x \in \R^d: x_1 = 0\right\}$, we have
	    \begin{equation}
	    \label{eq:reflectedDistanceEst}
	        |\Pi(x)-y|^2 = |x_1+y_1|^2+\sum_{i=2}^{d}|x_i-y_i|^2 \ge |x-y|^2
	    \end{equation}
	    when $x_1$ and $y_1$ have the same sign and hence when $x$ and $y$ are in the same half-space. Combining the integrals on each region, we have the desired estimate.
	\end{proof}
	
	The next lemma is a typical key estimate for obtaining comparison principles as we will be able to transfer this estimate to the time derivative of two solutions. However, notably in this case, we restrict to the half-space $H_{+}$ which requires additional considerations compared to the full domain $\R^d$.
    For this lemma we note that we can approximate the positive part of the sign function $\sign_0^{+}$ by a sequence $\left(q_{M}\right)_{M > 0}$ where
    \begin{displaymath}
        q_M(s) = \begin{cases}
            1& \text{ if } s \ge \frac{1}{M},\\
            Ms& \text{ if } 0 < s < \frac{1}{M},\\
            0&  \text{ if } s \le 0.
        \end{cases}
    \end{displaymath}
 
	\begin{lemma}
    	\label{lem:comparisonOmega1Lemma}
    	Let $\Pi:\R^d\rightarrow \R^d$ denote reflection around a hyperplane $H$, splitting $\R^d$ into $H_{+}$ and $H_{-}$. Suppose $u \in L^{1\cap\infty}$ such that $u^m \in D((-\Delta_p)^s_{1\cap\infty})$ and let $\hat{u}(x):= u(\Pi(x))$ for all $x\in \R^d$.
        We consider $q \in C^{1}(\R)$ satisfying $0 \le q \le 1$, $q(s) = 0$ for $s \le 0$ and $0 < q'(s) < M$ for all $s > 0$ given $M > 0$. Then
    	\begin{equation}
    	\label{eq:approxPointwiseComparisonOmega1}
    		\int_{H_{+}}\left((-\Delta_p)^s u^m-(-\Delta_p)^s \hat{u}^m\right)q(u-\hat{u})\dx \ge 0.
    	\end{equation}
    	%Further supposing $(-\Delta_p)^s u^m \in L^1$,
        In particular, one has that
        \begin{equation}
    	\label{eq:pointwiseComparisonOmega1}
    		\int_{H_{+}}\left((-\Delta_p)^s u^m-(-\Delta_p)^s \hat{u}^m\right)\mathds{1}_{\set{u>\hat{u}}}\dx \ge 0.
    	\end{equation}
	\end{lemma}
	
	\begin{proof}
	    We let $v = u^m$, $\hat{v} = \hat{u}^m$.
        Note that
	    \begin{align*}
	        [q(v-\hat{v})]_{s,p}& = \frac{1}{p}\left(\int_{\R^d}\int_{\R^d} \frac{|q(v(x)-\hat{v}(x))-q(v(y)-\hat{v}(y))|^p}{|x-y|^{d+sp}}\dx\dy\right)^{\frac{1}{p}}\\
	        & \le M\left(\int_{\R^d}\int_{\R^d} \frac{|v(x)-\hat{v}(x)-(v(y)-\hat{v}(y))|^p}{|x-y|^{d+sp}}\dx\dy\right)^{\frac{1}{p}}\\
	        & \le 2M[v]_{s,p}.
	    \end{align*}
	    Let $Q(x)$ denote $q(v(x)-\hat{v}(x))$ for $x \in \R^d$. Since $q(s) = 0$ for $s \le 0$, we have either $q(v(x)-v(\Pi(x))) = 0$ or $q(v(\Pi(x))-v(x)) = 0$, and so $Q$ satisfies~\eqref{eq:reflectionProperty}. Then by Lemma~\ref{lem:restrictionToHalfSpace},  $[Q\mathds{1}_{H_{+}}]_{s,p} < \infty$ and so we apply the variational formulation,
		\begin{align*}
		&\int_{\R^d}\left((-\Delta_p)^s v(x)-(-\Delta_p)^s \hat{v}(x)\right)Q(x)\mathds{1}_{H_{+}}\dx\\
		&\hspace{2cm} = \int_{\R^d}\int_{\R^d}\frac{(v(x)-v(y))^{p-1}-(\hat{v}(x)-\hat{v}(y))^{p-1}}{|x-y|^{d+sp}}\times\\
		&\hspace{5cm}\qquad (Q(x)\mathds{1}_{H_{+}}(x)-Q(y)\mathds{1}_{H_{+}}(y))\dy\dx.
		\end{align*}
		
		To show that this integral is non-negative, consider $x$
                and $y$ in the two regions $H_{+}$ and $H_{-}$ on either
                side of the hyperplane $H$. First note that in
                $H_{+}\times H_{+}$, we have
		\begin{equation}
		    \label{eq:omega1case}
    		\int_{H_{+}}\int_{H_{+}}\frac{\left((v(x)-v(y))^{p-1}-(\hat{v}(x)-\hat{v}(y))^{p-1}\right)
                  (Q(x)-Q(y))}{|x-y|^{d+sp}}\dy\dx.
		\end{equation}
		
		Considering the sign of the integrand in
                \eqref{eq:omega1case} and noting that $q$ is
                non-decreasing, we have that 
		\begin{align*}
		&\sign\left((v(x)-v(y))^{p-1}-(\hat{v}(x)-\hat{v}(y))^{p-1}\right)\\
                  &\qquad = \sign\left(v(x)-\hat{v}(x)-(v(y)-\hat{v}(y))\right)\\
		& \qquad = \sign(Q(x)-Q(y))
		\end{align*}
                provided $Q(x) \neq Q(y)$.  So~\eqref{eq:omega1case} is
                non-negative, while in the case of $H_{-}\times H_{-}$,
                both indicator functions are zero. For the cross terms
                $H_{+}\times H_{-}$ and $H_{-}\times H_{+}$, we apply a
                change of variables, mapping $y$ to $\Pi(y)$ and $x$ to
                $\Pi(x)$, respectively. Note that
                $|x-\Pi(y)| = |\Pi(x)-y|$. For the sum of these terms,
                we have
		\begin{align*}
		    &\int_{H_{+}}\int_{H_{-}}\frac{\left((v(x)-v(y))^{p-1}-(\hat{v}(x)-\hat{v}(y))^{p-1}\right)Q(x)}{|x-y|^{d+sp}}\dy\dx\\
    		&\hspace{2cm} -\int_{H_{-}}\int_{H_{+}}\frac{\left((v(x)-v(y))^{p-1}-(\hat{v}(x)-\hat{v}(y))^{p-1}\right)Q(y)}{|x-y|^{d+sp}}\dy\dx\\
    		&\quad = \int_{H_{+}}\int_{H_{+}}\frac{\left((v(x)-\hat{v}(y))^{p-1}-(\hat{v}(x)-v(y))^{p-1}\right)(Q(x)+Q(y))}{|x-\Pi(y)|^{d+sp}}\dy\dx.
		\end{align*}
		If this is non-negative we are done, so consider the domain on which this integrand may be negative. 
        In particular, since $Q$ is only non-zero when $v-\hat{v} > 0$, this can only occur when $v(x) > \hat{v}(x)$ and $v(y) < \hat{v}(y)$ or when $v(y) > \hat{v}(y)$ and $v(x) < \hat{v}(x)$. 
        By symmetry of $x$ and $y$, we may consider only the first case. Let
        \begin{displaymath}
          \Omega_x = H_{+}\cap\set{v(x) > \hat{v}(x)}\quad\text{ and
          }\quad \Omega_y = H_{+}\cap\set{v(y) < \hat{v}(y)}.
        \end{displaymath}
        Then for $(x,y) \in \Omega_x \times \Omega_y$, we have that
		\begin{align*}
		    &\int_{\Omega_x}\int_{\Omega_y}\frac{\left((v(x)-\hat{v}(y))^{p-1}-(\hat{v}(x)-v(y))^{p-1}\right)(Q(x)+Q(y))}{|x-\Pi(y)|^{d+sp}}\dy\dx\\
    		&\quad \ge -\int_{\Omega_x}\int_{\Omega_y}\frac{\left((v(x)-v(y))^{p-1}-(\hat{v}(x)-\hat{v}(y))^{p-1}\right)(Q(x)-Q(y))}{|x-y|^{d+sp}}\dy\dx
		\end{align*}
		where we also use that $|x-\Pi(y)| \ge |x-y|$ for $x$, $y \in H_{+}$ as in~\eqref{eq:reflectedDistanceEst}. 
        Hence the sum of both cases are bounded by the original $H_{+}\times H_{+}$ term,~\eqref{eq:omega1case}. 
        So we have the desired non-negativity,
		\begin{displaymath}
		    \int_{\R^d}\left((-\Delta_p)^s v(x)-(-\Delta_p)^s \hat{v}(x)\right)q(v(x)-\hat{v}(x))\mathds{1}_{H_{+}}\dx \ge 0.
		\end{displaymath}
		Letting $q$ converge to $[\sign_0]^{+}$ and noting that $r^m$ is strictly increasing for $r \in \R$, we can replace $v$ by $u^m$, giving~\eqref{eq:pointwiseComparisonOmega1}.
	\end{proof}
	
	We finally apply this half-space estimate to~\eqref{eq:1} and prove Theorem~\ref{thm:aleksandrovReflection}.
	
	\begin{proof}[Proof of Theorem \ref{thm:aleksandrovReflection}]
	    Using the translation and rotational invariance of the operator $(-\Delta_p)^s \cdot^m$, we may assume that $H_{+} = \set{x \in \R^d: x_1 > 0}$. Denote $\hat{u}(x) := u(\Pi(x))$ for all $x \in \R^d$. 
        First suppose that $u$ is a strong distributional solution to~\eqref{eq:1} with initial data $u_0 \in L^1 \cap L^{\infty}$ and that $(q_M)_{M > 0}$ converges to the positive indicator function $[\sign_0]^{+}$. Then $\hat{u}(x,t) := u(\Pi(x),t)$ is also a strong distributional solution with initial data $\hat{u}_0$. By the chain rule and applying Lemma~\ref{lem:comparisonOmega1Lemma},
		\begin{align*}
		\frac{\td }{\td t}&\int_{H_{+}}(u(t)-\hat{u}(t))^{+}\dmu = \int_{H_{+}}(u'(t)-\hat{u}'(t))\mathds{1}_{\set{u>\hat{u}}}\dmu\\
		& = -\int_{H_{+}}\left((-\Delta_{p})^{s}u^m(t)-(-\Delta_{p})^{s}\hat{u}^m(t)\right)\mathds{1}_{\set{u>\hat{u}}}\dx\\
		& \le 0
		\end{align*}
        noting that the composition $(-\Delta_{p})^{s}u^m(t) \in L^{1}$.
		Then we have,
		\begin{equation}
		\label{eq:gronwallAleksandrovComp}
    		\begin{split}
    		\int_{H_{+}}(u(t)-\hat{u}(t))^{+}\dmu& \le \int_{H_{+}}(u_{0}-\hat{u}_{0})^{+}\dmu.
    		\end{split}
		\end{equation}
        We approximate $L^1$ initial data by functions $u_{0,k} \in L^{1}$ such that $u_{0,k}^m \in D((-\Delta_{p})^{s}_{1\cap\infty})$ to obtain \eqref{eq:gronwallAleksandrovComp} for mild solutions $u$ with initial data $u_0 \in L^1$. Applying~\eqref{eq:gronwallAleksandrovComp} when $u_0(x) \le u_0(\Pi(x))$ then gives the result.
	\end{proof}

As an important consequence, we can prove now Aleksandrov's symmetry
principle (Theorem~\ref{cor:aleksandrov}).
% are the following two corollaries giving us
% radial symmetry of solutions with radially symmetric initial data.

% \begin{corollary}[Aleksandrov's symmetry principle]
% 		\label{cor:aleksandrov} Suppose $0 <s <1$, $p > 1$, $d
% \ge 1$ and $m > s$. Then mild solutions $u$ in $L^1$ of Cauchy
% problem~\eqref{eq:abstract-CP} with non-negative, compactly supported
% initial data $u_{0}$ in a ball $B_R(0)$ are radially decreasing in
% space for all $|x| \ge R$ and $t \ge 0$.
% \end{corollary}
	
\begin{proof}[Proof of Theorem~\ref{cor:aleksandrov}] The statement of
  this theorem follows
from Theorem~\ref{thm:aleksandrovReflection} by taking hyperplanes
perpendicular to the radial direction and with distance at least $R$
from the origin. Then data closer to the origin is non-negative and
the data away from the origin zero.
\end{proof}

The next corollary follows similarly.

\begin{corollary}
  \label{cor:aleksandrov2} Suppose $0 <s <1$, $p > 1$, $d \ge 1$ and
$m > s$.  Then mild solutions $u$ of~\eqref{eq:abstract-CP} with
non-negative, radially symmetric and radially decreasing initial data
are radially symmetric and radially decreasing in space for all $x \in
\R^d$ and $t \ge 0$.
\end{corollary}

\section{Barrier construction}
\label{sec:barrier}
	
In this section, we produce global barriers for solutions $u$ to
\eqref{eq:1} with initial data which is bounded with compact support.
These barrier functions are radially symmetric, decreasing in $x$ and
have sufficient decay at infinity to be integrable in space for all
$t \ge 0$.  Moreover, we will use these barrier functions to construct
Barenblatt solutions to~\eqref{eq:1}.  Such barrier functions have
been proven for the fractional $p$-Laplacian case with $m = 1$ in
\cite{MR4114983} and \cite{MR4280519}.  We apply the same methods for
more general $m > s$.
    
    We show the existence of such global barriers to~\eqref{eq:1} in the range $0 < s < 1$, $m > s$ and $p > p_c$, defined by~\eqref{eq:pcDef}, corresponding to the homogeneity condition
    \begin{displaymath}
        m(p-1) > 1-\frac{sp}{d}.
    \end{displaymath}
    In particular, this matches the range of the self-similar scaling transformation~\eqref{eq:transform1} and associated self-similar solutions~\eqref{eq:selfSimilarForm} with $\beta > 0$ defined by~\eqref{eq:scalingFactors}.
    
    When considering these barriers, we split the range of $p$ into three regions. 
    Define $p_1$ to be the positive solution to
    \begin{displaymath}
        \frac{msp_1(p_1-1)}{1-m(p_1-1)} = d,
    \end{displaymath}
    corresponding to the homogeneity condition
    \begin{equation}
    \label{eq:globalBarrierCondition}
        m(p_1-1) = \frac{d}{d+sp_1}.
    \end{equation}
    Note that $1 < p_{m,c}< p_1 < 1+\frac{1}{m}$. 
    
    We separately consider the upper region where $p > p_1$, the critical case $p = p_1$ and the lower sublinear region $p_{m,c}< p < p_1$. 
    We find barriers with different rates of decay at infinity in each region. 
    Hence the critical case $p_1$ is the transition point between the decay regimes of $|z|^{-\frac{sp}{1-m(p-1)}}$ (for $p < p_1$) and $|z|^{-(d+sp)}$ (for $p > p_1$) since 
    \begin{displaymath}
        d+sp_1 = \frac{sp_1}{1-m(p_1-1)}.
    \end{displaymath}

    We first define the decay function $g$ by
    \begin{equation}
    \label{eq:decayg}
        g(r) := \begin{cases}
        r^{-d-sp}& \text{if $p > p_1$,}\\
        r^{-d-sp}\log(r)& \text{if $p = p_1$,}\\
        r^{-\frac{sp}{1-m(p-1)}}& \text{if $p_{m,c}< p < p_1$.}
        \end{cases}
    \end{equation}
    for all $r > 0$.
    Then for positive constants $A$, $C_1$, $C_2$, $R_1$ and $R_2$, we introduce the barrier function $H:\R^d \times [0,\infty) \rightarrow [0,\infty)$ in the following way. For $p \ge p_1$,
    \begin{equation}
		\label{eq:barrierfunH}
		H(x,t) = \begin{cases}
		A(t+1)^{-d\beta}& \text{if } |x(t+1)^{-\beta}| \le R_1,\\
		C_1 |x|^{-d}& \text{if } R_1 < |x(t+1)^{-\beta}| \le R_2,\\
		C_2 (t+1)^{-d\beta}g(|x|(t+1)^{-\beta})& \text{if } |x(t+1)^{-\beta}| > R_2.
		\end{cases}
    \end{equation}

    For $p_{m,c}< p < p_1$,
    \begin{equation}
		\label{eq:barrierfunHSubcritical}
		H(x,t) = \begin{cases}
		A(t+1)^{-d\beta}& \text{if } |x(t+1)^{-\beta}| \le R_2,\\
		C_2 (t+1)^{-d\beta}g(|x|(t+1)^{-\beta})& \text{if } |x(t+1)^{-\beta}| > R_2.
		\end{cases}
    \end{equation}

    We then have existence of barrier functions with the following theorem.
    
    \begin{theorem}[Global barrier]
	    \label{thm:globalbarrier}
		Let $d \ge 1$, $0 < s < 1$, $m > s$ and $p > 1$. 
        Suppose $u$ is a mild solution to~\eqref{eq:abstract-CP} with initial data $u_0 \in L^1(\R^d)$ which is bounded and has compact support. 
        Define $H:\R^d \times [0,\infty) \rightarrow [0,\infty)$ by~\eqref{eq:barrierfunH} or~\eqref{eq:barrierfunHSubcritical} as appropriate.		
        Then there exist positive constants $A$, $C_1$, $C_2$, $R_1$ and $R_2$ such that $|u(t)| \le H(t)$ a.e.~in $\R^d$ for $t \ge 0$.
	\end{theorem}
    
    To prove this we apply the methods presented in \cite{MR4114983} in the superlinear case and \cite{MR4280519} in the sublinear case.
    The proofs follow essentially by taking $\hat{p} := m(p-1)+1$, $\hat{s} := \frac{sp}{\hat{p}}$ and applying the relevant calculations to $(-\Delta_{\hat{p}})^{\hat{s}}$ in place of $(-\Delta_p)^s \cdot^m$.
    We can apply such a replacement despite not being able to directly combine the exponent $m$ with $p-1$ due to the nature of the estimates required for these barrier functions.
    In particular, focusing on the spatial coordinates, for a barrier function $H$ and a point $y_0 \in \R^d$, we will only need to estimate $(-\Delta_p)^s H(y_0)$ from below.
    So we can typically ignore the first term of the difference, given that $H$ is non-negative, i.e.
    \begin{displaymath}
        (H(y_0)^m - H(y)^m)^{p-1} \ge - H(y)^{m(p-1)}.
    \end{displaymath}
    Meanwhile the issue of $\hat{s} \ge 1$ can be resolved by instead applying the smoothness of $H^m$ and using the associated boundedness properties of $(-\Delta_p)^s H^m$, in particular those given by \cite{MR4030247}.

    The outline in each case will be to find super-solutions to the self-similar profile equation~\eqref{eq:selfSimilarEq},
    rescale these to obtain super-solutions to~\eqref{eq:1} and then apply a comparison principle on $\R^d\times[0,\infty)$.

    \subsection{Super-solutions to the profile equation}
    To obtain pointwise estimates of~\eqref{eq:selfSimilarEq}, we use the Cauchy principal value of $(-\Delta_p)^s$ which will agree with Definition~\ref{def:realization-in-L1-of-NDN} when both exist but will not necessarily have the same domain. 
    In particular, we introduce the following approximation for $\varepsilon > 0$,
    \begin{equation}
    \label{eq:approximateSingularIntegral}
        (-\Delta_p)^s_{\varepsilon} u(z) := \int_{\R^d\setminus B_{\varepsilon}(z)}\frac{(u(z)-u(y))^{p-1}}{|z-y|^{d+sp}}\dy
    \end{equation}
    for all $z \in \R^d$. Then $(-\Delta_p)^s_{\varepsilon}$ converges to the principal value as $\varepsilon \rightarrow 0$. 
    
    The following lemma provides a global bound with decay of order $|x|^{-d}$. 
    The idea is to refine this estimate for large $|x|$ such that the global barrier is integrable in space.
	
	\begin{lemma}
    \label{lem:intermediateBound}
    	Let $d \ge 1$, $p > 1$, $0 < s < 1$ and $m > s$. 
        Suppose $u$ is a mild solution to~\eqref{eq:abstract-CP} with initial data $u_0 \in L^1(\R^d)$ which is bounded with compact support. 
        Then for $C > 0$ sufficiently large, depending on $d$ and $u_0$, we have
    	\begin{equation}
    	\label{eq:aprioriDecay2}
    	    |u(x,t)| \le C|x|^{-d}
    	\end{equation}
    	for all $x \in \R^d$ and $t > 0$.
	\end{lemma}
	
	\begin{proof}
	    Given a ball $B$ centered at the origin and containing the support of $u_0$, we bound $|u_0|$ pointwise a.e.~by a smooth, radially symmetric and radially decreasing function $\hat{u}_0$ with compact support defined such that
        \begin{displaymath}
            \hat{u}_0(x) = \norm{u_0}_{\infty}
        \end{displaymath}
        for all $x \in B$. 
        %By \cite{MR4030247} and the smoothness of $\hat{u}_0^m$ in $\R^d$, $\hat{u}_0^m \in D((-\Delta_p)^s_{L^{1\cap\infty}})$.
        We let $\hat{u}(x,t)$ be the associated mild solution to~\eqref{eq:abstract-CP} in $L^1$. 
        Since $\hat{u}_0$ is radially symmetric and radially decreasing, by Corollary~\ref{cor:aleksandrov2}, $\hat{u}(x,t)$ is radially decreasing in $\R^d$ for all $t \ge 0$. 
        Since $\hat{u}(x,t)$ is monotone decreasing in the radial variable, there exists $M > 0$, depending only on $d$, such that
	    \begin{align*}
	        \hat{u}(x,t)|x|^d& \le M\int_{\left\{y \in \R^d:|y| < |x|\right\}}\hat{u}(y,t)\dy
	         \le M\norm{\hat{u}_0}_1
	    \end{align*}
	    for a.e.~$x \in \R^d$, where we have used the standard accretive growth estimate in $L^1(\R^d)$. Hence $\hat{u}(x,t) \le C|x|^{-d}$ in $\R^d\times[0,\infty)$ for some $C > 0$ depending on $d$ and $u_0$. 
     Similarly, we have that $-\hat{u} \ge -C|x|^{-d}$. So comparing $u$, $-u$ with $\hat{u}$ and $-\hat{u}$ by~\eqref{eq:comparisonPrinciple}, we have \eqref{eq:aprioriDecay2}.
	\end{proof}
    
    For $p \ge p_1$ we choose a profile barrier function such that the intermediate region will be bounded by Lemma~\ref{lem:intermediateBound}. 
    We choose
    \begin{equation}
    \label{eq:barrierfunGeneral}
        G(z) = \begin{cases}
        	A& \text{if } |z| \le R_1,\\
        	C_1 |z|^{-d}& \text{if } R_1 < |z| \le R_2,\\
        	C_2 g(|z|)& \text{if } |z| > R_2,
    	\end{cases}
    \end{equation}
    where $g(r)$ is defined by~\eqref{eq:decayg}.  
    In particular, $g$ is smooth, decreasing and such that $G(z)$ is integrable at infinity.
    We glue the separate regions by matching constants. 
    In particular, we require that
    \begin{equation}
	\label{eq:matchingConstants}
    	A = C_1 R_1^{-d}\quad \text{and}\quad
        C_2 = \frac{C_1}{R_2^{d} g(R_2)}.
     \end{equation}

    We then require the super-solution condition in the remaining regions. We first consider the near region since this is independent of $g$.

    \begin{lemma}[Near region for $p > 1$]
    \label{lem:nearRegionSuperlinear}
         Let $d \ge 1$, $p > 1$, $0 < s < 1$ and $m > 0$.
        Define $G:\R^d \rightarrow [0,\infty)$ by~\eqref{eq:barrierfunGeneral}
		where the positive constants $C_1$ and $R_1$ are chosen to satisfy~\eqref{eq:nearCondition}, $R_2 \ge R_1$ and $A$, $C_2$ are given by~\eqref{eq:matchingConstants}.
        Then 
        \begin{equation}
        \label{eq:supersolEquation}
    		(-\Delta_p)^s_{\varepsilon} G(z)^m-\beta \nabla\cdot(z G(z)) \ge 0
        \end{equation}
        for all $\varepsilon > 0$ sufficiently small and all $z \in \R^d$ with $|z| \le R_1$.
    \end{lemma}

    \begin{proof}
    We have that
    \begin{equation}
    \label{eq:GgradientEst}
        -\beta \nabla\cdot(z G(z)) = -\beta d A
    \end{equation}
    for $|z| \le R_1$.
    Since the integrand is non-negative,
    \begin{align*}
        (-\Delta_p)^s_{\varepsilon} G(z)^m & \ge \int_{\{|z-y|\ge 3R_1\}\cap\{|y| \ge 2R_1\}} \frac{(A^m-G(y)^m)^{p-1}}{|z-y|^{d+sp}}\dy
    \end{align*}
    for all $0 < \varepsilon < R_1$ and all $z \in \R^d$ satisfying $|z| \le R_1$. 
    Then using~\eqref{eq:matchingConstants},
    \begin{align*}
        (A^m-G(y)^m)^{p-1}& \ge A^{m(p-1)}(1-2^{-dm})^{p-1}.
    \end{align*}
    Integrating and relabelling $C > 0$ as needed,
    \begin{align*}
        (-\Delta_p)^s_{\varepsilon} G(z)^m & \ge C A^{m(p-1)} \int_{3R_1}^{\infty} r^{-sp-1}\dr\\
        & = CA^{m(p-1)}R_1^{-sp}.
    \end{align*}
    So to satisfy the super-solution equation~\eqref{eq:supersolEquation} in this region, comparing with~\eqref{eq:GgradientEst}, relabelling $C$ and rearranging, we want to satisfy
    \begin{displaymath}
    A^{m(p-1)-1} \ge CR_1^{sp}
    \end{displaymath}
    where $C > 0$ depends on $s$, $p$, $d$ and $m$. 
    Applying~\eqref{eq:matchingConstants} gives
    \begin{equation}
    \label{eq:nearCondition}
    C_1^{m(p-1)-1} \ge CR_1^{d(m(p-1)-1)+sp}.
    \end{equation}
    \end{proof}
	
	We now prove the super-solution condition in the most involved region, the far region. For this we first consider the case $p > p_1$, defining
    \begin{equation}
    \label{eq:barrierfunSuperlinear}
        G(z) = \begin{cases}
        	A& \text{if } |z| \le R_1,\\
        	C_1 |z|^{-d}& \text{if } R_1 < |z| \le R_2,\\
        	C_2 |z|^{-(d+sp)}& \text{if } |z| > R_2,
    	\end{cases}
    \end{equation}
    with matching conditions
    \begin{equation}
	\label{eq:matchingConstants1}
    	A = C_1 R_1^{-d}\quad \text{and}\quad   	C_2 = C_1 R_2^{sp}.
     \end{equation}
    Here we will make room around the singularity in~\eqref{eq:approximateSingularIntegral} at $y=z$ by considering $|z| \ge 2R_2$.
    The intermediate region where $G(z) = C_1|z|^{-d}$ is important to be able to do this.
	\begin{lemma}[The far region for $p > p_1$]
		\label{lem:farRegion}
		Let $d \ge 1$, $0 < s < 1$, $m > 0$ and $p > p_1$.
        Define $G:\R^d \rightarrow [0,\infty)$ by~\eqref{eq:barrierfunSuperlinear} where positive constants $C_1$, $R_1$, $R_2$ satisfy~\eqref{eq:farRegionInequalities} and $A$, $C_2$ are given by~\eqref{eq:matchingConstants1}. 
        Then
        \begin{equation}
        \tag{\ref{eq:supersolEquation}}
    		(-\Delta_p)^s_{\varepsilon} G(z)^m-\beta \nabla\cdot(z G(z)) \ge 0
        \end{equation}
        for all $\varepsilon > 0$ sufficiently small and all $z \in \R^d$ with $|z| \ge 2R_2$.
	\end{lemma}
	
	\begin{proof}
		Fix $z \in \R^d$ with $|z| \ge 2R_2$. Applying \eqref{eq:matchingConstants1}, we have that
		\begin{equation}
		\label{eq:gradientTermEstimate}
        \begin{split}
		-\beta\nabla\cdot(z G(z)) = \beta sp C_1 R_2^{sp} |z|^{-d-sp} > 0
        \end{split}
		\end{equation}
		which we will use to control the fractional $p$-Laplacian term. 
        Hence we estimate $(-\Delta_p)^s_{\varepsilon} G(z)^m$ by separating~\eqref{eq:approximateSingularIntegral} into regions corresponding to $G$.
		
		In the region $\left\{y \in \R^d:|y| \le R_1\right\}$, noting that $G \ge 0$ on $\R^d$ and $|z| \le 2|z-y|$, we have
		\begin{align*}
		\int_{\left\{y \in \R^d:|y| \le R_1\right\}} \frac{(G(z)^m-G(y)^m)^{p-1}}{|z-y|^{d+sp}}\dy& \ge -CR_1^d A^{m(p-1)}|z|^{-d-sp}\\
        & = -CC_1^{m(p-1)} R_1^{-d(m(p-1)-1)} |z|^{-d-sp}
		\end{align*}
		for some $C > 0$ depending on $d$.
		
		In the region $\left\{y \in \R^d : R_1 < |y| \le R_2\right\}$, relabelling the constant $C$ as needed, we similarly have
		\begin{displaymath}
		\begin{split}
		\int_{\left\{y \in \R^d:R_1 < |y| \le R_2\right\}} &\frac{(G(z)^m-G(y)^m)^{p-1}}{|z-y|^{d+sp}}\dy \ge -C\int_{R_1}^{R_2}\frac{(C_1 r^{-d})^{m(p-1)}}{(|z|/2)^{d+sp}}r^{d-1}\dr\\
		  &\hspace{1.85cm} = -CC_1^{m(p-1)} \int_{R_1}^{R_2}r^{-dm(p-1)+d-1}\dr\, |z|^{-d-sp}
		\end{split}
		\end{displaymath}
		where $C > 0$ depends on $d$, $s$, $p$ and $m$.
  
		In the region $\left\{y \in \R^d : |y| > |z|\right\}$, the integrand will be positive so this can be ignored. However the intermediate region $\left\{y \in \R^d : R_2 \le |y| \le |z|\right\}$ contains the singularity at $y = z$. Hence we further split this region by integrating the ball $B_{|z|/2}(z)$ centered at $z$ separately.
		
		We treat the annulus region
                \begin{displaymath}
                  D = \left\{y \in \R^d : R_2 \le |y| \le |z|\right\}\setminus B_{|z|/2}(z)
                \end{displaymath}
as in the previous cases, noting that $|z-y| \ge |z|/2$. Then again relabelling $C$,
		\begin{align*}
		\int_{D} \frac{(G(z)^m-G(y)^m)^{p-1}}{|z-y|^{d+sp}}\dy& \ge -\int_{D} \frac{\left(C_2 |y|^{-(d+sp)}\right)^{m(p-1)}}{(|z|/2)^{d+sp}}\dy\\
		& \ge -CC_1^{m(p-1)}R_2^{-d(m(p-1)-1)} |z|^{-d-sp}
		\end{align*}
		since $(d+sp)m(p-1) > d$ by assumption~\eqref{eq:globalBarrierCondition}, where $C > 0$ depends on $d$, $s$, $p$ and $m$. In particular,~\eqref{eq:globalBarrierCondition} ensures that this and the following estimate on $B_{|z|/2}(z)$ do not grow relative to~\eqref{eq:gradientTermEstimate} as $|z| \rightarrow \infty$.
		
		Finally, we consider the ball  centered around $z$, $B_{|z|/2}(z)$. By \cite[Lemma 3.6]{MR4030247}, since $G^m$ is $C^2$ without critical points in $B_{|z|/2}(z)$, the integral
        \begin{displaymath}
        \int_{B_{|z|/2}(z)\setminus B_{\varepsilon}(z)}\frac{(G(z)^m-G(y)^m)^{p-1}}{|z-y|^{d+sp}}\dy
        \end{displaymath}
        is bounded independently of $\varepsilon > 0$. To determine the dependence on $|z|$, we apply a rescaling, finding that
        \begin{displaymath}
        \begin{split}
        \left|\int_{B_{|z|/2}(z)\setminus B_{\varepsilon}(z)}\frac{(G(z)^m-G(y)^m)^{p-1}}{|z-y|^{d+sp}}\dy\right|& \le C C_2^{m(p-1)}|z|^{-(d+sp)m(p-1)-sp}\\
        & \le C C_1^{m(p-1)}R_2^{-d(m(p-1)-1)}|z|^{-(d+sp)}.
        \end{split}
        \end{displaymath}
        where $C > 0$ depends on $d$, $s$, $p$ and $m$, again using that $(d+sp)m(p-1) > d$ by assumption. 
        Note that this matches the previous estimate.
	
		For $G(z)$ to be a super-solution we require that~\eqref{eq:gradientTermEstimate} bounds all these terms together for all $|z| \ge 2R_2$. 
        We then multiply the estimate in each region by $|z|^{d+sp}$ and apply~\eqref{eq:matchingConstants} to $R_2$ to reduce variables. 
        Then the super-solution condition~\eqref{eq:supersolEquation} holds given that
		\begin{equation}
		\label{eq:farRegionInequalities}
		\begin{split}
    		C_1^{1-m(p-1)}&\ge CR_1^{-d(m(p-1)-1)}R_2^{-sp},\\
    		C_1^{1-m(p-1)}&\ge C\int_{R_1}^{R_2}r^{-d(m(p-1)-1)-1}\dr R_2^{-sp},\\
            C_1^{1-m(p-1)}&\ge CR_2^{-d(m(p-1)-1)-sp}
		\end{split}
		\end{equation}
		for $C > 0$ depending on $d$, $s$, $p$ and $m$.
	\end{proof}

We prove the critical case $p = p_1$ very similarly to the case $p > p_1$, but now with a correction term. Define $G:\R^d \rightarrow [0,\infty)$ by
\begin{equation}
\label{eq:barrierfunCritical}
    G(z) = \begin{cases}
        A& \text{if } |z| \le R_1,\\
        C_1 |z|^{-d}& \text{if } R_1 < |z| \le R_2,\\
        C_2 |z|^{-(d+sp_1)}\log(|z|)^{\gamma}& \text{if } |z| > R_2,
    \end{cases}
\end{equation}
for all $z \in \R^d$ where $\gamma = \frac{1}{1-m(p-1)}$. 
We then have the matching conditions,
\begin{equation}
\label{eq:matchingConstantsCritical}
    A = C_1 R_1^{-d}\quad \text{and}\quad   	C_2 = \frac{C_1 R_2^{sp_1}}{\log(R_2)^{\gamma}}.
 \end{equation}

\begin{lemma}[The far region for $p = p_1$]
    Let $d \ge 1$, $0 < s < 1$, $m > 0$ and $p = p_1$.
    Define $G:\R^d\rightarrow [0,\infty)$ by~\eqref{eq:barrierfunCritical} where positive constants $C_1$, $R_1$, $R_2$ satisfy~\eqref{eq:farRegionInequalitiesCritical} and $A$, $C_2$ are given by \eqref{eq:matchingConstantsCritical}. 
    Then
    \begin{equation}
    \tag{\ref{eq:supersolEquation}}
        (-\Delta_p)^s_{\varepsilon} G(z)^m-\beta \nabla\cdot(z G(z)) \ge 0
    \end{equation}
    for all $\varepsilon > 0$ sufficiently small and all $z \in \R^d$ with $|z| > \min\set{2R_2,2e^{\frac{\gamma}{sp_1}}}$.
\end{lemma}

\begin{proof}
    In radial coordinates with $r = |z|$,
    \begin{displaymath}
        -\beta r^{1-d}(r^d G(y))_r = C_2\beta\left(sp_1-\frac{\gamma}{\log(r)}\right)r^{-(d+sp_1)}\log(r)^{\gamma}.
    \end{displaymath}
    For $r > 2e^{\frac{\gamma}{sp_1}}$, we have
    \begin{displaymath}
        -\beta r^{1-d}(r^d G(z))_r \ge \frac{C C_1 R_2^{sp_1}}{\log(R_2)^{\gamma}} |z|^{-(d+sp_1)}\log(|z|)^{\gamma}
    \end{displaymath}
    for some $C > 0$ depending on $m$, $s$ and $d$.
    As usual we now estimate the fractional $p$-Laplacian term $(-\Delta_{p_1})^s_{\varepsilon} G(z)^m$ for $z \in \R^d$ with $|z| > \min\set{2R_2,2e^{\frac{\gamma}{sp_1}}}$. 
    These estimates are done similarly to those in Lemma~\ref{lem:farRegion}. 
    For $\{y \in \R^d:|y| \le R_2\text{ or } |y| \ge |z|\}$, the calculations are the same.
    
    However, in the intermediate region $\left\{y \in \R^d : R_2 \le |y| \le |z|\right\}$ we must account for the new decay.
    For the annulus region 
    \begin{displaymath}
    	D = \left\{y \in \R^d : R_2 \le |y| \le |z|\right\}\setminus B_{|z|/2}(z),
    \end{displaymath}	 
    we integrate by parts to obtain,
    \begin{displaymath}
    \int_{D} \frac{(G(z)^m-G(y)^m)^{p_1-1}}{|z-y|^{d+sp_1}}\dy \ge -CC_2^{m(p_1-1)}|z|^{-d-sp}\log(|z|)^{\gamma}
    \end{displaymath}
    since $(d+sp_1)m(p_1-1) = d$, where $C > 0$ depends on $d$, $s$ and $m$.
    
    For the ball $B_{|z|/2}(z)$, we again apply \cite[Lemma 3.6]{MR4030247} and a rescaling, finding that
    \begin{displaymath}
    \left|\int_{B_{|z|/2}(z)\setminus B_{\varepsilon}(z)}\frac{(G(z)^m-G(y)^m)^{p-1}}{|z-y|^{d+sp}}\dy\right| \le C C_2^{m(p-1)}|z|^{-d-sp}\log(|z|)^{\gamma-1}
    \end{displaymath}
    again noting that $(d+sp)m(p_1-1) = d$.
    
    For $G(z)$ to be a super-solution we require that~\eqref{eq:gradientTermEstimate} bounds all these terms together for all $|z| \ge 2R_2$. 
    Then, relabelling $C > 0$ as needed, we can multiply each condition by $|z|^{d+sp_1}\log(|z|)^{-\gamma}$ and rearrange to obtain the conditions
    \begin{equation}
    \label{eq:farRegionInequalitiesCritical}
        C_1 \ge CR_2^{d-sp_1\gamma},\quad
        C_1 \ge C\log(R_2)^{\gamma}R_2^{-sp_1},\quad
        C_1 \ge CR_2^{-sp_1}.
    \end{equation}
\end{proof}

We now consider the lower sublinear case $p_{m,c}< p < p_1$. Here we use a barrier function of the form
\begin{equation}
\label{eq:lowSublinearProfileBound}
    G(z) = \begin{cases}
        A& \text{if } |z| \le R_2,\\
        C_2 |z|^{-\frac{sp}{1-m(p-1)}}& \text{if } |z| > R_2,
    \end{cases}
\end{equation}
for $z \in \R^d$ with the matching condition
\begin{equation}
\label{eq:matchingConditionSublinear}
    A = C_2R_2^{-\frac{sp}{1-m(p-1)}}.
\end{equation}
Note that $p > p_c$ ensures that $G$ is integrable at infinity.

\begin{lemma}[The far region for $p_{m,c}< p < p_1$]
    Let $d \ge 1$, $0 < s < 1$, $m > 0$ and $p_{m,c}< p < p_1$. 
    Define $G:\R^d\rightarrow [0,\infty)$ by \eqref{eq:lowSublinearProfileBound} where $C_2 > 0$ satisfies~\eqref{eq:farRegionInequalitiesSublinear}, $R_2 > 0$ and $A$ is given by~\eqref{eq:matchingConditionSublinear}. 
    Then
    \begin{equation}
    \tag{\ref{eq:supersolEquation}}
        (-\Delta_p)^s_{\varepsilon} G(z)^m-\beta \nabla\cdot(z G(z)) \ge 0
    \end{equation}
    for all $\varepsilon > 0$ sufficiently small and all $z \in \R^d$ with $|z| > R_2$.
\end{lemma}

\begin{proof}
    In radial coordinates with $r = |z|$,
    \begin{displaymath}
        -\beta r^{1-d}(r^d G(y))_r = \frac{C_2}{1-m(p-1)}r^{-\frac{sp}{1-m(p-1)}} > 0.
    \end{displaymath}
    We use this to compensate for the (possibly negative) fractional $p$-Laplacian term. 
    Evaluating $(-\Delta_p)^s_{\varepsilon} G^m$, first consider $y_0 \in \R^d$ such that $|y_0| = 1$. Then,
    \begin{displaymath}
        (-\Delta_p)^s_{\varepsilon} G(y_0)^m = C_2^{m(p-1)}\int_{\R^d\setminus B_{\varepsilon}(y_0)}\frac{\left(1-r^{-\frac{spm}{1-m(p-1)}}\right)^{p-1}}{|y_0-y|^{d+sp}}\dy.
    \end{displaymath}
    The integrand is positive for $r > 1$, so this region can be ignored.
    For $p < p_1$, $G^{m(p-1)} \in L^1(B_1(0))$ so the integrand is also bounded in $L^1$ near $y = 0$.
    Hence it remains to estimate near $y = y_0$.
    Since $G(y)^m$ is a $C^2$ function without critical points, we can also apply \cite{MR4030247} so that the principal value is bounded in a small ball around $y_0$ with bound independent of $y_0$.
    So as in \cite{MR4280519}, there is a finite constant $k(d,s,p,m)$ such that
    \begin{displaymath}
        (-\Delta_p)^s G(y_0)^m = -kC_2^{m(p-1)}.
    \end{displaymath}
    We evaluate for $|y_0| \neq 1$ by the spatial scaling transformation
    \begin{displaymath}
        v_h(y) := h^{\gamma}v(hy)
    \end{displaymath}
    for $h > 0$ with $\gamma = \frac{sp}{1-m(p-1)}$.
    This leaves the profile bound \eqref{eq:lowSublinearProfileBound} invariant and scales the profile equation~\eqref{eq:selfSimilarProfileRadial} by $r^{-\gamma}$ such that
    \begin{displaymath}
        (-\Delta_p)^s G(y)^m -\beta r^{1-d}(r^d G(y))_r = \left(\frac{C_2}{1-m(p-1)}-kC_2^{m(p-1)}\right)r^{-\gamma}.
    \end{displaymath}
    Then we satisfy the super-solution condition \eqref{eq:supersolEquation} for
    \begin{equation}
    \label{eq:farRegionInequalitiesSublinear}
        C_2 \ge \left(k(1-m(p-1))\right)^{1/(1-m(p-1))}.
    \end{equation}
\end{proof}

\subsection{A comparison principle}

We now convert profile functions back to the standard coordinates of~\eqref{eq:1} and apply a comparison principle.

\begin{lemma}[Comparison principle]
\label{lem:comparisonSupersolution}
    Let $u$ be a mild solution to~\eqref{eq:1} with initial data $u_0 \in L^1$. Suppose $H:\R^d\times[0,\infty) \rightarrow \R$ satisfies the super-solution condition
    \begin{equation}
    \label{eq:partialSupersolution}
        (H'(t)+(-\Delta_p)_{\varepsilon}^s H^m(t))\mathds{1}_{\set{u>H}} \ge 0
    \end{equation}
    a.e.~in $\R^d$ for all $t > 0$ and all $\varepsilon > 0$ sufficiently small. 
    Suppose $|u_0| \le H(0)$ a.e.~in $\R^d$. 
    Then for all $t \ge 0$, $u(t) \le H(t)$ almost everywhere in $\R^d$.
\end{lemma}

\begin{proof}
    By Theorem \ref{thm:well-posedness}, we can approximate $u_0$ by smooth data to get strong distributional solutions converging to $u$ in $L^1(\R^d)$. 
    Hence first suppose that $u$ is a strong distributional solution of \eqref{eq:1}.

    We intend to estimate the signed difference, noting that $H$ is piecewise smooth,
    \begin{align*}
        &\frac{\td }{\td t}\int_{\R^d}(u(t)-H(t))^{+}\dmu = \int_{\R^d}(u'(t)-H'(t))\mathds{1}_{\set{u>H}}\dmu\\
        &\hspace{2cm} \le -\int_{\R^d}((-\Delta_p)^s u(t)+(-\Delta_{p})_{\varepsilon}^{s}H^m(t))\mathds{1}_{\set{u>H}}\dmu.
    \end{align*}
    By the accretivity of $(-\Delta_p)^s \cdot^m$ in $L^1$ and using that $\cdot^m$ is strictly increasing so that $\mathds{1}_{\set{u^m>H^m}} = \mathds{1}_{\set{u>H}}$, we have the monotonicity estimate
    \begin{displaymath}
        \int_{\R^d}\left((-\Delta_{p})^{s}u^m(t)-(-\Delta_{p})^{s}H^m(t)\right)\mathds{1}_{\set{u>H}}\dmu \ge 0.
    \end{displaymath}
    Hence it remains to prove the convergence of $(-\Delta_{p})_{\varepsilon}^{s}H^m(t)$.
    We can approximate $\mathds{1}_{u > H}$ by $q(u-H)$ for $q \in C^1(\R)$ as usual with $0 \le q \le 1$, $q(r) = 0$ for $r \le 0$, $0 < q'(r) < M$ for some $M > 0$ for all $r > 0$. 
    Then letting $Q(t) = q(u(t)-H(t))$ for $0 < t < \infty$ and noting that $H \in W^{s,p}$ and $Q \in W^{s,p}$, we can use the symmetry of the fractional p-Laplacian, fixing $t > 0$, to obtain
    \begin{align*}
       & \int_{\R^d}\left((-\Delta_p)^{s}_{\varepsilon}H^m\right)Q\dmu \\
      & \qquad = \int_{\R^{2d}\setminus\{|z-y|<\varepsilon\}}\frac{(H^m(z)-H^m(y))^{p-1}Q(z)}{|z-y|^{d+sp}}\td(y,z)\\
        & \qquad = \frac{1}{2}\int_{\R^{2d}\setminus\{|z-y|<\varepsilon\}}\frac{(H^m(z)-H^m(y))^{p-1}(Q(z)-Q(y))}{|z-y|^{d+sp}}\td(y,z)
    \end{align*}
    which converges to $\left((-\Delta_p)^{s}_{\varepsilon}H^m\right)Q$ as $\varepsilon \rightarrow 0$.
    Taking $q$ to approximate $\sign_0^{+}$ and integrating the previous difference estimate, we obtain,
    \begin{align*}
        \int_{\R^d}(u(t)-H(t))^{+}\dmu \le \int_{\R^d}(u(0)-H(0))^{+}\dmu.
    \end{align*}
    
    Similarly, comparing $-u$ and $-H$ and noting that $-H$ will be a
    sub-solution to~\eqref{eq:1}, we have that
    \begin{displaymath}
      |u(t)| \le H(t)\qquad\text{a.e. on
    $\R^d$ for $t \ge 0$ provided $|u_0| \le H(0)$.}
    \end{displaymath}
\end{proof}

\begin{proof}[Proof of Theorem~\ref{thm:globalbarrier}]
    Converting the barriers $G$ given by~\eqref{eq:barrierfunSuperlinear} for $p > p_1$,~\eqref{eq:barrierfunCritical} for $p = p_1$ and \eqref{eq:lowSublinearProfileBound} for $p_{m,c}< p < p_1$ back to $(x,t)$ coordinates gives us the barrier function
    \begin{displaymath}
    H(x,t) = (t+1)^{-d\beta}G(x(t+1)^{-\beta})
    \end{displaymath}
    defined by~\eqref{eq:barrierfunH}. Note that $H(x,0) = G(x)$.

    Hence in light of Lemma \ref{lem:comparisonSupersolution} and the previous pointwise estimates, we need only prove that such constants can be chosen for $G$ in each case so that $H$ bounds $|u|$ at $t = 0$ and the super-solution condition \eqref{eq:supersolEquation} is satisfied whenever $u > H$.

    First, in the case $p > p_1$, we want to satisfy the conditions of Lemma \ref{lem:nearRegionSuperlinear} and Lemma \ref{lem:farRegion} while also ensuring that $H$ bounds $u$ elsewhere.
    In particular, this requires~\eqref{eq:nearCondition},~\eqref{eq:farRegionInequalities} and $C_1 \ge C$ for some $C > 0$ with this last inequality coming from the pointwise bound in the intermediate region $R_1 \le |z| \le 2R_2$.
    For this we consider cases on the homogeneity of the operator $(-\Delta_p)^s \cdot^m$.

    If $m(p-1) > 1$, we require
    \begin{displaymath}
        C_1 \ge CR_1^{d+\frac{sp}{m(p-1)-1}},\quad 
        C_1\le CR_1^{d}R_2^{\frac{sp}{m(p-1)-1}},\quad
        C_1 \ge C
    \end{displaymath}
    and so we can take $R_1$ large enough to contain the support of $u_0$, $C_1$ large enough to satisfy the first and last inequalities while also ensuring that $A \ge \norm{u_0}_{\infty}$ and $R_2$ large enough for the remaining inequality.

    In the case $m(p-1) = 1$, we require
    \begin{displaymath}
        R_1^{sp} \le C,\quad 
        R_2^{sp} \ge C,\quad
        R_2^{sp} \ge C\log(R_2),\quad
        C_1 \ge C.
    \end{displaymath}
    So we can take $R_1$ sufficiently small and $R_2$ sufficiently large to satisfy the first three inequalities.
    Then by taking $C_1$ large enough, we can ensure that $H(0)$ bounds $|u_0|$.

    In the case $m(p-1) < 1$, we require
    \begin{displaymath}
        C_1 R_1^{\frac{sp}{1-m(p-1)}-d} \le C,\quad
        C_1 R_2^{\frac{sp}{1-m(p-1)}-d}\ge C,\quad
        C_1 \ge C.
    \end{displaymath}
    Note that $p_c< p < 1+\frac{1}{m}$ implies that $\frac{sp}{1-m(p-1)} > d$.
    So fix $R_2 > 0$ such that $B_{R_2}$ contains the support of $u_0$. Then take $C_1$ large enough to satisfy the last two inequalities and ensure that $H$ bounds $u_0$ at $t = 0$. Finally take $R_1$ small enough to satisfy the first inequality.

    In the case $p = p_1$, we require~\eqref{eq:nearCondition},~\eqref{eq:farRegionInequalitiesCritical} and $C_1 \ge C$ for some $C > 0$. As in the previous case, we can take $R_2$ large, $C_1$ large and $R_1$ small to satisfy all conditions.

    Finally, in the case $p_{m,c}< p < p_1$, we require~\eqref{eq:farRegionInequalitiesSublinear} and, by applying Lemma~\ref{lem:nearRegionSuperlinear} with $R_1 = R_2$,
    \begin{displaymath}
        C_2 R_2^{\frac{sp}{1-m(p-1)}-d} \le C.
    \end{displaymath}
    so we can satisfy all conditions by taking $R_2$ small and $C_2$ large.    
	\end{proof}

\section{Existence of a Barenblatt solution}
\label{section:existence-of-Barenblatt}

In this section, we aim to prove the existence of a Barenblatt
solution which tends to a Dirac delta as $t \rightarrow 0+$. Here we
use a rescaling technique from~\cite{MR2286292}, see
also~\cite{MR4114983} where this is proven for the fractional
$p$-Laplacian.

\begin{theorem}
    \label{thm:BarenblattExistence1}
    Let $d \ge 1$, $0 < s < 1$, $m \ge 1$ and $p > p_c$ such that $m(p-1) \neq 1$.
    Then for any mass $M > 0$ there exists a strong distributional solution $\Gamma$ to~\eqref{eq:1} in $L^1$ which is nonnegative, radially symmetric, decreasing radially in space with decay given by  $g(|x|)$ and decays in time with order $t^{-\frac{1}{m(p-1)-\frac{p}{q_s}}}$ uniformly in $x$.
\end{theorem}

\begin{proof}
	For a given smooth, positive, radially symmetric initial datum $u_{0}\in L^{1\cap \infty}$ with compact support contained in the open unit ball $B_{1}(0)$ and having mass $\int_{\R^{d}}u_{0}\,\dx=1$, let $u$ be the corresponding positive strong distributional solution of~\eqref{eq:1} provided by Theorem~\ref{thm:regularisation-effect}.
    By applying the scaling transformation~\eqref{eq:transform1} to $u$, one obtains that for every integer $k \ge 1$,
	\begin{equation}
	\label{eq:scalingTrans}
	    u_k(x,t): = k^d u(kx, k^{d(m(p-1)-1)+sp}t),\quad x\in \R^{d},\;t\ge0,
	\end{equation}
	is a strong distributional solution of the doubly nonlinear nonlocal diffusion equation~\eqref{eq:heat-eq} with 
	initial datum $u_{k,0}:=k^{d}\,u_{0}(k\cdot)$ representing a nascent $\delta$-function in the sense that $u_{k,0}$ converges to the Dirac delta $\delta_{0}$ function in the sense of distributions as $k\to \infty$. Therefore, it remains to study the existence and properties of the limit function
	\begin{displaymath}
	    \lim_{k\to\infty} u_{k}(x,t)\qquad\text{for every $x\in \R^d$ and $t>0$.}
	\end{displaymath}
    In particular, we will take the limit of a subsequence of $(u_k)_{k \ge 1}$ to define the Barenblatt solution $\Gamma(x,t)$.
	By Proposition~\ref{prop:growthEstimate}, 
	\begin{equation}
	\label{eq:non-inc-uk-in-L1}
	    \norm{u_k(t)}_1 \le \norm{u_{k,0}}_{1}=1
	\end{equation}
	for every $t\ge 0$ and every $k\ge 1$. From the
    $L^1-L^{\infty}$ regularization estimate~\eqref{eq:aprioriL1Linf}, 
    we have that
    \begin{equation}
      \label{eq:uk-Linfty-time}
       \sup_{k\ge 1}\;\norm{u_{k}(t)}_{\infty}\le C\,
       t^{-\alpha}\qquad\text{for every $t>0$,}
    \end{equation}   
    where $\alpha > 0$ is defined as in Proposition~\ref{prop:L1LinfRegularity} and so \eqref{eq:non-inc-uk-in-L1} yields that
    \begin{equation}
        \label{eq:uniform-Lp-estimate-on-ukm}
         \sup_{k\ge 1}\,\norm{u_k^m(t)}_p 
         \le C^{(m-1)}\,t^{-(m-1)\alpha}\qquad
         \text{for every $t>0$.}
    \end{equation}
    Further, applying \eqref{eq:aprioriTime}
    and~\eqref{eq:energyEstimate} to 
    $(u_k)_{k\ge 1}$ gives that
      \begin{align}
      \label{eq:dtuk-L1-time}
       \sup_{k\ge 1}\;\norm{\partial_{t}u_{k}(t)}_{1} &\le 
       \frac{2}{|m(p-1)-1|}\frac{1}{t}\qquad\text{for every $t>0$,}\\ 
     \intertext{and}  
       \label{eq:wspuk-L1-time}
        \sup_{k\ge 1}\,[u^{m}_{k}(t))]^{p}_{s,p} &\le C\,t^{-(1+m\alpha)}
        \hspace{2.35cm}
        \text{for every $t>0$.}
    \end{align}
    Therefore, for every $\delta>0$,  
    $(\frac{\partial u_k}{\partial t})_{k\ge 1}$ is bounded in
    $L^{\infty}(\delta,\infty;L^{1})$. 
    In particular, \eqref{eq:uniform-Lp-estimate-on-ukm} and \eqref{eq:wspuk-L1-time} yield that the sequence $(u^m_k)_{k\ge 1}$ is bounded in $L^{\infty}(\delta,\infty;W^{s,p})$. 
    Since the previous estimates \eqref{eq:non-inc-uk-in-L1}-\eqref{eq:wspuk-L1-time} remain valid on any compact subset $K$ of $\R^d$, and $W^{s,p}(K)$ is compactly embedded into $L^1(K)$ by the Rellich-Kondrachov theorem for fractional Sobolev spaces (see~\cite[Theorem~2.1]{MR4108436}), it follows from Simon's compactness result \cite[Theorem 1]{MR916688} that $(u_k^m)_{k\ge 1}$ is relatively compact in $C([t_1,t_2];L^1(K))$ for every $0<t_1<t_2$. 
    
    Now, for every integer $n\ge 1$, let $I_{n}=[1/n,n]$ and $K_n
    =\{x\in\R^d : \abs{x}\le n\}$. Then $(I_{n})_{n\ge 1}$ is an increasing sequence of compact intervals, approximating the
    positive open real line $(0,\infty)$, and $(K_n)_{n\ge 1}$
    is an increasing sequence of compact subsets of $\R^d$
    approximating $\R^d$. Since both sequences $(I_{n})_{n\ge 1}$
    and $(K_n)_{n\ge 1}$ are countable, a standard diagonal argument
    yields the existence of a subsequence $(u_{\varphi(k)})_{k\ge
    1}$ of $(u_k)_{k\ge 1}$ and an element
    $\Gamma\in C((0,\infty);L^1_{loc})$ such that
    \begin{equation}
        \label{eq:cont-limit-in-L1K}
        \lim_{k\to\infty}u_{\varphi(k)}=\Gamma\qquad
        \text{in $C([t_1,t_2];L^1(K))$}
    \end{equation}
    for every $0<t_1<t_2$ and every compact subset $K$ of $\R^d$. 
    Further, since each $u_{\varphi(k)}$ is given by \eqref{eq:scalingTrans}, we can apply the global bound for $u$ given by~\eqref{eq:barrierfunH} and~\eqref{eq:barrierfunHSubcritical}, so that for all $R > R_2$ and $t \ge 0$, we have
    \begin{align*}
        \norm{u_{\varphi(k)}(t)}_{L^1(\{\abs{x}\ge R\})}
        &\le C_2(t+\varphi(k)^{-\frac{1}{\beta}})^{-d\beta}\int_{\{\abs{x}\ge R\}} g\left(|x|(t+\varphi(k)^{-\frac{1}{\beta}})^{-\beta}\right)\,\dx\\
        &= C_2\int_{\{\abs{z}\ge R(t+\phi(k)^{-\frac{1}{\beta}})^{-\beta}\}}g(|z|)\,\dz.
   \end{align*}
    
    Since $g(|z|)$ is integrable in $\R^d$ in each case, for any given $\varepsilon>0$, there is an $R>0$ such that for every $0<t_1<t_2$,
    \begin{displaymath}
        \sup_{k\ge 1} \sup_{t\in [t_1,t_2]}\norm{u_{\varphi(k)}(t)}_{L^1(\{\abs{x}\ge R\})}<\varepsilon.
    \end{displaymath}
    Combining this with \eqref{eq:cont-limit-in-L1K} gives
    \begin{displaymath}
        \lim_{k\to\infty}u_{\varphi(k)}=\Gamma\qquad
        \text{in $C([t_1,t_2];L^1)$}
    \end{displaymath}
    for every $0<t_1<t_2$, and by \eqref{eq:aprioriL1Linf},
    \begin{equation}
        \label{eq:cont-limit-in-L1}
        \lim_{k\to\infty}u_{\varphi(k)}=\Gamma\qquad
        \text{in $C([t_1,t_2];L^q)$}
    \end{equation}
    for every $1\le q<\infty$.
   % By Corollary \ref{cor:aleksandrov2}, 
   % $u(\cdot,k^{d(m(p-1)-1)+sp}\,t)$ is radially decreasing on
   % the exterior ball $\{\abs{y}\ge k\,R\}$ for every $t>0$. 
    %Thus and, one has that
    %\begin{displaymath}
    %    \sup_{k\ge 1}\norm{u_k(t)}_{L^1(\{\abs{x}\ge R\})}
    %     \le \sup_{k\ge 1}\norm{u(k^{d(m(p-1)-1)+sp}
    %\,t)}_{L^1(\{\abs{y}\ge R\})}.
    %\end{displaymath}
    %In particular, \eqref{eq:non-inc-uk-in-L1} implies that %$U(t)\in L^1$ satisfying
    %\begin{displaymath}
    %    \norm{U(t)}_1\le 1\qquad \text{for every $t>0$.}
    %\end{displaymath}
    Next, let $0<t_1<t_2$ and $K$ be a compact subset of $\R^d$. Then by \eqref{eq:dtuk-L1-time} and \eqref{eq:cont-limit-in-L1},
    one sees that
	\begin{align*}
		\norm{\Gamma(t_{2})-\Gamma(t_1)}_{L^1(K)}& \le \lim\limits_{k\rightarrow\infty}\norm{u_k(t_2)-u_k(t_1)}_{L^1(K)}\\
		& \le \lim\limits_{k\rightarrow\infty}\int_{t_1}^{t_2}\lnorm{\partial_t u_k(t)}_{L^1(K)}\,\dt\\
        & \le \lim\limits_{k\rightarrow\infty}\int_{t_1}^{t_2}\lnorm{\partial_t u_k(t)}_{1}\,\dt\\
		& \le 2\frac{\log t_2-\log t_1}{\abs{m(p-1)-1}}.
	\end{align*}
    Hence, applying an increasing sequence $(K_n)_{n\ge 1}$ of
    compact subsets of $\R^d$ to the preceding inequality and
    subsequently sending $n\to\infty$ yields that
    \begin{equation}
    \label{eq:absolute-continuity-in-L1}
        \norm{\Gamma(t_{2})-\Gamma(t_1)}_{1}\le  2\,
        \frac{\log t_2-\log t_1}{\abs{m(p-1)-1}}
    \end{equation}
    for every $0<t_1<t_2$. Therefore $\Gamma\in C((0,\infty);L^1)$ and is locally absolutely continuous with values in $L^1$.
    Moreover, by \eqref{eq:absolute-continuity-in-L1}, $\Gamma$ is in $W^{1,\infty}(\delta,T;L^1)$ satisfying
    \begin{displaymath}
    \lim_{h\to0+}\lnorm{\frac{\Gamma(t+h)-\Gamma(t)}{h}}_{1}\le  2\,
       \lim_{h\to0+} \frac{\log( t+h)-\log t}{h\abs{m(p-1)-1}}
        =\frac{2}{t\,\abs{m(p-1)-1}}
    \end{displaymath}
    
    Next, since $u_{\varphi(k)}$ is a strong distributional solution of \eqref{eq:heat-eq}, we may multiply~\eqref{eq:heat-eq} by $u_{\varphi(k)}^m$ and subsequently integrate over $(t_1,t_2)$ for given $0<t_1<t_2$. Then, one obtains that
    \begin{equation}
        \label{eq:integrated-approx-eq}
        \tfrac{1}{m+1}\norm{u_{\varphi(k)}(t_2)}_{m+1}^{m+1}+
        \int_{t_1}^{t_2}[u^{m}_{\varphi(k)}(t))]^{p}_{s,p}\,\dt=
        \tfrac{1}{m+1}\norm{u_{\varphi(k)}(t_1)}_{m+1}^{m+1}.
    \end{equation}
    Integrating~\eqref{eq:wspuk-L1-time} over $(t_1,t_2)$, we also have that
    \begin{equation}
        \label{eq:uniform-bdd-Wsp}
        %\tfrac{1}{m+1}\norm{u_k(t_2)}_{m+1}^{m+1}+
        \int_{t_1}^{t_2}[u^{m}_{\varphi(k)}(t))]^{p}_{s,p}\,\dt
        \le \tfrac{C}{m\alpha}t_{1}^{-m\alpha}.
    \end{equation}
    Thanks to the two estimates \eqref{eq:uniform-Lp-estimate-on-ukm} and \eqref{eq:uniform-bdd-Wsp}, we have that for every $0<t_1<t_2$, $(u^m_{\varphi(k)})_{k\ge 1}$ is bounded in $L^p(t_1,t_2;W^{s,p})$.
    From here, we proceed similarly to the proof of Theorem \ref{thm:regularisation-effect} in order to show that $\Gamma$ is a strong distributional solution to~\eqref{eq:1}. Letting $0<t_1<t_2$, the space
    $L^{p}(t_1,t_2;W^{s,p})$ is reflexive, so
    by~\eqref{eq:cont-limit-in-L1} we can conclude that $\Gamma^{m}
    \in L^{p}(t_1,t_2;W^{s,p})$. After possibly passing
    to another subsequence of $(u^m_{\varphi(k)})_{k\ge 1}$, we have that
    \begin{equation}
        \label{eq:weak-lim-LpWsp-Barenblatt}
        \lim_{k\to\infty}u_{\varphi(k)}^m=
        \Gamma^m\qquad\text{weakly in $L^{p}(t_1,t_2;W^{s,p})$.}
    \end{equation}
    In particular, one has that the sequence $(\mathcal{A}_{s,p}
    (u_{\varphi(k)}^{m}))_{k\ge 1}$ of linear bounded functionals on $L^{p}
    (t_1,t_2;W^{s,p})$ given by
    \begin{align*}
        &\langle \mathcal{A}_{s,p}(u_{\varphi(k)}^{m}),
        \xi\rangle_{L^{p^{\mbox{}_{\prime}}}
        (t_1,t_2;W^{-s,p^{\mbox{}_{\prime}}}),L^{p}(t_1,t_2;W^{s,p})}\\
        &\quad 
        =\frac{1}{2}\int_{t_{1}}^{t_{2}}\int_{\R^{2d}}
        \frac{
        (u_{\varphi(k)}^{m}(x)-u_{\varphi(k)}^{m}(y))^{p-1}(\xi(x)-\xi(y))}{
        \abs{x-y}^{d+sp}}\,\td(x,y)\,\dt
    \end{align*}
    for every $\xi\in L^{p}(t_1,t_2;W^{s,p})$, is bounded in
    $L^{p^{\mbox{}_{\prime}}}(t_1,t_2;W^{-s,p^{\mbox{}_{\prime}}})$.
    Therefore, there is an $\chi\in L^{p^{\mbox{}_{\prime}}}(t_1,t_2;W^{-s,p^{\mbox{}_{\prime}}})$ such that after
    possibly passing to a subsequence, one has that 
    \begin{equation}
        \label{eq:limit-of-chiX-2}
        \lim_{n\to\infty}\mathcal{A}_{s,p}(u_{\varphi(k)}^{m})=\chi\qquad
        \text{weakly$\mbox{}^{\ast}$ in $L^{p^{\mbox{}_{\prime}}}
        (t_1,t_2;W^{-s,p^{\mbox{}_{\prime}}})$.}
    \end{equation}
    Further, since $u_{\varphi(k)}$ is a strong distributional solution of \eqref{eq:heat-eq}, it follows from
    \begin{equation}
        \label{eq:Asp-heat-eq}
        \partial_{t}u_{\varphi(k)}+\mathcal{A}_{s,p}(u_{\varphi(k)}^{m})=0
    \end{equation}
    that $(\partial_{t}u_{\varphi(k)})_{k\ge 1}$ is bounded in $L^{p^{\mbox{}_{\prime}}}(t_1,t_2;W^{-s,p^{\mbox{}_{\prime}}})$. 
    Thus and since $\Gamma\in L^1$, it follows from~\eqref{eq:cont-limit-in-L1} that $\partial_t \Gamma\in L^{p^{\mbox{}_{\prime}}}(t_1,t_2;W^{-s,p^{\mbox{}_{\prime}}})$ and, after possibly passing to another subsequence of $(u_{\varphi(k)})_{k\ge 1}$, that
    \begin{equation}
        \label{eq:limit-of-partial-t-u-varphi-k}
        \lim_{k\to\infty}\partial_{t}u_{\varphi(k)} =\partial_{t}\Gamma\qquad
        \text{weakly$\mbox{}^{\ast}$ in $L^{p^{\mbox{}_{\prime}}}(t_1,t_2;W^{-s,p^{\mbox{}_{\prime}}})$.}
    \end{equation}
    Hence, multiplying~\eqref{eq:Asp-heat-eq} by a test function $\xi \in C_c^{\infty}((t_1,t_2)\times \R^d)$ and sending $k\to \infty$ yields that
    \begin{align*}
        \braket{\partial_t \Gamma(t)+X^*,\xi}_{L^{p^{\mbox{}_{\prime}}}(t_1,t_2;W^{-s,p^{\mbox{}_{\prime}}}),L^p(t_1,t_2;W^{s,p})} = 0,
    \end{align*}
    which yields that 
    \begin{equation}
      \label{eq:evol-eq-Gamma}
        \partial_{t}\Gamma(t)+\chi=0\qquad\text{in $W^{-s,p^{\mbox{}_{\prime}}}$ for a.e. $t\in (t_1,t_2)$.}
    \end{equation}
    Now, we are ready to prove that $\Gamma$ is a distributional solution of~\eqref{eq:1} on $(t_1,t_2)$. To do this, it remains to show that $\chi=\mathcal{A}_{s,p}(\Gamma^{m})$, where $\mathcal{A}_{s,p}$ is the lifted operator $\mathcal{A}_{s,p} : 
    L^{p}(t_1,t_2;W^{s,p})\to L^{p^{\mbox{}_{\prime}}}(t_1,t_2;W^{-s,p^{\mbox{}_{\prime}}})$ given by 
    \begin{align*}
        &\langle \mathcal{A}^{s}_{p}(v),\xi\rangle_{L^{p^{\mbox{}_{\prime}}}(t_1,t_2;W^{-s,p^{\mbox{}_{\prime}}}),L^{p}(t_1,t_2;W^{s,p})}\\
        &\qquad 
        =\frac{1}{2}\int_{t_{1}}^{t_{2}}\int_{\R^{2d}}
        \frac{\abs{v(x)-v(y)}^{p-2}
        (v(x)-v(y))(\xi(x)-\xi(y))}{\abs{x-y}^{d+sp}}\,\td(x,y)\,\dt
    \end{align*}
    for every $v$, $\xi\in L^{p}(t_1,t_2;W^{s,p})$. 
    %For this, we apply another time the monotonicity trick to
    %the lifted operator $\mathcal{A}^{s}_{p}$. First, we note 
    For this, note that multiplying~\eqref{eq:evol-eq-Gamma} by $\Gamma^{m}$ yields that
    \begin{displaymath}
        \tfrac{1}{m+1}\norm{\Gamma}_{m+1}^{m+1}\Big\vert_{t_1}^{t_2}+\langle \chi,\Gamma^{m}\rangle_{L^{p^{\mbox{}_{\prime}}}(t_1,t_2;W^{-s,p^{\mbox{}_{\prime}}}),L^{p}(t_1,t_2;W^{s,p})}=0.
    \end{displaymath}
    Thus~\eqref{eq:cont-limit-in-L1} and \eqref{eq:integrated-approx-eq} yield that
    \begin{equation}
        \label{eq:lim-lift-uvarphi-k-m}
        \lim_{n\to\infty}\int_{t_1}^{t_2}[u_{\varphi(k)}^{m}(t)]_{p,s}^{s}\,\dt=
        \langle \chi,\Gamma^{m}\rangle_{L^{p^{\mbox{}_{\prime}}}(t_1,t_2;W^{-s,p^{\mbox{}_{\prime}}}),L^{p}(t_1,t_2;W^{s,p})}.
    \end{equation}
    By the monotonicity of $\mathcal{A}_{s,p}$, one has that
    \begin{align*}
      0\le &\langle \mathcal{A}^{s}_{p}(u_{\varphi(k)}^{m})-\mathcal{A}_{s,p}(\xi),u_{\varphi(k)}^{m}-\xi\rangle_{L^{p^{\mbox{}_{\prime}}}(t_1,t_2;W^{-s,p^{\mbox{}_{\prime}}}),L^{p}(t_1,t_2;W^{s,p})}\\
        =&\int_{t_1}^{t_2}[u_{\varphi(k)}^{m}(t)]_{s,p}^{p}\,\dt- \langle \mathcal{A}^{s}_{p}(u_{\varphi(k)}^{m}),\xi
        \rangle_{L^{p^{\mbox{}_{\prime}}}(t_1,t_2;W^{-s,p^{\mbox{}_{\prime}}}),L^{p}(t_1,t_2;W^{s,p})}\\
        &\hspace{2cm} -
        \langle\mathcal{A}_{s,p}(\xi),u_{\varphi(k)}^{m}-\xi\rangle_{L^{p^{\mbox{}_{\prime}}}(t_1,t_2;W^{-s,p^{\mbox{}_{\prime}}}),L^{p}(t_1,t_2;W^{s,p})}
    \end{align*}
    for every $k\ge 1$. Thus, sending $k\to\infty$ in the last inequality and by using~\eqref{eq:weak-lim-LpWsp-Barenblatt}, \eqref{eq:limit-of-chiX-2}, and \eqref{eq:lim-lift-uvarphi-k-m}, one obtains that
    \begin{displaymath}
         0\le \langle \chi-\mathcal{A}_{s,p}(\xi),u^{m}-\xi\rangle_{L^{p^{\mbox{}_{\prime}}}(t_1,t_2;W^{-s,p^{\mbox{}_{\prime}}}),L^{p}(t_1,t_2;W^{s,p})}
    \end{displaymath}
    for every $\xi\in L^{p}(t_1,t_2;W^{s,p})$. Now, by proceeding as in the proof of Proposition~\ref{prop:approxStrongSols}, we can choose $\xi=u^{m}-\mu\,\zeta$, taking $\mu \rightarrow 0+$, to conclude that $\chi=\mathcal{A}_{s,p}(\Gamma^{m})$.

    Thus $\Gamma$ is a distributional solution of \eqref{eq:heat-eq} and it remains to show that $\Gamma$ is differentiable with values in $L^1_{loc}$ at almost everywhere $t\in (0,\infty)$. Since the argument is exactly the same as the one given in the proof of Theorem~\ref{thm:well-posedness}, we omit this part.
    
    %$u_{\varphi(k)}^m(t)\rightharpoonup U^m(t)$ weakly in $W^{s,p}$. Hence and  %\eqref{eq:wspuk-L1-time} imply that
    %\begin{displaymath}
    %   [U^{m}(t))]^{p}_{s,p} \le \liminf_{k\to\infty}[u^{m}_{k}
    %(t))]^{p}_{s,p} \le C\,t^{-\delta}.
    %\end{displaymath}
	
	We can estimate $\Gamma$ by applying the global bounds~\eqref{eq:barrierfunH} and~\eqref{eq:barrierfunHSubcritical} to $u$ so that there exists $C_2 > 0$ and $R_2 > 0$ such that
	\begin{displaymath}
	    u(x,t) \le C_2(t+1)^{-d\beta}g(|x|(t+1)^{-\beta}),
	\end{displaymath}
 for a.e. $x \in \R^d$ and $t \ge 0$ satisfying $|x| > (t+1)^{\beta}R_2$. Then for $k \ge 1$,
	\begin{displaymath}
	    u_k(x,t) \le C_2 (t+k^{-\frac{1}{\beta}})^{-d\beta}g\left(|x|(t+k^{-\frac{1}{\beta}})^{-\beta}\right)
	\end{displaymath}
	for all $|x| > (t+k^{-\frac{1}{\beta}})^{\beta}R_2$. Hence there exists $C > 0$ such that
	\begin{equation}
    \label{eq:GammaBound}
	    \Gamma(x,t) \le Ct^{-d\beta}g(|x|t^{-\beta})
	\end{equation}
	for all $|x| > Ct^{\beta}$. So $\Gamma$ tends to a Dirac delta as $t \rightarrow 0+$.
 
    It remains to prove the radial symmetry properties of $\Gamma$. Since $u_{\phi(k)}^m$ converged to $\Gamma^m$ weakly in $L^p(t_1,t_2;W^{s,p})$, we have a subsequence, relabelled as $u_k$, which converges pointwise almost everywhere in $[t_1,t_2]\times \R^d$ for $0 < t_1 < t_2$. 
    Then applying Corollary~\ref{cor:aleksandrov2}, $u_k$ is radially symmetric and radially decreasing so that the same holds for $\Gamma$ almost everywhere in $(0,\infty)\times\R^d$.
    Moreover,~\eqref{eq:GammaBound} gives the decay of $\Gamma$ in space and~\eqref{eq:uk-Linfty-time} gives the decay in time.
    \end{proof}

    \begin{proposition}
        The strong distributional solution $\Gamma$ given by Theorem~\ref{thm:BarenblattExistence1} is self-similar with the form~\eqref{eq:selfSimilarForm} where $F:\R^d \rightarrow \R$ satisfies the properties of Theorem~\ref{thm:BarenblattSolutions}.
    \end{proposition}
    
    \begin{proof}
    We apply a rescaling to self-similar variables,
    \begin{equation}
    \label{eq:selfSimilarRescaling}
        V(y,\tau) = t^{d\beta}\Gamma(x,t)\quad \text{for } t > 0 \text{ and } x\in \R^d
    \end{equation}
    where $y = xt^{-\beta}$ and $\tau = \log t$. 
    Then since $\Gamma \in W^{1,\infty}(\delta,\infty;L^1)\cap L^{\infty}(\delta,\infty;W^{s,p})$ for all $\delta > 0$, we can apply those estimates such that $V \in W^{1,\infty}(\delta,\infty;L^1)\cap L^{\infty}(\delta,\infty;W^{s,p})$ for all $\delta > 0$. 
    Moreover, since $\Gamma$ is a strong distributional solution, $V$ is a strong distributional solution to
    \begin{displaymath}
    \partial_{\tau}V-\beta\nabla\cdot(yV)+(-\Delta_p)^s V = 0
    \end{displaymath}
    for $\tau \in \R$ and $y \in \R^d$.
    In particular, since $\partial_{t}\Gamma \in L^{\infty}(\delta,\infty;L^1)$, $\partial_{\tau}V \in L^{\infty}(\delta,\infty;L^1)$ and
    \begin{displaymath}
        \partial_{t}\Gamma = t^{-d\beta-1}\left(\partial_{\tau}V-\beta\nabla\cdot (yV)\right)
    \end{displaymath}
    by the rescaling~\eqref{eq:selfSimilarRescaling}, we have that $\nabla\cdot(yV) \in L^{\infty}(\delta,\infty;L^1)$ for all $\delta > 0$.

    We have the following lemma from \cite{MR4114983} wherein $u_k$ is defined by~\eqref{eq:scalingTrans} and $u_1$ is the positive strong distributional solution corresponding to a given radially symmetric initial datum $u_0 \in L^{1\cap\infty}$ with compact support contained in the unit ball and having mass 1 as in the proof of Theorem~\ref{thm:BarenblattExistence1}.
    
    \begin{lemma}[{\cite[Lemma 6.1]{MR4114983}}]
    \label{lem:rescalingLem}
        If $v_1$ is the rescaled function from $u_1$ and $v_k$ from $u_k$ according to~\eqref{eq:selfSimilarRescaling}, then
        \begin{displaymath}
        v_k(y,\tau) = v_1(y,\tau+\log(k)).
        \end{displaymath}
    \end{lemma}
    
    Since (relabelling by an appropriate subsequence) $u_k$ converges pointwise almost everywhere to $\Gamma$ in $\R^d\times(0,\infty)$, $v_k(y,\tau)$ also converges pointwise almost everywhere to $V(y,\tau)$ in $\R^d\times(-\infty,\infty)$.
    By Lemma \ref{lem:rescalingLem}, we know that $v_k(y,\tau) = v_1(y,\tau+\log(k))$ so that, for $h > 0$, pointwise we have that
    \begin{align*}
    V(y,\tau+h)& = \lim\limits_{k\rightarrow\infty}v_{k}(y,\tau+h)\\
    & = \lim\limits_{k\rightarrow\infty}v_{k+e^h}(y,\tau)\\
    & = V(y,\tau)
    \end{align*}
    almost everywhere in $\R^d\times(-\infty,\infty)$. Hence $V$ is independent of $\tau$ and $V_{\tau} = 0$ almost everywhere.
    Hence $F = V$ is an appropriate profile function. Moreover, carrying the radial symmetry properties of $U$ to $F$ gives us the properties of Theorem~\ref{thm:BarenblattSolutions}.
    \end{proof}

\section{\texorpdfstring{Mass conservation}{}}
\label{sec:massConservation}

In this section we prove a stronger version of the standard $L^1$ decay for accretive operators in $L^1$. In particular, we prove the following. 
Recall the definition of $p_1$ given by \eqref{eq:globalBarrierCondition}.

\begin{theorem}
    Suppose $0 < s < 1$, $m > s$ and $p \ge p_1$. Let $u$ be a mild solution to \eqref{eq:1} with $u_0 \in L^1$. 
    Then
    \begin{equation}
    \label{eq:massConservation}
        \int_{\R^d}u(x,t_2)\dx = \int_{\R^d}u(x,t_1)\dx
    \end{equation}
    for all $0 < t_1 < t_2 < \infty$.
\end{theorem}

\begin{proof}
    By density, we approximate $u_0 \in L^1$ by initial data in $L^{1\cap\infty}$ with compact support, giving strong distributional solutions to \eqref{eq:1}. 
    Proving~\eqref{eq:massConservation} for such initial data, we then have the result for initial data in $L^1$ by the comparison estimate in $L^1$ given by \eqref{eq:comparisonPrinciple}.
    
    The main idea of the proof is to multiply \eqref{eq:1} by a sequence of positive test functions $(\phi_n)_{n\ge 1} \subset C_c^{\infty}(\R^d)$ such that $\phi_n \rightarrow \mathds{1}_{\R^d}$ pointwise in $\R^d$ as $n \rightarrow \infty$.
    In particular,  consider $\phi$ such that $\phi = 1$ for $|x| \le 2$ and $\phi = 0$ for $|x| \ge 3$. Then set $\phi_n:= \phi(x/n)$.
    Subsequently, integrating over $(t_1, t_2)$ for $0 < t_1 < t_2$, we have
    \begin{equation}
    \label{eq:massConservationMainEstInter}
    \begin{split}
        &\left|\int_{\R^d}\left(u(x,t_2)-u(x,t_1)\right)\phi_n(x) \dx\right|\\  
        &\quad\le  \int_{t_1}^{t_2} \int_{\R^{2d}}\frac{\left|(u^m(x,t)-u^m(y,t))^{p-1}(\phi_n(x)-\phi_n(y))\right|}{|x-y|^{d+sp}}\td(x,y)\dt.
    \end{split}
    \end{equation}
    Note that
    \begin{equation}
    \label{eq:phiNEstimate}
        [\phi_n]_{s,p} = n^{\frac{d-sp}{p}}[\phi]_{s,p}.
    \end{equation}
    So when $d < sp$, we estimate \eqref{eq:massConservationMainEstInter} by
    \begin{displaymath}
        \int_{t_1}^{t_2} [u^m]_{s,p}^{p-1} \dt\, (t_2-t_1) [\phi_n]_{s,p}.
    \end{displaymath}
    Then we apply Lemma~\ref{lem:generalEnergyEst} with $r = m$,
    \begin{align*}
        \int_{t_1}^{t_2}[u^m]_{s,p}^{p}\dt \le \frac{1}{m+1}\left(\norm{u(t_1)}_{m+1}^{m+1}-\norm{u(t_2)}_{m+1}^{m+1}\right)
    \end{align*}
    implying the boundedness of $\int_{t_1}^{t_2} [u^m]_{s,p}^{p-1} \dt$ by Jensen's inequality and the standard decay estimate.
    So taking $n \rightarrow \infty$ in \eqref{eq:massConservationMainEstInter} gives \eqref{eq:massConservation}.
    
    When $d \ge sp$, we estimate 
    \eqref{eq:massConservationMainEstInter} more precisely by considering regions of $\R^{2d}$ depending on $n$, in particular so that we can use the cancellation of $\phi_n$ terms when $|x| \le 2n$ and $|y| \le 2n$.    
    We first consider the region 
    \begin{displaymath}
        A_n = \set{x, y \in \R^d | |x| \ge n, |y| \ge n}.
    \end{displaymath}
    We estimate \eqref{eq:massConservationMainEstInter} on $(t_1,t_2)\times A_n$ by
    \begin{equation}
    \label{eq:massConservationAnHolder}
        \left(\int_{t_1}^{t_2}\int_{A_n}\frac{|u^m(x,t)-u^m(y,t)|^p}{|x-y|^{d+sp}}\td(x,y)\dt\right)^{\frac{p-1}{p}}\left(\int_{t_1}^{t_2}[\phi_n]_{s,p}^p\dt\right)^{\frac{1}{p}}.
    \end{equation}
    Here we use the key inequality, for $a$, $b > 0$ and $0 < \varepsilon < 1$,
    \begin{displaymath}
        |a-b|\left(|a|+|b|\right)^{-1+\varepsilon} \le C_{\varepsilon}|a^{\varepsilon}-b^{\varepsilon}|
    \end{displaymath}
    where $C_{\varepsilon} > 0$ depends only on $\varepsilon$. 
    Applying this to estimate $|a-b|^{p}$,
    \begin{align*}
        |a-b|^{p}& \le C_{\varepsilon}|a-b|^{p-1}(|a|+|b|)^{1-\varepsilon}|a^{\varepsilon}-b^{\varepsilon}|
    \end{align*}
    We note that when $a = u^m(x,t)$ and $b = u^m(y,t)$ with $(x,y) \in A_n$, we can estimate
    \begin{align*}
        (|u^m(x,t)|+|u^m(y,t)|)^{1-\varepsilon} \le 2\norm{u}_{L^{\infty}(t_1,t_2;L^{\infty}(\R^d\setminus B_n))}^{m(1-\varepsilon)}
    \end{align*}
    for all $t_1 \le t \le t_2$ where $B_n$ is the ball with radius $n$.
    By the global estimate of Lemma~\ref{lem:intermediateBound}, we have that
    \begin{align*}
        \norm{u}_{L^{\infty}(t_1,t_2;L^{\infty}(\R^d \setminus B_n))} \le Cn^{-d}
    \end{align*}
    where $C > 0$ depends on $d$ and $u_0$.
    Hence
    \begin{align*}
        &|u^m(x,t)-u^m(y,t)|^{p}\\
        &\quad\le C_{\varepsilon}n^{-d}|u^m(x,t)-u^m(y,t)|^{p-1}|u^{m\varepsilon}(x,t)-u^{m\varepsilon}(y,t)|,
    \end{align*}
    relabelling $C_{\varepsilon} > 0$ which now depends on $\varepsilon$, $d$ and $u_0$.
    Applying Lemma~\ref{lem:generalEnergyEst} with $q = m\varepsilon$,
    \begin{align*}
        &\int_{t_1}^{t_2}\int_{\R^{2d}}\frac{|u^m(x,t)-u^m(y,t)|^{p-1}|u^{m\varepsilon}(x,t)-u^{m\varepsilon}(y,t)|}{|x-y|^{d+sp}} \td(x,y)\dt\\
        &\quad\le \norm{u(t_1)}_{m\varepsilon+1}^{m\varepsilon+1}.
    \end{align*}
    Together with \eqref{eq:phiNEstimate}, we then estimate \eqref{eq:massConservationAnHolder} by
    \begin{displaymath}
        C_{\varepsilon}n^{-dm(1-\varepsilon)\frac{p-1}{p}}n^{\frac{d-sp}{p}}(t_2-t_1)^{\frac{1}{p}}\left(\norm{u(t_1)}_{m\varepsilon+1}^{m\varepsilon+1}\right)^{\frac{p-1}{p}}[\phi]_{s,p}.
    \end{displaymath}
    
    So taking $n \rightarrow \infty$, we have convergence to zero whenever 
    \begin{displaymath}
        (d+sp)m(1-\varepsilon)\frac{p-1}{p} > \frac{d-sp}{p}.
    \end{displaymath}
    In particular, we can choose a suitably small $\varepsilon > 0$ whenever
    \begin{displaymath}
        (d+sp)m(p-1) > d-sp
    \end{displaymath}
    which is satisfied since $p \ge p_1$.
    We look at the remaining regions. Whenever $|x|,|y| \le 2n$, $\phi_n(x)-\phi_n(y) = 0$. 
    So it remains to consider $D_n := \set{(x,y)||x|\le n, |y| \ge 2n}$ and $E_n := \set{(x,y)||x|\ge 2n, |y| \le n}$ which follows by symmetry.
    
    In $D_n$, $|x-y| \ge n$, so
    \begin{align*}
        &\int_{D_n}\frac{|u^m(x,t)-u^m(y,t)|^{p-1}|\phi_n(x)-\phi_n(y)|}{|x-y|^{d+sp}}\td(x,y)\\
        &\quad\le n^{-d-sp}\int_{\set{|x| \le n}}\int_{\set{|y| \ge 2n}}|u^m(x,t)-u^m(y,t)|^{p-1}\dy\dx\\
        &\quad \le 2n^{-d-sp}\left(\int_{\set{|x| \le n}}|u(x,t)|^{m(p-1)}\dx+\int_{\set{|y| \ge 2n}}|u(y,t)|^{m(p-1)}\dy\right)\\
        &\quad \le 2n^{-d-sp}\left(Cn^d\norm{u_0}_{L^{\infty}}^{m(p-1)}+\int_{2n}^{\infty}r^{-d(m(p-1)-1)-1}\dr\right).
    \end{align*}
    We have
    \begin{displaymath}
        \int_{\set{|x| \le n}}|u(x,t)|^{m(p-1)}\dx \le Cn^d\norm{u_0}_{L^{\infty}}^{m(p-1)}.
    \end{displaymath}
    For the second term we apply Theorem~\ref{thm:globalbarrier} with $n$ sufficiently large to estimate $|u(y,t)|$ by $C_2 (t+1)^{sp\beta}|x|^{-d-sp}$.
    Then since $p > p_1$, we find that the integral converges.
    Hence the $y$ term is bounded for all $n$ sufficiently large and so taking $n \rightarrow \infty$, this term also converges to zero.

    Since $\phi_n$ converges to $\mathds{1}$ pointwise in $\R^d$, \eqref{eq:massConservation} follows by taking $n \rightarrow \infty$ in \eqref{eq:massConservationMainEstInter}.
\end{proof}

\section{\texorpdfstring{Refined dissipation in $L^1$}{}}
\label{sec:dissipationInL1}

Due to the $L^1$ setting that we have introduced for $(-\Delta_p)^s \cdot^m$ and, in particular, the $T$-accretivity properties, we find dissipation of differences $(u_1-u_2)^{+}$ in $L^1$ for solutions to \eqref{eq:1}.
This follows the case of the evolution fractional $p$-Laplacian given by \cite{MR4114983}.
We note that this is a refinement of \eqref{eq:comparisonPrinciple}, the difference being that here we keep the full $(-\Delta_p)^s u^m$ terms when estimating rather than applying $T$-accretivity to estimate these by zero.
Due to this we cannot use the $T$-accretivity directly.

We note that this provides a strict dissipative effect for $(u_1 - u_2)^{+}$ whenever we have a non-zero set of times $t$ with $|\set{u_1 > u_2}| > 0$ and $|\set{u_1 \le u_2}| > 0$.

\begin{theorem}
\label{thm:L1dissipation}
    Suppose $m > 0$.
    We consider $u_1$, $u_2$ to be two strong distributional solutions to \eqref{eq:1} in $L^1$.
    Then for all $0 < t_1 < t_2 < T$,
    \begin{equation}
    \label{eq:dissipativeEffectEst}
    \begin{split}
         \norm{[u_1-u_2]^{+}(t_2)}_1& \le \norm{[u_1-u_2]^{+}(t_1)}_1 - \frac{1}{2}\int_{t_1}^{t_2} I(\tau)\dtau
    \end{split}
    \end{equation}
    where for all $\tau \in [0,T]$, omitting $\tau$ in $u_1(x,\tau)$ etc., $I(\tau)$ is given by
    \begin{equation}
    \label{eq:dissipationI}
        \int_{D(\tau)}\frac{\left|(u_1^m(x)-u_1^m(y))^{p-1}-(u_2^m(x)-u_2^m(y))^{p-1}\right|}{|x-y|^{d+sp}}\td(x,y)
    \end{equation}
    with the domain of integration 
    \begin{displaymath}
    \begin{split}
        D(\tau) = &\set{u_1(x,\tau) > u_2(x,\tau), \,u_1(y,\tau) \le u_2(y,\tau)}\\
        &\cup \set{u_1(x,\tau) \le u_2(x,\tau),\, u_1(y,\tau) > u_2(y,\tau)}.
    \end{split}
    \end{displaymath}
\end{theorem}

\begin{proof}
We multiply \eqref{eq:1} by $\sign_0^{+}(u_1-u_2)$ and integrate over $\R^d\times (t_1, t_2)$.
Then
\begin{align*}
    &\norm{[u_1-u_2]^{+}(t_1)}_1 - \norm{[u_1-u_2]^{+}(t_2)}_1\\
    &\quad= \int_{t_1}^{t_2}\int_{\R^d} \left((-\Delta_p)^s u_1^m - (-\Delta_p)^s u_2^m\right)\mathds{1}_{\set{u_1 > u_2}}\dx\dtau
\end{align*}
We evaluate the fractional $p$-Laplacian term as in \cite[Section 5]{MR4114983}, noting that $r \mapsto r^m$ strictly increasing implies that $\mathds{1}_{\set{u_1 > u_2}} = \mathds{1}_{\set{u_1^m > u_2^m}}$.
Here we omit the time dependence of $u_1(x,\tau)$ and $u_2(x,\tau)$ for brevity of notation.
We have, for $\tau \in (0,T)$,
\begin{equation}
\label{eq:keyDissipativeEstimate}
\begin{split}
    &\int_{\R^d} \left((-\Delta_p)^s u_1^m - (-\Delta_p)^s u_2^m\right)\mathds{1}_{u_1^m > u_2^m}\dx\\
    &\quad= \frac{1}{2}\int_{\R^{2d}}\frac{(u_1^m(x)-u_1^m(y))^{p-1}-(u_2^m(x)-u_2^m(y))^{p-1}}{|x-y|^{d+sp}}\\
    &\hspace{2.5cm}\times(\mathds{1}_{\set{u_1^m(x) > u_2^m(x)}}-\mathds{1}_{\set{u_1^m(y) > u_2^m(y)}})\td(x,y).
\end{split}
\end{equation}
Looking at the integrand, we notice that
\begin{equation}
\label{eq:dissipationSignTerm}
    \mathds{1}_{\set{u_1^m(x) > u_2^m(x)}}-
    \mathds{1}_{\set{u_1^m(y) > u_2^m(y)}}
\end{equation}
is $1$ if and only if both $u_1(x) > u_2(x)$ and $u_1(y) \le u_2(y)$. 
Similarly, \eqref{eq:dissipationSignTerm} is $-1$ if and only if $u_1(x) \le u_2(x)$ and $u_1(y) > u_(y)$.
Otherwise, this difference term is zero.
So in both cases where \eqref{eq:dissipationSignTerm} is non-zero, the integrand of \eqref{eq:keyDissipativeEstimate} will be positive and in fact we can write
\begin{displaymath}
\label{eq:keyDissipativeEstimate2}
\begin{split}
    &\int_{\R^d} \left((-\Delta_p)^s u_1^m - (-\Delta_p)^s u_2^m\right)\mathds{1}_{\set{u_1 > u_2}}\dx = I(\tau)
\end{split}
\end{displaymath}
with $I$ given by \eqref{eq:dissipationI}.

Importantly, this provides a dissipative effect for this difference $u_1 - u_2$ whenever $|\set{u_1 > u_2}| > 0$ and $|\set{u_1 \le u_2}| > 0$.
\end{proof}

The following important lemma allows us to compare solutions via a mass difference analysis and in particular, to obtain uniqueness of certain limits.

\begin{lemma}
\label{lem:L1dissipation}
    Suppose $m > 0$.
    We consider $u_1$, $u_2$ to be two strong distributional solutions to \eqref{eq:1} in $L^1$. Suppose $0 < t < \infty$ such that
    \begin{displaymath}
        \int_{\R^d}u_1(x,t)\dx = \int_{\R^d}u_2(x,t)\dx.
    \end{displaymath}
    %with initial data $u_{0,1}$ and $u_{0,2}$, respectively, such that $\norm{u_{0,1}}_1 = \norm{u_{0,2}}_1$.
    Then %for all $0 < t < \infty$,
    \begin{displaymath}
        \norm{\left(u_1(t)-u_2(t)\right)^{+}}_1
    \end{displaymath}
    is strictly decreasing unless $u_1(t) = u_2(t)$ a.e.~on $\R^d$.
\end{lemma}

\begin{proof}
    Since
    \begin{displaymath}
        \int_{\R^d} u_1(x,t)-u_2(x,t)\dx = 0,
    \end{displaymath}
    either $u_1(t) = u_2(t)$ a.e.~in $\R^d$ or
    \begin{displaymath}
        \int_{\R^d} (u_1(x,t)-u_2(x,t))^{+}\dx > 0
    \end{displaymath}
    and
    \begin{displaymath}
        \int_{\R^d} (u_2(x,t)-u_1(x,t))^{+}\dx > 0.
    \end{displaymath}
    Applying Theorem~\ref{thm:L1dissipation}, this second case implies that
    \begin{displaymath}
        \int_{\R^d} (u_1(x,t)-u_2(x,t))^{+}\dx
    \end{displaymath}
    must be decreasing.
\end{proof}

\section{Uniqueness of the Barenblatt solution}
\label{sec:barenblattUniqueness}

We obtain uniqueness of the self-similar profile using a method of mass difference analysis, applied to the evolution fractional $p$-Laplacian in \cite{MR4114983} (see also \cite{MR2286292}). 
This completes the proof of Theorem~\ref{thm:BarenblattSolutions}.

\begin{theorem}
    The Barenblatt solution given by Theorem~\ref{thm:BarenblattExistence1} is unique for each $M > 0$.
\end{theorem}

\begin{proof}
    Let $\Gamma_1$ and $\Gamma_2$ be two Barenblatt solutions to~\eqref{eq:heat-eq} with profile functions $F_1$ and $F_2$, respectively as defined in \eqref{eq:selfSimilarForm}.
    Suppose $F_1$ and $F_2$ have the same mass $M$. 
    Since $F_1$ and $F_2$ are both positive and
    \begin{displaymath}
        \int_{\R^d} F_1(z)\dz = \int_{\R^d} F_2(z)\dz = M,
    \end{displaymath}
    one either has that $F_1 = F_2$ a.e.~on $\R^d$, or
    \begin{displaymath}
        \int_{\R^d} (F_1(z) - F_2(z))^{+}\dz > 0
    \end{displaymath}
    and
    \begin{displaymath}
        \int_{\R^d} (F_2(z) - F_1(z))^{+}\dz > 0.
    \end{displaymath}
    Then \eqref{eq:selfSimilarForm} implies that for $t > 0$, the difference $\Gamma_1(t)-\Gamma_2(t)$ of the corresponding Barenblatt solutions $\Gamma_1$ and $\Gamma_2$ is also sign-changing.
    In particular, $\Gamma_1$ and $\Gamma_2$ satisfy the $L^1$ dissipation inequality for differences, Theorem~\ref{thm:L1dissipation}, with strict inequality.
    Hence using the form \eqref{eq:selfSimilarForm} for $\Gamma_1$ and $\Gamma_2$,
    \begin{displaymath}
        t_2^{-d\beta}\norm{[F_1(xt_2^{-\beta})-F_2(xt_2^{-\beta})]^{+}}_1 < t_1^{-d\beta}\norm{[F_1(xt_1^{-\beta})-F_2(xt_1^{-\beta})]^{+}}_1
    \end{displaymath}
    for $0 < t_1 < t_2$.
    Now a change of variable implies that
    \begin{displaymath}
        \norm{[F_1-F_2]^{+}}_1 < \norm{[F_1-F_2]^{+}}_1
    \end{displaymath}
    giving a contradiction.
\end{proof}

\end{document}